\documentclass[reqno]{amsart}

\theoremstyle{plain}
\newtheorem{thm}{Theorem}[section]
\newtheorem{lem}[thm]{Lemma}

\theoremstyle{definition}
\newtheorem{df}[thm]{Definition}
\newtheorem{exa}[thm]{Example}

\theoremstyle{remark}
\newtheorem{rem}[thm]{Remark}
\newenvironment{acknowledgements}{\section*{Acknowledgements}}{}

\usepackage{amsmath,amssymb}
\usepackage{graphicx}       %% for pdflatex and lualatex
\graphicspath{{./graph/}}
\usepackage{color}
\usepackage{here}%for fix graphics files
\usepackage{placeins} % for graphics float barrier
\usepackage{listings}
\usepackage{xcolor} % for setting colors
\usepackage{bm}
\usepackage{multirow}
\usepackage{mathtools}
\usepackage{bbold}
\usepackage[ruled,vlined,linesnumbered]{algorithm2e}
\SetKwComment{Comment}{/* }{ */}
\usepackage{algorithmic}
\usepackage{enumerate}
\usepackage{nicematrix} % for matrix with comments\\
\usepackage{enumitem}

\usepackage{xcolor}

\DeclareMathOperator{\Kernel}{Ker}
\DeclareMathOperator{\card}{card}
\DeclareMathOperator{\radius}{rad}
\DeclareMathOperator{\weight}{wt}
\DeclareMathOperator{\support}{supp}
\DeclareMathOperator{\reduced}{red}
\DeclareMathOperator{\rank}{rank}

\usepackage[most]{tcolorbox}
\newtcolorbox{highlighted}{colback=yellow,coltext=black,breakable}

\SetAlgoLined
\SetKwProg{MyStruct}{Struct}{ contains}{end}
\SetKwProg{MyFunction}{Function}{}{end}

\usepackage{xcolor}
\newcommand{\rev}[1]{#1}

\usepackage[hidelinks]{hyperref}
\numberwithin{equation}{section}
\numberwithin{figure}{section}
\numberwithin{table}{section}
\usepackage{mathptmx}      % use Times fonts if available on your TeX system
\usepackage{graphicx} % used for the example graphics below 
\usepackage{latexsym, amsfonts,amssymb,amsmath} %some typical packages 

\usepackage{etoolbox}

\makeatletter
\newcommand*{\JELname}{JEL Classification}
\newcommand*{\JEL}[1]{\gdef\@JEL{#1}}
\let\@JEL\@empty
\def\@setJEL{{\itshape\JELname.}\enspace\@JEL\@addpunct.}
\patchcmd{\@adminfootnotes}
  {\ifx\@empty\@keywords}
  {\ifx\@empty\@JEL\else\@footnotetext{\@setJEL}\fi\ifx\@empty\@keywords}
  {}{\g@addto@macro\@adminfootnotes
      {\ifx\@empty\@JEL\else\@footnotetext{\@setJEL}\fi}}
\makeatother

\begin{document}

\title[High-order recombination for weak approximation of SDE{\lowercase{s}}]{A high-order recombination algorithm for weak approximation of stochastic differential equations}

\author[Ninomiya]{Syoiti Ninomiya}
\address{Department of Mathematics, Institute of Science Tokyo}
\email{ninomiya@math.sci.isct.ac.jp}

\author[Shinozaki]{Yuji Shinozaki}
\address{School of Business Administration, Hitotsubashi University}
\email{yshinozaki@hub.hit-u.ac.jp}

\thanks{This work was supported by JSPS KAKENHI Grant Numbers JP21K03365(Scientific Research (C)) and JP25K17302(Early-Career Scientists).}

\JEL{C02, C63, G13}

\begin{abstract}
  This paper presents an algorithm for applying the high-order
  recombination method, originally introduced by Litterer and Lyons in
  \textit{``High-order recombination and an application to
  cubature on Wiener space''} (2012),
  to practical problems in mathematical finance.
  A refined error analysis is developed, yielding a sharper condition for space partitioning.
  Building on this condition, a computationally feasible recursive partitioning algorithm is constructed.
  Numerical examples are provided, demonstrating that the
  proposed algorithm effectively mitigates the explosive growth in the cardinality
  of the support as the required approximation order increases.
\end{abstract}

\maketitle

\newcommand{\cubthree}{\mathord{\mathsf{Cub3}}}
\newcommand{\NV}{\mathord{\mathsf{NV}}}
\newcommand{\Alg}{\mathord{\mathsf{Alg}}}
\section{Introduction}
\noindent
In Litterer~\cite{LittererDthesis}
as well as
Litterer and Lyons~\cite{litterer2012high}, it is demonstrated that the
probability measure arising in weak approximations of a stochastic
differential equation (SDE) can avoid exponential growth in the
cardinality of its support with respect to the approximation order,
by means of a method that transforms the measure of a given finite
probability space into one with smaller support, while preserving
moments up to a prescribed degree.
This method, referred to as the high-order recombination method, is the main subject of the present paper.
The aim here is to present a theoretically guaranteed
and practically feasible algorithm for the
high-order recombination method,
based on a refined error analysis and an improved condition
for space partitioning, and to demonstrate its applicability
to real problems in mathematical finance while maintaining the
theoretically predicted performance.

\subsection{Weak approximation problem of SDEs}
Let $(\Omega,\mathcal{F},P)$ be a probability space,
and let
$(B^1_t, \dots, B^d_t)$
denote the standard Brownian motion in $\mathbb{R}^d$.
For notational convenience,
we define
$B^0_t=t$.
We denote by
$C^\infty_b(\mathbb{R}^N;\mathbb{R}^N)$
the set of $\mathbb{R}^N$-valued smooth functions on $\mathbb{R}^N$
whose derivatives of all orders are bounded.
Let $I_N$ denote the identity
map on $\mathbb{R}^N$.
Let
$X_{t}^{x}$ be an $\mathbb{R}^N$-valued
diffusion process defined as the solution to the Stratonovich-form SDE:
\begin{equation}\label{eq:SDE}
  X_{t}^{x}=x+\sum_{i=0}^d\int_0^t V_iI_N(X_{s}^{x})\circ \mathrm{d}B^i_{s},
\end{equation}
where
$x\in\mathbb{R}^N$ and
$V_0,\dots,V_d$ are tangent vector fields on $\mathbb{R}^N$,
each having coefficients in $C^\infty_b(\mathbb{R}^N;\mathbb{R}^N)$. 
To simplify notation, we omit the dependence on the initial value $x$ whenever it is clear from the context. 
For each $i\in\{0,\dots,d\}$,
$V_i I_N =
\begin{pmatrix}
  (V_iI_N)^1 & \dots & (V_iI_N)^N
\end{pmatrix}
\in C_b^\infty(\mathbb{R}^N; \mathbb{R}^N)$,
and
the identity
\begin{equation*}
  V_ig(y)=
    \sum_{j=1}^N (V_iI_N)^j(y)\frac{\partial g}{\partial x_j}(y)
\end{equation*}
holds for all $g\in C^\infty(\mathbb{R}^N)$ and $y\in\mathbb{R}^N$.
Furthermore, the SDE~\eqref{eq:SDE} can also be written in It\^o form:
\begin{equation*}
  X_t=x+\int_0^t \tilde{V_0}I_N(X_{s})\,\mathrm{d}s
  +\sum_{i=1}^d\int_0^t V_i I_N(X_{s})\, \mathrm{d}B^i_{s},
\end{equation*}
where
\begin{equation*}
  (\tilde{V}_0I_N)^k=(V_0I_N)^k
  +\frac{1}{2}\sum_{i=1}^dV_i(V_iI_N)^k
\end{equation*}
for all $k\in\{1,\dots,N\}$.
This identity is known as the It\^o--Stratonovich transformation.
\par
In this paper, we consider the weak approximation problem for SDEs;
that is, the numerical evaluation of $E[f(X_T)]$,
where $X_{t}$ denotes the value of an $N$-dimensional diffusion
process at time $t$ and $f$ is a $\mathbb{R}$-valued function defined
on $\mathbb{R}^{N}$ with appropriate regularity.
This type of computation is crucial in many areas of applied science and
engineering---for example, in the pricing and hedging of
financial derivatives.
There are two main approaches to solving this problem numerically:
one based on partial differential equations and the other on simulation.
We focus on the latter, in particular on a so-called
tree-based simulation method, which is described in a later subsection.
\subsection{High-order cubature on Wiener space method}
The high-order cubature on Wiener space methods of Lyons and Victoir~\cite{lyons2004cubature}, 
together with the formulation provided by Ninomiya and Shinozaki~\cite{NinomiyaShinozaki2019}, 
are adopted.
In these methods, one approximates the
infinite-dimensional
Wiener measure by a finite-dimensional measure on a given partition of
the time interval $[0,T]$.
We denote this partition by $\Delta=\{0=t_0<\dots<t_{\sharp\Delta}=T\}$
and set $s_j=t_{j}-t_{j-1}$. 
Here, $\sharp\Delta$ denotes the number of subintervals contained in
$\Delta$.
\subsubsection{$\cubthree$: a cubature of degree $3$}\label{ss:cub3}
As a simple illustrative example of a cubature measure on Wiener space,
we outline a degree~3 cubature method.
This method is called $\cubthree$ and is described in  Lyons and Victoir
~\cite{lyons2004cubature} and
Ninomiya and Victoir~\cite{ninomiya2008weak}.

\begin{align*}
  X^{(\cubthree,\Delta)}_{t_0}(\eta)
    &= x, \notag\\
  X^{(\cubthree,\Delta)}_{\lambda t_{j+1}+(1-\lambda)t_j}(\eta)
    &= \notag\\
  \exp &\bigg(
    \lambda\bigg(
      s_{j+1}V_0+\sum_{i=1}^d\sqrt{s_{j+1}}\eta_{i,j+1}V_i
    \bigg)
  \bigg)X^{(\cubthree,\Delta)}_{t_j}(\eta) 
\end{align*}
where $0\leq\lambda\leq 1$,
$\exp(V)x$ denotes the solution at time $1$
of the ordinary differential equation (ODE)
\begin{equation}\label{eq:ODE}
  \frac{\mathrm{d}z_t}{\mathrm{d}t} = VI_N (z_t),\;\; z_0=x,
\end{equation}
for a tangent vector field $V$ on $\mathbb{R}^N$, and
\begin{equation*}
\eta=(\eta_{i,j})_{\substack{i\in\{1,\dots, d\}\\j\in\{1,\dots,\sharp\Delta\}}}
\end{equation*}
denotes a $d\sharp\Delta$-dimensional vector whose components are
independent and identically distributed random variables satisfying $P[\eta_{i,j}=\pm 1]=1/2$.
We also remark that for any $s\in\mathbb{R}$, the expression
$\exp(sV)(x)$ is well-defined by virtue of the chain rule and
solves the ODE~\eqref{eq:ODE}.
\par
We then obtain a cubature measure on Wiener space of degree~3,
denoted by $\mu^{(\cubthree,\Delta)}$,
whose support is given by the collection of paths
$\{
\omega_\eta\,:\,\eta\in\{-1, 1\}^{d\sharp\Delta}
\}$,
where each path is defined as
\begin{equation*}
  \omega_\eta =
  \big(X_s^{(\cubthree,\Delta)}(\eta)\big)_{s\in [0,T]}.
\end{equation*}
Every path $\omega_\eta$ carries an equal weight of $2^{-d\sharp\Delta}$.
Therefore, the computational cost of evaluating
\begin{align*}
  \mu_T^{(\cubthree,\Delta)}(f)
    &= E[f(X^{(\cubthree,\Delta)}_T)] \\
    &= \sum_{\eta\in\{-1,1\}^{d\sharp\Delta}}
       2^{-d\sharp\Delta}
       f(\omega_\eta|_{s=T}).
\end{align*}
grows exponentially with $\sharp \Delta$.
\subsubsection{Tree-based simulation}\label{tbsec}
Discretisation methods similar to the $\cubthree$ method are
referred to as tree-based simulations.
They proceed by drawing a direct product of discrete random variables
$\eta$.
This terminology stems from the structural interpretation that the support
$\{\omega_\eta\,\,:\,\eta\in\{-1,1\}^{d\sharp\Delta}\}$
of the measure can be viewed as a tree, with the root
$\omega_{\eta}\vert_{s=0}$ and the leaves
$
\{
\omega_{\eta}\vert_{s=T}\,:\,
\eta\in\{-1,1\}^{d\sharp\Delta}\}
$.
\par
We now introduce notation for tree-based simulations. 
Let 
$(X^{(\Alg, \Delta)}_{t_i})_{i\in\{0,\dots,\sharp\Delta\}}$
denote the discretised process obtained by a discretisation method
$\Alg$,
based on discrete random variables with respect to the partition
$\Delta$ of the interval $[0,T]$.
We denote by $\mu_{t_i}^{(\Alg, \Delta)}$
the probability measure on $\mathbb{R}^N$
induced by $X^{(\Alg, \Delta)}_{t_i}$ and write it as  
\begin{equation*}
  \mu_{t_i}^{(\Alg, \Delta)} = \sum_{j=1}^{M_i} w_j^{(i)} \delta_{p_j^{(i)}},
\end{equation*}
where  
$((p_j^{(i)}, w_j^{(i)}))_{j\in\{1,\dots,M_i\}}$  
is the set of pairs consisting of points
$p_j^{(i)} \in \mathbb{R}^N$
and the corresponding weights $w_j^{(i)}$.  
The symbol $\delta_x$ denotes the Dirac measure at $x\in\mathbb{R}^N$.
%%%%%%
We also denote the transition probability by
\begin{equation*}
  w^{(i)}(k,j) = P [X^{(\Alg, \Delta)}_{t_{i+1}}
  =p_j^{(i+1)} \vert\, X^{(\Alg, \Delta)}_{t_i}
  =p_k^{(i)} ].
\end{equation*}
We then define the one-step forward operator on the measure space as
\begin{equation*}
  \hat{Q}^{(\Alg, \Delta)}_{(s_{i+1})}\bigg(
  \sum_{j=1}^{M_i}\alpha_j\delta_{p_{j}^{(i)}} \bigg)
  =\sum_{k=1}^{M_i}\sum_{j=1}^{M_{i+1}}\alpha_k w^{(i)}(k,j)\delta_{p_{j}^{(i+1)}},
\end{equation*}
where $\alpha_i\in\mathbb{R}$.
We also define
\begin{equation*}
  w_{j}^{(i+1)}=\begin{cases}
  P[X^{(\Alg,\Delta)}_{t_1}=p^{(1)}_j] & \text{if $i=0$}, \\
    \sum_{k=1}^{M_i} w^{(i)}(k, j) w_{k}^{(i)}
  &
  \text{if $1\leq i\leq \sharp\Delta-1$.}
  \end{cases}
\end{equation*}
Finally, we obtain the measure
$\mu_T^{(\Alg,\Delta)}=
\hat{Q}^{(\Alg, \Delta)}_{(s_{\sharp \Delta})}\circ
\cdots \circ \hat{Q}^{(\Alg, \Delta)}_{(s_{1})}
(\delta_{x})$.
The process of approximating
$E[f(X_T)]$ by the expectation
\begin{equation}\label{tbcalc}
  \mu_T^{(\Alg,\Delta)}(f)
  =\sum_{j=1}^{M_{\sharp\Delta}}w_j^{(\sharp\Delta)}
  f(p_j^{(\sharp\Delta)}),
\end{equation}
is referred to as a weak approximation by tree-based simulation. 
The $\cubthree$ method $\mu_T^{(\cubthree,\Delta)}$, introduced in
subsection~\ref{ss:cub3}, is one of the simplest examples of
such tree-based simulation methods. The N--V method
(Ninomiya and Victoir~\cite{ninomiya2008weak}) $\mu_T^{(\NV,\Delta)}$,
which serves as
the subject of the numerical experiments in this paper, is introduced
in subsection~\ref{nvsec}.
\subsubsection{Order of discretisation}\label{ss:order}
The order of a discretisation method is defined as follows.
\begin{df}\label{df:pth-order}
  $\Alg$ is said to be a discretisation of order $p$
  or a $p$th order discretisation
  if there exists a constant $C_p>0$ such that for all $n\in\mathbb{N}$
  there exists a partition $\Delta$ of $[0,T]$ with $\sharp\Delta=n$
  satisfying
%\begin{equation*}
%  \Big\lvert
%  E\big[f(X_{T})\big]
%  - E\Big[f\big(X^{(\Alg, \Delta)}_{T}\big)\Big]
%  \Big\rvert < C_p (\sharp \Delta)^{-p}.
%\end{equation*}
\begin{equation*}
  \lvert
  E [f(X_{T}) ]
  - E [f (X^{(\Alg, \Delta)}_{T})]
  \rvert < C_p (\sharp \Delta)^{-p}.
\end{equation*}

\end{df}
In what follows, discretisation methods of order $2$ or higher are
referred to as high-order discretisation methods.
Several such high-order discretisation methods based on cubature on
Wiener space have been proposed;
see for example Ninomiya and Victoir~\cite{ninomiya2008weak},
Ninomiya and Ninomiya~\cite{ninomiya2009new}
and Shinozaki~\cite{Shinozaki2017}.
The particular method used in this paper is introduced in subsection~\ref{nvsec}.
\par
In addition, research on the construction of high-order cubature measures has
explored various avenues, including simple algebraic
construction by Gyurk{\'{o}} and Lyons~\cite{gyurko2011efficient},
randomised construction by Hayakawa and Tanaka~\cite{hayakawa2022monte}
and constructions based on Hopf algebras and unshuffle
expansions by
Ferrucci, Herschell, Litterer and Lyons~\cite{ferrucci2024high}.
%%%%
For ODEs, discretisation methods of fixed order can yield methods of
higher order when they are linearly combined using weights that may
include negative values.
This approach is commonly referred to as extrapolation. 
In the context of weak approximation
for SDEs, the extrapolation of the Euler--Maruyama method is studied
in Bally and Talay~\cite{bally1996law}.
Furthermore, as demonstrated
in Fujiwara~\cite{fujiwara}
as well as
in Oshima, Teichmann and Velu{\v{s}}{\v{c}}ek~\cite{oshima2012new},
extrapolation methods can be
combined with the method of Ninomiya and Victoir~\cite{ninomiya2008weak}
to construct
discretisations of arbitrary order.

%%%%%%%%%%%%%%%%%%%%%%%%%%%%%%%%%%%%%%%%%%%%%%%%%%
\subsection{High-order recombination}
When implementing high-order cubature measures on Wiener space,
one often encounters the problem of exponential growth
in the cardinality of the $\support(\mu_T^{(\Alg,\Delta)})$
as the number of time partitions $\sharp \Delta$ increases.
%%%%%%%%%%%%%%%%%%%%%
For example, as seen in subsection~\ref{ss:cub3},
the cardinality of the support of
$\mu_T^{(\cubthree,\Delta)}$ is $2^{d \sharp \Delta}$.
That of $\mu_T^{(\NV,\Delta)}$ is \rev{$(2\times 3^d)^{\sharp\Delta}$},
where $\NV$ refers the discretisation method of order $2$
defined in subsection~\ref{nvsec}.
%%%%%%%%%%%%%%%%%%%%%%%%%%%%%%%%%%%%%%%%%%%%%%%%%%
In what follows, we refer to this issue as the support explosion problem. 

%% \begin{table}[H]
  \begin{table}\normalsize
    \centering
    \begin{tabular}{|l||c|c|}
      \hline
      \begin{tabular}{c}
      Algorithm\\ (Discretisation method)
      \end{tabular}
      & \begin{tabular}{c}
          Number of $(\hat{\eta})$ \\ in $\hat{Q}_{(s_i)}$
        \end{tabular}
      & \begin{tabular}{c}
        Cardinality of the support \\of
        the measure of \\
        $\hat{Q}_{(s_n)} \circ \cdots \circ \hat{Q}_{(s_1)}$
        \end{tabular}
      \\
      \hline\hline
      $\cubthree$ & $d$ & $2^{dn}$ \\
      $\NV$ & $d$ & $(2 \times 3^{d})^{n}$ \\
      N--N \cite{ninomiya2009new} & $2d$ & $3^{2dn}$ \\
      $Q_{(s)}^{(7,2)}$ ($d=2$) \cite{Shinozaki2017} & $7$ & $5^{7n}$ \\
      \hline
    \end{tabular}
    \caption{Number of random variables and support sizes for
      $\hat{X}^{(\text{method})}_{s,x}$ with $\sharp\Delta = n$}
    \label{tab:NUMS}
  \end{table}
\par
At present the most widely used approach to overcome this problem is
partial sampling such as Monte Carlo and quasi-Monte Carlo methods.
Among these the technique known as TBBA is particularly effective
and well suited to the tree-based simulation
methods Crisan and Lyons~\cite{dan2002minimal},
Ninomiya~\cite{ninomiya2003partial}
and
Ninomiya~\cite{ninomiya2010application}.
Furthermore in Crisan and Manolarakis~\cite{Crisan2012}, 
as well as in their 2014 study~\cite{Crisan2014}
and in Ninomiya and Shinozaki~\cite{NinomiyaShinozaki2019},
tree-based simulation with TBBA is successfully applied
to the numerical calculation of forward-backward SDEs.
\par
High-order recombination is proposed in
Litterer and Lyons~\cite{litterer2012high} as an
alternative approach to resolving the support explosion problem.
There the authors show that the cardinality of the support can be kept polynomial in the required approximation accuracy.
Several applications of high-order recombination are presented
in Lee and Lyons~\cite{lee2016adaptive}.
In Tchernychova~\cite{tchernychova2016caratheodory},
a number of algorithms for obtaining reduced measures are discussed.
The notion of reduced measure plays a significant role
in the construction of recombination measures, and will be introduced
in the subsequent subsection~\ref{def:redmeas}.
\par
This subsection presents an overview of the procedure for constructing
high-order recombination measures.
These measures form the main subject of this paper.
\subsubsection{Reduced measure}\label{def:redmeas}
We first introduce the concept of a reduced measure.
It is a measure with reduced support that preserves relations defined
by a set of test functions.
%%%%
Let $\mu$ be a discrete probability measure on a space
$\tilde{\Omega}$,
and let \rev{$G=\{g_1,\ldots,g_{\card(G)}\}$} be a finite set of test
functions.
Throughout this paper, \(\card(S)\) denotes
the cardinality of a set \(S\).
Each \(g_i:\tilde{\Omega}\to\mathbb{R}\) is assumed to be measurable.
%%%%
\begin{df}\label{redmeas}
  A discrete probability measure $\tilde{\mu}$ is said to be a reduced measure from $\mu$ with respect to $G$ if it satisfies the following three conditions:
  \begin{enumerate}[label=\textnormal{(\roman*)}]
  \item $\support (\tilde{\mu} ) \subset
    \support (\mu )$,
  \item $\int g(x)\,\tilde{\mu}(\mathrm{d}x) = \int g(x)\,\mu(\mathrm{d}x)
    \quad\text{for all}\quad g \in G$,
  \item $\card(\support (\tilde{\mu}))
    \leq \card(G)$.
  \end{enumerate}
  The set of all reduced measures from \(\mu\) with respect to \(G\)
  is denoted by \(\reduced_G(\mu)\).
\end{df}
\noindent
In Litterer and Lyons~\cite{litterer2012high},
the existence of a reduced measure is
shown for a finite Borel measure $\mu$ on a Polish space
$\tilde{\Omega}$.
\par
For a discretisation \rev{algorithm defined by an $m$-moment similar family
(see Definition~\ref{MSF}), that is, of degree $m$, and with
$\tilde{\Omega} = \mathbb{R}^N$, we take as $G$ the set of all monomials
of degree at most $m$.}
Then $\card(G)=\binom{N+m}{m}$, and $G$ spans the space of polynomials
of degree at most $m$.
We number the monomials so that $g_1(x) = 1$, and hence
$g_2(0) = g_3(0) = \dots = g_{\card(G)}(0) = 0$.
Note that condition {\rm (ii)} is determined by $g_1(x)=1$.
It ensures that the reduced measure $\tilde{\mu}$ is also
a probability measure.
\subsubsection{High-order recombination measures and their construction procedures}\label{recsubsec}
The high-order recombination measures
$(\mu^{(\Alg,\Delta,\rm{Rec})}_{t_i})_{i\in\{0,\dots,\sharp\Delta\}}$
in tree-based simulation with respect to
$(\Alg, \Delta)$ are constructed recursively as follows:
\begin{align}
  \mu^{(\Alg,\Delta,\rm{Rec})}_{t_0} &= \delta_{x_0}, \notag\\
  \nu^{(i)} &= \hat{Q}_{(s_i)}^{(\Alg,\Delta)}(\mu_{t_{i-1}}^{(\Alg,\Delta,\rm{Rec})}), \notag\\
  \mu^{(\Alg,\Delta,\rm{Rec})}_{t_i} &= \sum_{k=1}^{l_i} \nu^{(i)}(U^{(i)}_k) \tilde{\nu}^{(i)}_k. \label{eq:high-order-recombination}
\end{align}
%
% Revised in November 2026: no recombination at the last step, as in
% Litterer and Lyons.
\rev{The second relation is applied for $i\in\{1,\dots,\sharp\Delta\}$ and
the third for $i\in\{1,\dots,\sharp\Delta-1\}$.
At the last step no recombination is performed, that is,
$\mu^{(\Alg,\Delta,\rm{Rec})}_{t_{\sharp\Delta}}=\nu^{(\sharp\Delta)}$,
as in Litterer and Lyons~\cite{litterer2012high}.
Recombination at the last step would not reduce the computational cost,
since $\nu^{(\sharp\Delta)}$ is only used to evaluate the expectation of
the payoff.
For $i\in\{1,\dots,\sharp\Delta-1\}$,
$(U^{(i)}_k)_{k\in\{1,\dots,l_i\}}$} is a set of subsets of
$\support (\nu^{(i)})$ satisfying
\begin{equation*}
  \bigcup_{k=1}^{l_i} U^{(i)}_k = \support
  (\nu^{(i)})
  \quad \text{and} \quad U^{(i)}_j \cap U^{(i)}_k = \emptyset
  \quad \text{for all } j \ne k.
\end{equation*}
Moreover,
\begin{equation*}
  \tilde{\nu}^{(i)}_k \in
  \reduced_G(\nu^{(i)}\vert_{U^{(i)}_k})
  \quad \text{for } 1 \leq k \leq l_i,
\end{equation*}
\begin{equation}\label{eq:patch-measure}
  \nu^{(i)}\vert_{U_{k}^{(i)}}
  = \frac{1}{\nu^{(i)}(U_k^{(i)})}
  \sum_{p_{j}^{(i)} \in U_{k}^{(i)}}
  \hat{w}_{j}^{(i)} \delta_{p_{j}^{(i)}},
\end{equation}
which is the probability measure on
\(U^{(i)}_k\) induced by
\(\nu^{(i)} = \sum_{j=1}^{M_i} \hat{w}_j^{(i)} \delta_{p^{(i)}_j}\).
%%
%%%%%%%%%%%%%%%%%%%%%%%%%%%%%%%%%%%%%%%%%%%%%%%%%%
Henceforth a family $(U_k)_{k\in\{1,\dots,l_i\}}$ satisfying
\begin{equation*}
  \bigcup_{k=1}^{l_i} U_k  = \support(\nu^{(i)})\quad\text{and}\quad
  U_j\cap U_k = \emptyset\quad\text{for all $j\neq k$}
\end{equation*}
is referred to as a patch division of $\support(\nu^{(i)})$.
This procedure is interwoven with the tree-based simulation and
is carried out simultaneously at each time step except the last one.
The total weight assigned to each patch remains invariant
throughout the recombination procedure.
We therefore refer to this as the weight of the patch
\(U_k^{(i)}\) and denote it by \(\weight(U_k^{(i)})\).
The following identity summarises the quantities introduced above:
\begin{equation*}
  \weight(U^{(i)}_k)=\nu^{(i)}(U^{(i)}_k)
  =\mu_{t_i}^{(\Alg,\Delta,\rm{Rec})}(U^{(i)}_k).
\end{equation*}
A simplified illustration of these procedures can be found in
Figure~\ref{illust} and  a toy example of the
reduced measure is provided in Appendix~\ref{appendix2}.
\begin{figure}\normalsize
  \begin{center}
    \caption{Illustration of recombination method} \label{illust}
    \includegraphics [width=8cm]{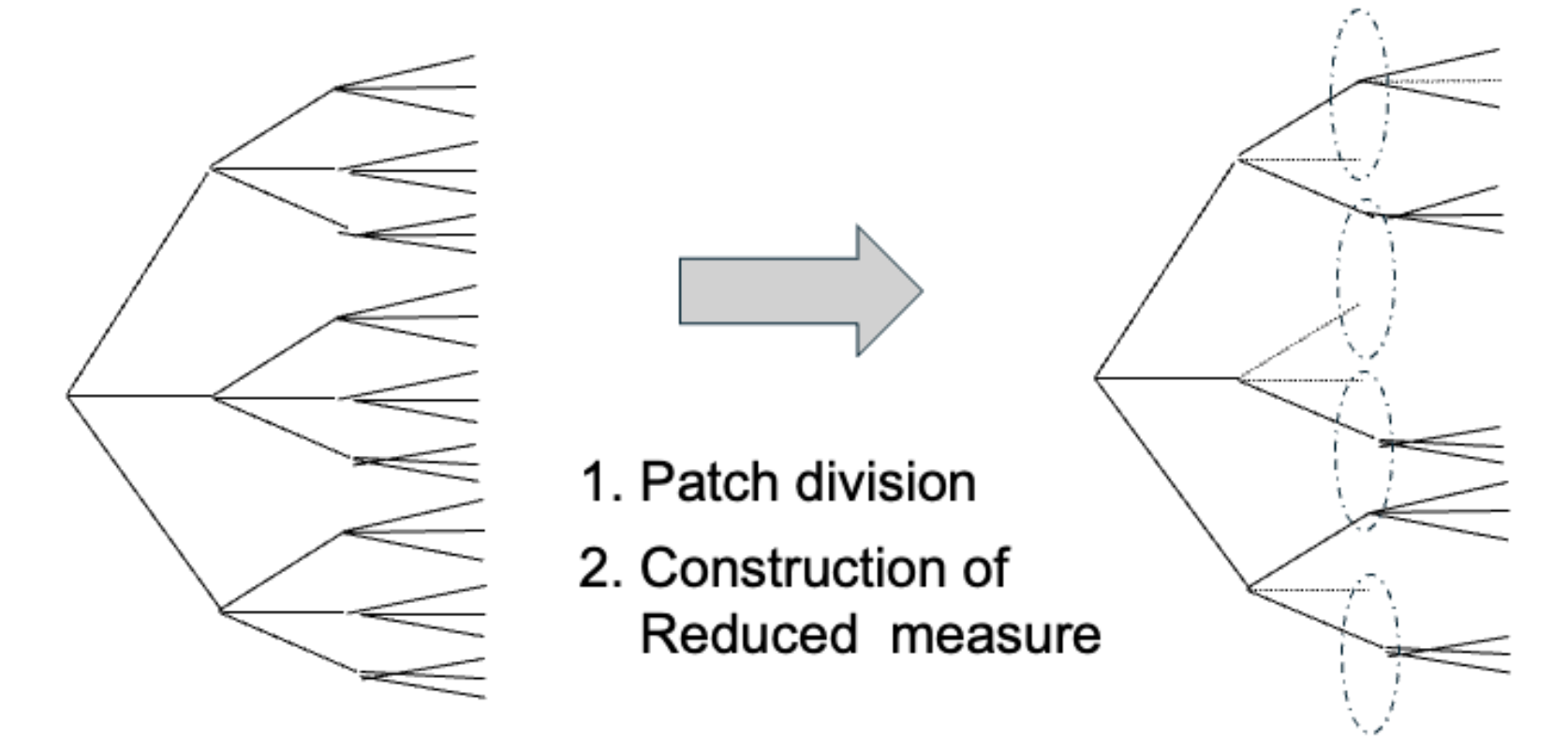}    
  \end{center}
\end{figure}
%%%%
It is shown in Litterer and Lyons~\cite{litterer2012high}
that high-order recombination
can be performed in such a way that the order of discretisation
is preserved by choosing appropriate patch radii
\begin{equation*}
  \radius(U_j)
  = \max\{
  \lvert p_k - p_{k^\prime}\rvert
  \,:\, p_k, p_{k^\prime} \in U_j
  \},
\end{equation*}
where $\lvert\cdot\rvert$ denotes the Euclidean norm on
$\mathbb{R}^N$.
Numerical examples of the high-order recombination method applied to
practical problems appear only in Uchino~\cite{uchino2020} and
in Ninomiya and Shinozaki~\cite{rec2021}.
It is shown
numerically in those papers that, even when ad hoc patch division is
employed---which has not been proved to avoid the support explosion
problem---the recombination method can still be effective for certain
problems.
\par
This paper provides a detailed explanation of the algorithms
corresponding to the aforementioned procedure.
In particular:
\begin{itemize}
\item The explicit construction of $\hat{Q}_{(s_i)}^{(\NV,\Delta)}$
  is presented in subsection~\ref{nvsec}.
\item An algorithm for the patch division of $\support(\nu^{(i)})$
  into $(U^{(i)}_k)_{k\in\{1,\dots,l_i\}}$
  is proposed in subsection~\ref{patchdivsec}.
\item The algorithm for finding
  $\tilde{\nu}^{(i)}_k \in
  \reduced_G(\nu^{(i)}\vert_{U^{(i)}_k})$
  ---which is identical to the one used in
  Ninomiya and Shinozaki~\cite{rec2021}---
  is provided in
  Appendix~\ref{appendix1}.
\end{itemize}
\subsection{Contribution of this paper}
This paper presents a recursive patch division algorithm for
high-order recombination that addresses the support explosion problem
and thereby makes the method practically feasible for applications in
mathematical finance.
The contributions of this paper are summarised as follows:
\begin{itemize}[label=\textbullet]
\item An improved version of the error estimation
  in Litterer and Lyons~\cite{litterer2012high}
  is developed which incorporates the patch
  weights $\weight (U^{(i)}_k)$.
\item A recursive patch division algorithm is proposed which divides
  the support set
  while taking into account both the radii and the weights of
  the patches.
  The algorithm includes a control parameter that can be varied to change
  the degree of recombination.
  Its role is discussed in Example~\ref{pel}, and it is used as a tuning
  parameter in the numerical investigation of the recombination error
  in Section~\ref{subsec:RecErr}.
  Its recursive structure
  allows efficient and adaptive patch division based on the refined
  error evaluation.
\item
  A numerical example is presented to demonstrate the application of
  the proposed algorithm to realistic option pricing problems.
  The example shows that the algorithm successfully avoids the support
  explosion problem and illustrates how adjusting the control parameter
  affects the degree of recombination.
\end{itemize}
The remainder of this paper is organised as follows.
Section~\ref{sec:preliminary}
introduces the necessary notation and preliminary results on
high-order discretisation schemes and high-order recombination.
Section~\ref{sec:refined} presents the refined error evaluation
and the novel recursive patch division algorithm.
These constitute the main contributions of this paper.
Section~\ref{sec:numerical} provides numerical results
that demonstrate the performance of the proposed algorithm.
\section{Preliminary}\label{sec:preliminary}
\noindent
This section introduces the notation and prior work on high-order
discretisation and the recombination algorithm.
\subsection{Free Lie algebras}
Let $A=\{0, 1, \dots , d\}$ and define $A^+=\bigcup_{k=1}^\infty A^k$
as the set of all words consisting of elements of $A$, and define
$A^\ast=\{\emptyset\}\cup A^+$.
$A^\ast$ becomes a semigroup with identity $\emptyset$ when endowed
with the concatenating product `$\cdot$'. 
For arbitrary elements $(a_1,\dots,a_m)$ and
$(b_1,\dots,b_n)$ in $A^\ast$,
\begin{equation*}
  (a_1,\dots,a_m)\cdot(b_1,\dots,b_n)
  = (a_1,\dots,a_m,b_1,\dots,b_n)
\end{equation*}
and
\begin{equation*}
  \emptyset\cdot(a_1,\dots,a_m)
  = (a_1,\dots,a_m)\cdot\emptyset
  = (a_1,\dots,a_m).
\end{equation*}
To simplify the notation, we identify $(a_1, a_2,\ldots,a_m)\in A^\ast$
with
$a_1\cdot a_2 \cdots  a_m$.
If there is no risk of confusion, we also omit the symbol `$\cdot$'.
We prefer to write $a_1a_2\ldots a_m$ instead of
$(a_1, a_2, \ldots, a_m)$ or $a_1\cdot a_2\cdots \cdot a_m$.
For $\alpha=a_{1}a_{2}\dots a_{k} \in A^\ast$, we define the length of
$\alpha$ by $\lvert\alpha\rvert=k$,
and set $\lvert\emptyset\rvert=0$.
We also define the order of
$\alpha$ as
\begin{equation*}
  \|\alpha\|=
  \lvert\alpha\rvert+\card(\{1\leq i\leq k\;:\;a_{i}=0\}).
\end{equation*}
Let $A^\ast_1=A^\ast\setminus \{\emptyset, 0\}$ and define
\begin{equation*}
  A^\ast(m)=\{\alpha \in A^\ast\;:\;\|\alpha\|\leq m\}\quad\text{and}
  \quad A^\ast_1(m)=\{\alpha \in A^\ast_1\;:\;\|\alpha\|\leq m\}
\end{equation*}
for $m\in \mathbb{N}$.
\par   
Let $\mathbb{R}\langle A \rangle$ be the free
$\mathbb{R}$-algebra with basis
$A^\ast$ and $\mathbb{R}\langle\langle A\rangle\rangle$ be its completion.
These are given by
\begin{equation*}
  \mathbb{R} \langle A \rangle =
  \bigg\{\sum_{k=1}^{l} r_{k}\alpha_{k}\;:\;
    l\in\mathbb{N},\,r_{k} \in \mathbb{R},\,\alpha_{k} \in A^\ast \bigg\}
\end{equation*}
and
\begin{equation*}
  \mathbb{R}\langle\langle A\rangle\rangle=
  \bigg\{\sum_{\alpha \in A^\ast} r_{\alpha}\alpha \;:\;
  r_{\alpha} \in \mathbb{R} \bigg\}.
\end{equation*}
The latter is the set of all $\mathbb{R}$-coefficient formal series
with basis
$A^\ast$.
For $m \in \mathbb{Z}_{\geq 0}$, let $j_{m}$ be a truncation operator
defined by
\begin{equation*}
  j_{m}\bigg(\sum_{\alpha \in A^\ast} r_{\alpha}\alpha\bigg)=
  \sum_{\|\alpha\| \leq m} r_{\alpha}\alpha.
\end{equation*}
The multiplication in $\mathbb{R}\langle\langle A \rangle\rangle$
is defined by extending the product in $\mathbb{R}\langle A \rangle$
continuously.
More precisely, for such
$ R, S \in \mathbb{R}\langle\langle A\rangle\rangle$ written as
\begin{equation*}
  R = \sum_{\alpha\in A^\ast} r_\alpha \alpha, \quad
  S = \sum_{\beta\in A^\ast} s_\beta \beta,
\end{equation*}
their product is given by
%%% \begin{equation*}
%%%   RS = \sum_{\gamma\in A^\ast}
%%%   \bigg(\sum_{\substack{\alpha,\beta\in A^\ast \\ \alpha\beta=\gamma}}
%%%   r_\alpha s_\beta\bigg)\gamma.
%%% \end{equation*}
%%% %%
\begin{equation*}
   RS = \sum_{\gamma\in A^\ast}
   \bigg(\sum_{(\alpha,\beta\in A^\ast)\wedge (\alpha\beta=\gamma)}
   r_\alpha s_\beta\bigg)\gamma.
\end{equation*}
This product makes $\mathbb{R}\langle\langle A \rangle\rangle$
a unital non-commutative algebra over $\mathbb{R}$,
in which $\mathbb{R}\langle A \rangle$ is a subalgebra.
The Lie bracket is defined by $[R, S] = RS-SR$ as usual.
\par
${\mathcal L}_{\mathbb{R}}(A)$ denotes the set of Lie polynomials over
$A$.
It is the smallest $\mathbb{R}$-submodule of
$\mathbb{R}\langle A\rangle$ that contains $A$ and is closed under the
Lie bracket $[,]$.
\par
Let $\Phi$ be the homomorphism from $\mathbb{R} \langle A\rangle$ to
the $\mathbb{R}$-algebra of smooth differential operators on
$\mathbb{R}^{N}$, defined by
\begin{align*}
  \Phi(\emptyset) &= \mathrm{Id},\notag\\
  \Phi(\alpha) &= V_{a_{1}} V_{a_{2}}\cdots V_{a_{k}}
  \quad \text{for $\alpha = a_{1} a_{2} \dots a_{k} \in A^+$.}
\end{align*}
$\mathrm{Id}$ denotes the identity operator.
%%%%%%
\par
We remark that $\Phi$ canonically induces a Lie algebra homomorphism from
$\mathcal{L}_\mathbb{R}(A)$ to the Lie algebra of vector fields generated by
$\{V_0, \dots, V_d\}$.
We denote $V_R = \Phi(R)$ for $R \in \mathbb{R}\langle A \rangle$,
and define $[\alpha] = [[[\dots[[a_1,a_2],a_3],\dots], a_{k-1}], a_k]$
for
$\alpha = a_1 a_2 \dots a_k \in A^*$.
That is, $[\alpha]$ is the left-normed Lie bracket of the letters
in $\alpha$.
%%%
% Pairing, 2-norm, dilation and flows (revised October 2026).
% The Taylor formula that stood here (with t R and O(t^{m+1})) was
% incorrect, since j_m truncates by the order of words; it is restated
% with the dilation below (revised November 2026).
For $R\in\mathbb{R}\langle\langle A\rangle\rangle$ and
$R'\in\mathbb{R}\langle A\rangle$ written as
$R=\sum_{\alpha}r_\alpha\alpha$ and $R'=\sum_{\alpha}r'_\alpha\alpha$,
we set $\langle R, R'\rangle=\sum_{\alpha}r_\alpha r'_\alpha$,
which is a finite sum.
In particular, $\langle R, \alpha\rangle=r_\alpha$ is the coefficient of
the word $\alpha$ in $R$, and we define
$\|R\|_{2}=\langle R, R\rangle^{1/2}$ for $R\in\mathbb{R}\langle A\rangle$.
\rev{For $s>0$, let
$\mathcal{D}_s\colon\mathbb{R}\langle\langle A\rangle\rangle
\to\mathbb{R}\langle\langle A\rangle\rangle$
be the linear map defined by
\begin{equation*}
  \mathcal{D}_s\bigg(\sum_{\alpha\in A^\ast}r_\alpha\alpha\bigg)
  =\sum_{\alpha\in A^\ast}s^{\|\alpha\|/2}r_\alpha\alpha,
\end{equation*}
which is denoted by $\Psi_s$ in Kusuoka~\cite{kusuoka2004approximation}.}
Since $\|\alpha\beta\|=\|\alpha\|+\|\beta\|$,
the map $\mathcal{D}_s$ is an algebra automorphism of
$\mathbb{R}\langle\langle A\rangle\rangle$.
It maps ${\mathcal L}_{\mathbb{R}}(A)$ onto itself, commutes with $j_m$,
and satisfies $\mathcal{D}_s(\exp(R))=\exp(\mathcal{D}_s(R))$
for $R$ without constant term.
Moreover, $\mathcal{D}_s\circ\mathcal{D}_{s'}=\mathcal{D}_{ss'}$
for $s,s'>0$.
We call the maps $\mathcal{D}_s$, $s>0$, dilations, since they form a
one-parameter group of algebra automorphisms and $\mathcal{D}_s$
multiplies the homogeneous part of order $n$ by $s^{n/2}$;
they reflect the scaling property of Brownian motion
(see Remark~\ref{rem:m-moment-similar}~(a)).
For $L\in{\mathcal L}_{\mathbb{R}}(A)$, $\Phi(L)$ is a vector field, and
$\exp(\Phi(L))\colon\mathbb{R}^N\to\mathbb{R}^N$ denotes its time-one flow,
that is, the solution map at time $1$ of the ODE~\eqref{eq:ODE}
with $V=\Phi(L)$.
For $R\in{\mathcal L}_{\mathbb{R}}(A)$,
$f\in C_b^{\infty}(\mathbb{R}^N;\mathbb{R})$ and $m\in\mathbb{N}$,
Taylor's formula along the flow (see Corollary~10 of
Kusuoka~\cite{kusuoka2004approximation}) gives the approximation
\begin{equation*}
  f\big(\exp(\Phi(\mathcal{D}_tR))(x)\big)
  =\Big(\Phi\big(j_m(\exp(\mathcal{D}_tR))\big)f\Big)(x)
  +O(t^{(m+1)/2}),\qquad t\downarrow 0,
\end{equation*}
since every word of order greater than $m$ in the expansion carries the
factor $t^{(m+1)/2}$ or a higher power of $t$.

%%%
\subsection{High-order discretisation scheme}
In this subsection, we introduce high-order discretisation schemes for SDEs
and related results.
\par
We use the following norms throughout the paper.
Let $C_b(\mathbb{R}^N)$ denote the space of bounded continuous functions
$f : \mathbb{R}^N \to \mathbb{R}$. For $f \in C_b(\mathbb{R}^N)$,
we define the uniform norm by
\[
  \|f\|_{\infty} = \sup_{x \in \mathbb{R}^N} \Big|f(x)\Big|.
\]
%%% The space C_0 and the L^p norms were used only in Theorem 2.3 of the
%%% published version, which was removed in October 2026.
%%% We also use the subspace
%%% \[
%%%   C_0(\mathbb{R}^N)
%%%   = \{
%%%     f \in C_b(\mathbb{R}^N)
%%%     \;:\;
%%%     \lim_{\|x\|\to\infty} f(x) = 0
%%%   \}.
%%% \]
%%% It is the space of continuous functions vanishing at infinity.
%%% For $1\le p <\infty$, we define the $L^p(\mathrm{d}x)$ norm for
%%% measurable functions $f : \mathbb{R}^N \to \mathbb{R}$ such that
%%% $|f|^p \in L^1(\mathbb{R}^N, \mathrm{d}x)$ by
%%% \[
%%%   \|f\|_{L^p(\mathrm{d}x)} = 
%%%   \bigg( \int_{\mathbb{R}^N} |f(x)|^p \, \mathrm{d}x \bigg)^{1/p}.
%%% \]
We denote by $C_{\mathrm{Lip}}(\mathbb{R}^N)$ the space of Lipschitz
continuous
functions on $\mathbb{R}^N$. For $f \in C_{\mathrm{Lip}}(\mathbb{R}^N)$,
the Lipschitz norm is denoted by
\[
  \|f\|_{\rm Lip}= 
  \sup_{x \neq y} \frac{|f(x)-f(y)|}{|x-y|}.
\]
%%%%%%
\begin{df}
  \label{dfpt}
  Let $X_{t}^{x}$ be the $\mathbb{R}^N$-valued diffusion process defined by~\eqref{eq:SDE}.
  For $f \in C_{\mathrm{Lip}}(\mathbb{R}^N; \mathbb{R})$ and $t \geq 0$,
  define
  \begin{equation*}
    (P^X_t f)(x) = E[f(X_{t}^{x})]
  \end{equation*}
  and obtain a semigroup of linear operators
  $(P^X_t)_{t \geq 0}$ on $C_{\mathrm{Lip}}(\mathbb{R}^N; \mathbb{R})$.
\end{df}
We begin by presenting an approximation theorem for the high-order
discretisation scheme based on a series of works by
Kusuoka~\cite{kusuoka2001approximation,kusuoka2003malliavin,kusuoka2004approximation}.
%%%%
\begin{df}\label{UFG}
  Let $\ell$ be an integer. 
  The vector fields $V_0, V_1, \dots, V_d$
  satisfy the uniformly finitely generated (UFG) condition if for all $\alpha \in A^\ast_1$
  there exists a family
  \begin{equation*}
    \{ \varphi_{\alpha, \beta}\,\,:\,\beta\in A^\ast_1(\ell)
    \}
    \subset C^{\infty}_b(\mathbb{R}^N; \mathbb{R})
  \end{equation*}
  such that
  \begin{equation*}
    V_{[\alpha]} = \sum_{\beta} \varphi_{\alpha, \beta} V_{[\beta]}.
  \end{equation*}
\end{df}
%%%
Hereafter, the integer $\ell$ in the definition above is referred to
as the degree of the UFG condition.
\par
%%% Theorem 2.3 of the published version (an L^p estimate from Kusuoka 2003)
%%% was removed in October 2026. It was used only in the proof of
%%% Theorem 3.1, which now uses Corollary 14 of Litterer and Lyons.
%%% We next introduce a theorem from Kusuoka~\cite{kusuoka2003malliavin},
%%% which concerns the smoothness of the semigroup
%%% $(P^X_t)_{t \geq 0}$
%%% in certain directions.
%%% %%
%%%
%%% \begin{thm}\label{sdeest}
%%%   Assume that the vector fields $V_0, V_1, \dots, V_d$ satisfy
%%%   the UFG condition and let $X$ be the diffusion process defined
%%%   in~\eqref{eq:SDE}.
%%%   Let $\alpha_1, \alpha_2, \dots, \alpha_k \in A^\ast_1$.  
%%%   Then there exists a constant $C > 0$ such that
%%%   \begin{equation*}
%%%     \|
%%%     V_{[\alpha_1]} \cdots V_{[\alpha_{k'}]} P^X_t V_{[\alpha_{k'+1}]}
%%%     \cdots
%%%     V_{[\alpha_k]} f
%%%     \|_{L^p(\mathrm{d}x)}
%%%     \leq C t^{-(\|\alpha_1\|
%%%       + \cdots + \|\alpha_k\|)/2} \|f\|_{L^p(\mathrm{d}x)}.
%%%   \end{equation*}
%%%   for $f \in C_{0}(\mathbb{R}^N), t\in (0,1]$, and $p\in [0,\infty]$ 
%%% \end{thm}
%%% %%% 
A condition on approximation operators is now introduced.
\begin{df}\label{MSF}
  Let $m\geq 2$.
  We say that a family $\{Q_{(s)}\}_{s>0}$ of Markov operators is
  \emph{$m$-moment similar} if there exist $K\in\mathbb{N}$ and random
  variables $L^{(1)},\dots,L^{(K)}$ with values in
  ${\mathcal L}_{\mathbb{R}}(A)$, defined on a common probability
  space and not depending on $s$, such that the following conditions
  hold:
  \begin{enumerate}[label=\textnormal{(\roman*)}]
  \item $j_m(L^{(k)})=L^{(k)}$ and
    $E[{\|L^{(k)}\|_{2}}^n]<\infty$
    for all $k\in\{1,\dots,K\}$ and $n\in\mathbb{N}$;
  \item $\sum_{k=1}^{K}\langle L^{(k)}, (0)\rangle=1$ almost surely;
  \item the moment condition
    \begin{equation*}
      E\big[j_{m}\big(\exp(L^{(1)})\exp(L^{(2)})\cdots
      \exp(L^{(K)})\big)\big]
      =E\bigg[j_{m}\bigg(\sum_{\alpha\in A^\ast}
      B^{\circ\alpha}(1)\,\alpha\bigg)\bigg]
    \end{equation*}
    holds;
  \end{enumerate}
  \rev{and, for all $s>0$, all continuous functions
  $f\colon\mathbb{R}^N\to\mathbb{R}$ of at most linear growth and all
  $x\in\mathbb{R}^N$,
  \begin{equation*}
    (Q_{(s)}f)(x)
    =E\Big[f\Big(\big(\exp(\Phi(\mathcal{D}_sL^{(K)}))\circ\cdots\circ
    \exp(\Phi(\mathcal{D}_sL^{(1)}))\big)(x)\Big)\Big].
  \end{equation*}}
  We refer to $m$ as the degree of the family.
\end{df}
\noindent
Here $B^{\circ\alpha}(t)$, $\alpha\in A^\ast$, are the iterated
Stratonovich integrals given by $B^{\circ\emptyset}(t)=1$ and
\begin{equation*}
  B^{\circ\alpha'a}(t)=\int_0^tB^{\circ\alpha'}(u)\circ\mathrm{d}B^a_u,
  \qquad \alpha'\in A^\ast,\ a\in A,
\end{equation*}
with $B^0_u=u$.
The flow of $\mathcal{D}_sL^{(1)}$, which corresponds to the leftmost
factor in condition~(iii), acts first
(see Corollary~12 of Kusuoka~\cite{kusuoka2004approximation}).
Since the coefficients of $V_0, V_1, \dots, V_d$ belong to
$C_b^{\infty}(\mathbb{R}^N;\mathbb{R}^N)$, the displacement of the
composed flows is bounded by a constant times
$\sum_{k=1}^{K}\|\mathcal{D}_sL^{(k)}\|_{2}$.
Hence, by condition~(i), the expectation defining $(Q_{(s)}f)(x)$ is
finite, and $Q_{(s)}f$ is again continuous and of at most linear
growth.
In particular, $Q_{(s)}f$ is defined for every Lipschitz continuous
function $f$, and compositions of such operators are well defined.
Definition~\ref{MSF} extends Definition~2 of
Kusuoka~\cite{kusuoka2004approximation}
(see also Definition~3 of Kusuoka~\cite{kusuoka2001approximation})
from a single exponential to a product of $K$ exponentials.
For $K=1$, it is equivalent to Definition~2 of
Kusuoka~\cite{kusuoka2004approximation} combined with the operator in
Theorem~3 of the same paper, with the perturbation term there set
to $0$.
% Added in November 2026: relation between degree and order.
\par
The degree $m$ is related to the order of
Definition~\ref{df:pth-order} as follows.
Under the UFG condition, Theorem~\ref{KLNVAppThm}~(b) below, applied
with the partitions $t_i=T(1-(1-i/\sharp\Delta)^{\gamma})$,
$\gamma>m-1$ (see Remark~\ref{regpart}), shows that the discretisation
algorithm defined by an $m$-moment similar family is of order
$(m-1)/2$.
For example, the N--V method has degree $5$ and order $2$
(see Remark~\ref{rem:m-moment-similar}~(d) and
Theorem~\ref{thm:N--V}).

\begin{rem}\label{rem:m-moment-similar}
  \begin{enumerate}[label=\textnormal{(\alph*)}]
  \item Since $\mathcal{D}_s$ is an algebra automorphism and,
    by the scaling property of Brownian motion,
    $(B^{\circ\alpha}(s))_{\alpha}$ has the same law as
    $(s^{\|\alpha\|/2}B^{\circ\alpha}(1))_{\alpha}$,
    condition~(iii) is equivalent to
    \begin{align*}
      &E\big[j_{m}\big(\exp(\mathcal{D}_sL^{(1)})\cdots
      \exp(\mathcal{D}_sL^{(K)})\big)\big]\\
      &\qquad=E\bigg[j_{m}\bigg(\sum_{\alpha\in A^\ast}
      B^{\circ\alpha}(s)\,\alpha\bigg)\bigg]
      \qquad\text{for every $s>0$.}
    \end{align*}
    Thus the random variables $\mathcal{D}_sL^{(k)}$,
    $k\in\{1,\dots,K\}$, satisfy the moment condition at every time step
    $s$, and $\sum_{k=1}^{K}\langle \mathcal{D}_sL^{(k)}, (0)\rangle=s$
    almost surely.
  \item Condition~(ii) restricts only the sum of the coefficients of
    the one-letter word $(0)$, so that a single factor may carry the
    whole drift, as in the N--V method.
    It corresponds to the normalisation of the coefficient of
    $(0)$ in Definition~2 of Kusuoka~\cite{kusuoka2004approximation}
    and in Definition~3 of Kusuoka~\cite{kusuoka2001approximation}.
    For $m\geq 4$, it follows from condition~(iii):
    with $c=\sum_{k=1}^{K}\langle L^{(k)}, (0)\rangle$,
    the coefficients of the words $(0)$ and $(0,0)$ in condition~(iii)
    give $E[c]=1$ and $E[c^2]=1$, and hence $c=1$ almost surely.
  \item The requirement $j_m(L^{(k)})=L^{(k)}$ in condition~(i) may be
    dropped if, in the definition of $Q_{(s)}$, each
    $\mathcal{D}_sL^{(k)}$ is replaced by its truncation
    $j_m(\mathcal{D}_sL^{(k)})$, as in Theorem~3 of
    Kusuoka~\cite{kusuoka2004approximation}.
  \item The N--V operator of Theorem~\ref{thm:N--V} is of this form
    with $K=d+1$, where each $L^{(k)}$ is a random multiple of a
    one-letter word.
    Let $\epsilon$ be a random variable with
    $P[\epsilon=1]=P[\epsilon=-1]=1/2$, independent of
    $(\eta_1,\dots,\eta_d)$, and choose $L^{(1)},\dots,L^{(d+1)}$ so
    that
    \begin{equation*}
      \big(\Phi(L^{(1)}),\dots,\Phi(L^{(d+1)})\big)=
      \begin{cases}
        (\eta_dV_d,\dots,\eta_1V_1,V_0) & \text{if $\epsilon=1$},\\
        (V_0,\eta_1V_1,\dots,\eta_dV_d) & \text{if $\epsilon=-1$}.
      \end{cases}
    \end{equation*}
    Since $\mathcal{D}_s$ multiplies $(0)$ by $s$ and the other
    one-letter words by $s^{1/2}$, the operator $Q_{(s)}$ of
    Definition~\ref{MSF} is then the operator~\eqref{NVOP}.
    Conditions~(i)--(iii) hold with $m=5$,
    both for standard normal $\eta_i$ and for the discrete
    $\hat{\eta}_i^{(\text{$5$th})}$, as one checks by expanding the
    exponentials (see also Ninomiya and Victoir~\cite{ninomiya2008weak}).
  \end{enumerate}
\end{rem}

%%%
%%%%%%%%%%%%%%%%%%%%%%%%%%%%%%%%%%%%%%%%%%%%%%%%%%
% Revised Theorem 2.5 of the published version (October 2026), with the
% semi-norm and the composition of operators introduced before it.
% Theorem 2.3 of the published version has been removed, and Remark 2.4
% and Lemma 2.6 are new. Hence this theorem keeps the number 2.5, and
% the later numbers in Section 2 increase by one (Theorem 2.9 for the
% N--V method, Definition 2.12 for the UH condition and Theorem 2.13
% for the estimate of Litterer and Lyons).
For $q\in\mathbb{N}$ and $f\in C_b^{\infty}(\mathbb{R}^N;\mathbb{R})$,
we define the semi-norm
\begin{equation*}
  \|f\|_{V,q}
  =\sum_{k=1}^{q}
  \sum_{\substack{\alpha_1,\dots,\alpha_k\in A^\ast_1\\
      \|\alpha_1\|+\dots+\|\alpha_k\|=q}}
  \|V_{[\alpha_1]}\cdots V_{[\alpha_k]}f\|_{\infty}
\end{equation*}
as in Kusuoka~\cite{kusuoka2004approximation}.
Since $A^\ast_1$ does not contain the one-letter word $(0)$, the vector
field $V_0$ does not appear by itself; it enters only through Lie
brackets.
Hereafter, we denote by
$(G_{1} \circ G_{2})g(x)$
the composition $(G_{1}(G_{2}g))(x)$ of $G_{1}$ and $G_{2}$,
where $G_{1}$ and $G_{2}$ are operators on the space of continuous
functions on $\mathbb{R}^N$ of at most linear growth,
such as the operators $Q_{(s)}$ of Definition~\ref{MSF}.
\par
The following theorem, which is essentially due to
Kusuoka~\cite{kusuoka2004approximation},
forms the foundation of high-order discretisation schemes;
see the comments after the theorem.
\begin{thm}
  \label{KLNVAppThm}
  Let $X$ be the diffusion process defined by~\eqref{eq:SDE},
  and let $\{Q_{(s)}\}_{s>0}$ be an $m$-moment similar family,
  with $K$ and $L^{(1)},\dots,L^{(K)}$ as in Definition~\ref{MSF}.
  \begin{enumerate}[label=\textnormal{(\alph*)}]
  \item There exist a constant $C>0$ and an integer $J\geq m+1$,
    depending only on $m$, $K$, the law of $(L^{(1)},\dots,L^{(K)})$
    and the vector fields $V_0, V_1, \dots, V_d$, such that
    \begin{equation*}
      \|P^{X}_s f - Q_{(s)} f\|_{\infty} \leq
      C\bigg(\sum_{q=m+1}^{J}s^{q/2}\|f\|_{V,q}
      +s^{(m+1)/2}\|\nabla f\|_{\infty}\bigg)
    \end{equation*}
    and
    \begin{equation*}
      \|P^{X}_s f - Q_{(s)} f\|_{\infty} \leq
      C s^{1/2}\|\nabla f\|_{\infty}
    \end{equation*}
    for all $s\in(0,1]$ and $f\in C_b^{\infty}(\mathbb{R}^N;\mathbb{R})$.
    Moreover, the second estimate holds with $\|\nabla f\|_{\infty}$
    replaced by $\|f\|_{\mathrm{Lip}}$ for all $s\in(0,1]$ and
    $f\in C_{\mathrm{Lip}}(\mathbb{R}^N;\mathbb{R})$.
  \item Assume in addition that the vector fields
    $V_0, V_1, \dots, V_d$ satisfy the UFG condition.
    Let $T>0$ be fixed, and consider a partition $\Delta$ of the interval
    $[0,T]$ with $\sharp\Delta=n$, given by $t_0=0<t_1<\dots<t_n=T$,
    and set $s_i=t_i-t_{i-1}$.
    Then there exists a positive constant $C'$, depending only on the
    constants in {\rm (a)}, the vector fields $V_0, V_1, \dots, V_d$
    and $T$, such that
    \begin{align*}
      &\|P^X_T f - Q_{(s_1)}\circ Q_{(s_2)}
        \circ\cdots\circ Q_{(s_n)} f\|_{\infty} \notag\\
      &\qquad \leq C'\|f\|_{\mathrm{Lip}}
        \bigg( {s_n}^{1/2} +\sum_{i=1}^{n-1}
        \frac{{s_i}^{(m+1)/2}}{(T - t_i)^{m/2}}\bigg)
    \end{align*}
    for all $f\in C_{\mathrm{Lip}}(\mathbb{R}^N;\mathbb{R})$.
  \end{enumerate}
\end{thm}
\noindent
In the terminology of Definition~1 of
Kusuoka~\cite{kusuoka2004approximation}, part~(a) states that
$\{Q_{(s)}\}_{s\in(0,1]}$ is $m$-similar; the third condition there,
$\|Q_{(s)}f\|_{\infty}\leq\exp(Cs)\|f\|_{\infty}$, holds for every
$C\geq 0$, since $\|Q_{(s)}f\|_{\infty}\leq\|f\|_{\infty}$ for the
Markov operator $Q_{(s)}$.
For $K=1$, part~(a) is Theorem~3 of
Kusuoka~\cite{kusuoka2004approximation}
with the perturbation term there set to $0$.
For $K\geq 2$, it follows from the same theorem applied to
$\xi=\log(\exp(L^{(1)})\cdots\exp(L^{(K)}))$, with the perturbation
term there given by the function $\zeta$ of Lemma~\ref{lem:product}
below.
The estimate for Lipschitz continuous functions holds because both
$E[\lvert X^x_s-x\rvert]$ and the expected displacement under $Q_{(s)}$
are bounded by a constant times $s^{1/2}$.
For $s\in(1,T]$, the estimates in part~(a) remain valid with constants
depending also on $T$, since then $s^{1/2}\geq 1$ and
$\|P^X_sf-Q_{(s)}f\|_{\infty}$ is bounded by a constant times
$\|\nabla f\|_{\infty}$.
Part~(b) follows from part~(a) by the telescoping argument in the proof
of Theorem~4 of Kusuoka~\cite{kusuoka2004approximation}.
The $i$th step with $i\leq n-1$ is estimated by applying part~(a) to
$P^X_{T-t_i}f$ and using
\begin{equation*}
  \|P^X_t f\|_{V,q}\leq C t^{-(q-1)/2}\|\nabla f\|_{\infty},
  \qquad t\in(0,T],\ q\in\{m+1,\dots,J\}.
\end{equation*}
For $t\in(0,1]$, this estimate is Theorem~2 of
Kusuoka~\cite{kusuoka2003malliavin};
see also Theorem~5 of Kusuoka~\cite{kusuoka2004approximation}.
For $t\in(1,T]$, it follows from the case $t=1$ applied to
$P^X_{t-1}f$, since $P^X_t=P^X_1P^X_{t-1}$ and
$\|\nabla P^X_{t-1}f\|_{\infty}\leq C\|\nabla f\|_{\infty}$,
with $C$ depending also on $T$.
If $s_i>T-t_i$, the second estimate in part~(a) is used instead.
This is why only the terms with the exponent $m/2$ appear in part~(b).
The last step gives the term ${s_n}^{1/2}$, and the extension from
$C_b^{\infty}(\mathbb{R}^N;\mathbb{R})$ to
$C_{\mathrm{Lip}}(\mathbb{R}^N;\mathbb{R})$ is obtained by mollification
and truncation.
%
% The proof of the following lemma was a sketch in the October 2026
% revision and is written out in full here (November 2026). It extends
% Proposition 9, Lemma 11 and Corollary 12 of Kusuoka (2004) from two
% flows to K flows; the last step is that of Lemma 14 and Corollary 15
% there.
\rev{\begin{lem}\label{lem:product}
  Let $L^{(1)},\dots,L^{(K)}$ be as in Definition~\ref{MSF}, and let
  \begin{equation*}
    \xi=\log\big(\exp(L^{(1)})\cdots\exp(L^{(K)})\big).
  \end{equation*}
  Then $\xi$ is $m$-$\mathcal{L}$-moment similar in the sense of
  Definition~2 of Kusuoka~\cite{kusuoka2004approximation}.
  Moreover, let
  \begin{equation*}
    \zeta(s,x)=\big(\exp(\Phi(\mathcal{D}_sL^{(K)}))\circ\cdots\circ
    \exp(\Phi(\mathcal{D}_sL^{(1)}))\big)(x)
    -\exp(\Phi(j_m(\mathcal{D}_s\xi)))(x)
  \end{equation*}
  for $s\in(0,1]$ and $x\in\mathbb{R}^N$.
  Then $\zeta$, regarded as a function of $s$, $x$ and the sample
  point, is measurable and continuous in $x$, and satisfies
  \begin{equation*}
    \sup_{s\in(0,1],\,x\in\mathbb{R}^N}s^{-(m+1)/2}
    E[\lvert \zeta(s,x)\rvert]<\infty.
  \end{equation*}
\end{lem}
\begin{proof}
  Since every word in $\exp(L^{(1)})\cdots\exp(L^{(K)})-1$ has order at
  least one, we have
  \begin{equation*}
    j_m(\xi)=j_m\bigg(\sum_{r=1}^{m}\frac{(-1)^{r-1}}{r}
    \Big(j_m\big(\exp(L^{(1)})\cdots\exp(L^{(K)})\big)-1\Big)^{r}\bigg).
  \end{equation*}
  Hence the coefficients of $j_m(\xi)$ are polynomials of degree at most
  $m$ in the coefficients of $L^{(1)},\dots,L^{(K)}$, and
  $\|j_m(\xi)\|_{2}\leq C(1+\sum_{k=1}^{K}\|L^{(k)}\|_{2})^{m}$.
  Here and below, $C$ denotes a constant depending only on $m$, $K$,
  $d$, $N$ and bounds on the derivatives of $V_0, V_1, \dots, V_d$.
  By condition~(i), $E[\|j_m(\xi)\|_{2}^{n}]<\infty$ for all
  $n\in\mathbb{N}$.
  By the Baker--Campbell--Hausdorff formula, $\xi$ is a Lie series, and
  $\xi-\sum_{k=1}^{K}L^{(k)}$ consists of iterated Lie brackets of
  $L^{(1)},\dots,L^{(K)}$, which contain only words of length at least
  two.
  Hence $\langle\xi,(0)\rangle=\sum_{k=1}^{K}\langle L^{(k)},(0)\rangle=1$
  almost surely by condition~(ii).
  Since $\exp(\xi)=\exp(L^{(1)})\cdots\exp(L^{(K)})$, condition~(iii) of
  Definition~\ref{MSF} is exactly the moment condition in Definition~2 of
  Kusuoka~\cite{kusuoka2004approximation}.
  Thus $\xi$ is $m$-$\mathcal{L}$-moment similar.
  \par
  For each sample point, the flows depend continuously on $s$ and on
  the initial point $x$, and for each $(s,x)$, $\zeta(s,x)$ is a random
  variable.
  Hence $\zeta$ is measurable and continuous in $x$.
  \par
  To estimate $\zeta$, we first extend Lemma~11 and Corollary~12 of
  Kusuoka~\cite{kusuoka2004approximation} from two flows to $K'$ flows.
  Let $z_1,\dots,z_{K'}\in{\mathcal L}_{\mathbb{R}}(A)$, and let
  $g\in C^{\infty}(\mathbb{R}^N;\mathbb{R})$ be such that $\Phi(\beta)g$
  is bounded for every word $\beta\in A^\ast$ with
  $\lvert\beta\rvert\geq 1$.
  We claim that, for every integer $n\geq 0$ and every
  $x\in\mathbb{R}^N$,
  \begin{align}
    &\bigg\lvert g\Big(\big(\exp(\Phi(z_{K'}))\circ\cdots\circ
      \exp(\Phi(z_1))\big)(x)\Big)
      -\sum_{\ell_1+\dots+\ell_{K'}\leq n}
      \frac{\big(\Phi(z_1^{\ell_1}\cdots z_{K'}^{\ell_{K'}})g\big)(x)}
      {\ell_1!\cdots\ell_{K'}!}\bigg\rvert \notag\\
    &\qquad\leq\sum_{1\leq k_1\leq\dots\leq k_{n+1}\leq K'}
      \|\Phi(z_{k_1}\cdots z_{k_{n+1}})g\|_{\infty}. \tag{$\ast$}\label{eq:Kflows}
  \end{align}
  For $K'=1$, this is Proposition~9 of
  Kusuoka~\cite{kusuoka2004approximation}, that is, Taylor's formula
  along the flow of $\Phi(z_1)$.
  For $K'\geq 2$, let $y$ be the image of $x$ under the composition of
  the first $K'-1$ flows.
  We apply Taylor's formula of order $n$ along the flow of
  $\Phi(z_{K'})$ at $y$, and then~\eqref{eq:Kflows} for the first $K'-1$
  flows, with $n$ replaced by $n-\ell_{K'}$, to the functions
  $\Phi(z_{K'}^{\ell_{K'}})g/\ell_{K'}!$, $\ell_{K'}=0,\dots,n$.
  Since $\Phi$ is a homomorphism, the main terms add up to the sum
  in~\eqref{eq:Kflows}.
  Each remainder is bounded by a term on the right-hand side
  of~\eqref{eq:Kflows} with a coefficient at most one, and each
  nondecreasing sequence $(k_1,\dots,k_{n+1})$ occurs at most once.
  \par
  Fix $l\in\{1,\dots,N\}$ and let $\pi_l(x)=x_l$.
  Since the coefficients of $V_0, V_1, \dots, V_d$ belong to
  $C_b^{\infty}(\mathbb{R}^N;\mathbb{R}^N)$, $\Phi(\beta)\pi_l$ is
  bounded for $\lvert\beta\rvert\geq 1$.
  We apply~\eqref{eq:Kflows} with $g=\pi_l$, $n=m$, $K'=K$ and
  $z_k=\mathcal{D}_sL^{(k)}$, where $s\in(0,1]$.
  By condition~(i), $z_k$ is a linear combination of words $\alpha$ with
  $1\leq\|\alpha\|\leq m$, with coefficients
  $s^{\|\alpha\|/2}\langle L^{(k)},\alpha\rangle$.
  Hence a product $z_{k_1}\cdots z_{k_{m+1}}$ is a linear combination of
  words $\beta$ with $\|\beta\|\geq m+1$ and
  $\lvert\beta\rvert\leq m(m+1)$, whose coefficients contain the factor
  $s^{\|\beta\|/2}\leq s^{(m+1)/2}$.
  The same holds for $z_1^{\ell_1}\cdots z_K^{\ell_K}-j_m(z_1^{\ell_1}\cdots
  z_K^{\ell_K})$ with $\ell_1+\dots+\ell_K\leq m$.
  Moreover, every word in $z_1^{\ell_1}\cdots z_K^{\ell_K}$ has order at
  least $\ell_1+\dots+\ell_K$, so that
  \begin{equation*}
    j_m\big(\exp(z_1)\cdots\exp(z_K)\big)
    =\sum_{\ell_1+\dots+\ell_K\leq m}
    \frac{j_m(z_1^{\ell_1}\cdots z_K^{\ell_K})}{\ell_1!\cdots\ell_K!}.
  \end{equation*}
  Since $\exp(z_1)\cdots\exp(z_K)=\mathcal{D}_s(\exp(\xi))
  =\exp(\mathcal{D}_s\xi)$, we obtain
  \begin{align*}
    &\pi_l\Big(\big(\exp(\Phi(\mathcal{D}_sL^{(K)}))\circ\cdots\circ
      \exp(\Phi(\mathcal{D}_sL^{(1)}))\big)(x)\Big)\\
    &\qquad=\Big(\Phi\big(j_m(\exp(\mathcal{D}_s\xi))\big)\pi_l\Big)(x)
      +R_1(s,x).
  \end{align*}
  Similarly, \eqref{eq:Kflows} with $K'=1$ and
  $z_1=j_m(\mathcal{D}_s\xi)=\mathcal{D}_s(j_m(\xi))$ gives
  \begin{equation*}
    \pi_l\big(\exp(\Phi(j_m(\mathcal{D}_s\xi)))(x)\big)
    =\Big(\Phi\big(j_m(\exp(\mathcal{D}_s\xi))\big)\pi_l\Big)(x)
    +R_2(s,x).
  \end{equation*}
  Here we have used the identity $j_m(\exp(j_m(Z)))=j_m(\exp(Z))$
  for every $Z\in\mathbb{R}\langle\langle A\rangle\rangle$ without
  constant term.
  It holds since the words of order greater than $m$ in $Z$ contribute
  to $\exp(Z)$ only through words of order greater than $m$.
  By the above observations on the words and their coefficients, the
  remainders satisfy
  \begin{align*}
    \lvert R_1(s,x)\rvert
    &\leq C s^{(m+1)/2}\Big(1+\sum_{k=1}^{K}\|L^{(k)}\|_{2}\Big)^{m+1},\\
    \lvert R_2(s,x)\rvert
    &\leq C s^{(m+1)/2}\big(1+\|j_m(\xi)\|_{2}\big)^{m+1}.
  \end{align*}
  The $l$th component of $\zeta(s,x)$ is $R_1(s,x)-R_2(s,x)$.
  Therefore, by the bound on $\|j_m(\xi)\|_{2}$,
  \begin{equation*}
    E[\lvert\zeta(s,x)\rvert]
    \leq C s^{(m+1)/2}
    E\bigg[\Big(1+\sum_{k=1}^{K}\|L^{(k)}\|_{2}\Big)^{m(m+1)}\bigg]
  \end{equation*}
  for all $s\in(0,1]$ and $x\in\mathbb{R}^N$, and the right-hand side is
  finite by condition~(i).
  For $K=2$, this estimate is contained in Corollary~15 of
  Kusuoka~\cite{kusuoka2004approximation}.
\end{proof}}

\subsubsection{Intuition behind the framework of high-order discretisation schemes}
The core idea of high-order weak approximation is to
match the stochastic Taylor expansion of the diffusion process up to a
prescribed degree~$m$. 
This is achieved by constructing approximation
operators~$Q_{(s)}$ whose moments coincide with those of the associated semigroup
$P_s^{X}$.
When these moment conditions are satisfied, the resulting operator
provides a weak approximation with controlled discretisation error.
The cubature on Wiener space method,
introduced in Kusuoka~\cite{kusuoka2001approximation} and
in Lyons and Victoir~\cite{lyons2004cubature},
is a precise and practically applicable framework based on this idea,
and its corresponding approximation result is formulated
in Theorem~\ref{KLNVAppThm}.
\par
It replaces the Wiener measure with a finitely supported probability measure that satisfies the moment conditions up to degree~$m$. 
To make this construction algebraically tractable, the free Lie algebra
$A^{\ast}$ is employed to index both the differential operators and the
iterated Wiener integrals appearing in the stochastic Taylor expansion.
The Lie polynomials $\mathcal{L}_{\mathbb{R}}(A)$ and their
exponential map represent the flows generated by
vector fields. 
Within this algebraic framework, one considers Lie
polynomial-valued random variables, denoted by $L^{(k)}$ in
Definition~\ref{MSF}, which may be interpreted as ODE-valued random variables.
By appropriately weighting and composing the associated exponential flows, one
constructs an $m$-moment similar family $\{Q_{(s)}\}_{s>0}$.\par
Based on this representation, numerical schemes are implemented by sampling random
Lie polynomials, computing the associated flows, and composing them along a time
partition.  
In this way, high-order weak schemes reduce SDE simulation to solving ODEs with
random coefficients specifically designed to preserve the moment structure.
The N--V method, introduced in the next subsection, is a leading practical example following this principle.
\begin{rem}
  Both the moment similar families of Kusuoka~\cite{kusuoka2001approximation},
  which consist of random variables,
  and the cubature on Wiener space
  of Lyons and Victoir~\cite{lyons2004cubature} aim to match the joint moments of iterated Wiener integrals up to a prescribed degree $m$.
They differ in that the former approximates the joint law of iterated integrals using suitably constructed random variables, whereas the latter replaces the Wiener measure with a discrete measure supported on bounded-variation paths.

% Revised in November 2026: the random variables L_{(s)} and
% L^{(k)}_{(s)} of the published version are written with the dilation.
Motivated by this perspective,
Kusuoka~\cite{kusuoka2004approximation} introduced a 
Lie-polynomial-based formulation, in which a single random Lie series
$\xi$ is scaled by the dilation $\mathcal{D}_s$ and the exponential map
$\exp(\mathcal{D}_s\xi)$ yields moment matching directly via the flow
(see Definition~2 and Theorem~3 there).
Following this development, practical high-order schemes such as the 
N--V \cite{ninomiya2008weak} and N--N \cite{ninomiya2009new} methods were constructed by composing multiple 
exponential maps,
\[
  \exp(\mathcal{D}_sL^{(1)})
  \exp(\mathcal{D}_sL^{(2)})
  \cdots
  \exp(\mathcal{D}_sL^{(K)}),
\]
where $L^{(1)},\dots,L^{(K)}$ are random variables as in
Definition~\ref{MSF}, chosen so that the moment condition~(iii) there
holds.
In this paper, we adopt the notion of
$m$-moment similar families of operators (Definition~\ref{MSF}) as a natural formulation for the approximation result stated in Theorem~\ref{KLNVAppThm}.
We also note that the third-order scheme 
developed in Shinozaki~\cite{Shinozaki2017} is obtained directly
via the construction of a
moment similar family of random variables in the sense of
Kusuoka~\cite{kusuoka2001approximation}, rather than through compositions
of exponential flows.
\end{rem}

\begin{rem}\label{regpart}
In the original works of Kusuoka
\cite{kusuoka2001approximation,kusuoka2004approximation}
and of Lyons and Victoir~\cite{lyons2004cubature}, 
the approximation theorem (Theorem~\ref{KLNVAppThm}~(b))
is formulated under the assumption  
that $f \in C_b^{\infty}(\mathbb{R}^N)$.
Therein, such a specific class of
non-uniform time partitions is also employed that
$t_i = T(1-(1-i/\sharp \Delta)^{\gamma})$,
with the parameter $\gamma$ being determined by
the prescribed approximation order. 
Furthermore, the smoothness requirement on $f$ can be relaxed to 
$f \in C_{\mathrm{Lip}}(\mathbb{R}^N)$, as explicitly discussed
in Section~3 (b) of Lyons and Victoir\cite{lyons2004cubature}.

Later numerical studies of
Ninomiya~\cite{ninomiya2003new};
Ninomiya and Victoir~\cite{ninomiya2008weak};
Ninomiya and Ninomiya~\cite{ninomiya2009new};
and Ninomiya\cite{ninomiya2010application}
suggest that such structural assumptions on the time partition are not essential in practice: 
even under uniform time partition, the approximation errors exhibit the expected 
convergence behaviour.

Motivated by these observations, Kusuoka~\cite{kusuoka2013gaussian}
further established 
a theoretical justification of  a specific class of schemes,
including N--V \cite{ninomiya2008weak} and N--N \cite{ninomiya2009new}
schemes, with a uniform time partition 
and for Lebesgue-measurable functions $f$, thereby extending the 
applicability of the framework well beyond the original smooth setting.
\end{rem}
\subsection{N--V method}\label{nvsec}
We introduce the N--V method from
Ninomiya and Victoir~\cite{ninomiya2008weak,kusuoka2013gaussian}.
It is a discretisation method of order $2$.
\begin{thm}\label{thm:N--V}
  Let $\eta_{1}, \eta_{2}, \dots , \eta_{d}$
  be independent standard normal
  random variables and define
  \begin{align}
    Q_{(s)}^{(\NV)}f(x)
    &= E\bigg[\frac{1}{2}
      f\big(\exp(sV_{0})\exp(\sqrt{s}\eta_{1}V_{1})
      \cdots
      \exp(\sqrt{s}\eta_{d}V_{d})(x)\big) \notag\\
    &\quad +\frac{1}{2}f\big(\exp(\sqrt{s}\eta_{d}V_{d})
      \cdots \exp(\sqrt{s} \eta_{1}V_{1})
      \exp(sV_{0})(x)\big)
      \bigg]. \label{NVOP}
  \end{align}
  Then $(Q^{(\NV)}_{(s)})_{s>0}$
  is a discretisation method of order $2$. 
\end{thm}
\noindent
Replacing the independent standard normal random variables
$(\eta_i)_{i\in\{1,\dots,d\}}$ in the definition above with the
independent discrete random variables
$(\hat{\eta}_i^{(\text{$5$th})})_{i\in\{1,\dots,d\}}$, we define the
discrete N--V method as in
Ninomiya and Shinozaki~\cite{NinomiyaShinozaki2019}.  
Each
$\hat{\eta}_i^{(\text{$5$th})}$ is a discrete random variable on
$(\hat{\Omega}, \hat{\mathcal{F}}, \hat{P})$ defined by
\begin{equation*}
  \hat{P}[\hat{\eta}_i^{(\text{$5$th})}=\pm \sqrt{3}]=
  \frac{1}{6},
  \;\;\hat{P}[\hat{\eta}_i^{(\text{$5$th})}=0]=\frac{2}{3}.
\end{equation*}
We note that the moments of $\hat{\eta}_i^{(\text{$5$th})}$ up to the
$5$th order match those of the standard normal distribution.
\par
As in subsection~\ref{tbsec}, we denote by $\mu_{t_i}^{(\Alg, \Delta)} = \sum_{j=1}^{M_i} w_{j}^{(i)}
\delta_{p_{j}^{(i)}}$ the probability measure on
$\mathbb{R}^{N}$ induced by $X^{(\Alg, \Delta)}_{t_i}$
, where
$\big((p_{j}^{(i)}, w_{j}^{(i)})\big)_{j\in\{1,\dots,M_i\}}$ is the set of pairs
consisting of points $p_j^{(i)}$ in $\mathbb{R}^N$ and weights
$w_j^{(i)}$ assigned to them.
%%%
Let $(\hat{\eta}^{(5{\rm th})}_{i,l})_{\substack{i\in\{1,\dots,
    \sharp \Delta\}\\l\in\{1,\dots, d\}}}$ be a set of independent and identically
distributed discrete random variables following the same distribution
as $\hat{\eta}_i^{(5{\rm th})}$.
%%%%%%%%
Let $\phi$ be a function from $\{\pm \sqrt{3}, 0\}^d$ to $[0,1]$ defined by
\begin{equation*}
  \phi(\eta)=\bigg(\frac{1}{6}\bigg)^{\card\{i\,:\, \eta_i =\pm\sqrt{3}\}}
  \bigg(\frac{2}{3}\bigg)^{\card\{i \,:\, \eta_i = 0\}}
\end{equation*}
for $\eta = (\eta_1, \dots, \eta_d) \in \{\pm \sqrt{3}, 0\}^d$.
We also define functions $y^{(+)}$ and $y^{(-)}$ from
$\mathbb{R}^{N} \times \{\pm \sqrt{3}, 0\}^d \times \mathbb{R}_{>0}$
to $\mathbb{R}^{N}$ as
\begin{align*}
  y^{(+)}(x, \eta, s) &=
                        \exp(s V_{0})
                        \exp(\sqrt{s} \eta_1 V_{1}) \cdots
                        \exp(\sqrt{s} \eta_d V_{d})(x),\\
  y^{(-)}(x, \eta, s) &=
                        \exp(\sqrt{s} \eta_d V_{d}) \cdots
                        \exp(\sqrt{s} \eta_1 V_{1})
                        \exp(s V_{0})(x).
\end{align*}
%%%%%%%%% 
Then, the operator $\hat{Q}_{(s_i)}^{(\NV, \Delta)}$ between
discrete measures is defined as 
\begin{align}
  &\hat{Q}_{(s_i)}^{(\NV, \Delta)}
    \bigg(\sum_{j=1}^{M_i} w_{j}^{(i)} \delta_{p_{j}^{(i)}}\bigg) \notag\\
  &\qquad=\sum_{j=1}^{M_i}\sum_{\eta\in\{\pm\sqrt{3}, 0\}^d} w_{j}^{(i)} \frac{\phi(\eta)}{2}
    \big( \delta_{y^{(+)}(p_j^{(i)},\eta, s_i)}+\delta_{y^{(-)}(p_j^{(i)}, \eta, s_i)}\big).
    \label{eq:nvalg}
\end{align}
%%%%%%%%% 
In this way, we construct a concrete example of
$\hat{Q}^{(\Alg, \Delta)}_{(s_{i})}$ used in this paper,
as discussed in subsection~\ref{tbsec}.
We then apply the recombination algorithm
described in subsection~\ref{recsubsec} to this operator and obtain
$\mu_{T}^{(\NV, \Delta, {\rm Rec})}$.
\par
We similarly define $\hat{Q}_{(s_i)}(\mu)$ for an $m$-moment similar
family $\{Q_{(s)}\}_{s>0}$, using discrete random variables
$\hat{\eta}^{(m\text{th})}$ whose moments up to the $m$th order match
those of the standard normal random variables, provided that the
family obtained in this way is again $m$-moment similar
(for the N--V method, see Remark~\ref{rem:m-moment-similar}~(d)).
In each case, the operator $\hat{Q}_{(s_i)}$ on discrete measures is
the dual of the operator $Q_{(s_i)}$ on functions, where $Q_{(s_i)}$ is
defined with the discrete random variables, that is,
\begin{equation*}
  \hat{Q}_{(s_i)}(\mu)(f)=\mu(Q_{(s_i)}f)
\end{equation*}
for every discrete probability measure $\mu$ on $\mathbb{R}^N$ and
every $f\in C_{\mathrm{Lip}}(\mathbb{R}^N;\mathbb{R})$.
For example, \eqref{eq:nvalg} is the dual of \eqref{NVOP} with $\eta_i$
replaced by $\hat{\eta}_i^{(\text{$5$th})}$.
In this paper, we consider the discretisation algorithm defined by an
$m$-moment similar family and its associated recombination measure.

\begin{rem}
  In Ninomiya and Victoir~\cite{ninomiya2008weak},
  the authors propose implementing the N--V
  method using Bernoulli random variables in addition to the standard
  normal random variables $\eta_{1}, \eta_{2}, \dots , \eta_{d}$
  in order to halve the number of ODEs to be solved.
  In contrast, the present paper employs the N--V method without utilising
  Bernoulli random variables instead calculating the expectation as shown in
  equation~\eqref{NVOP}.
  This approach is motivated by the need to
  minimise the number of support points in the discrete N--V method so
  as to regulate the overall computational cost, as noted in
  Remark~\ref{NVad}.
\end{rem}

%%%%%%%%%%%%%%%%%%%%%%%%%%%%%%%%%%%%%%%%%%%%%%%%%%
\begin{rem}
  \label{NVad}
  The N--V method offers the following advantages over other high-order discretisation methods, such as the methods in 
  Ninomiya and Ninomiya~\cite{ninomiya2009new}
  and Shinozaki~\cite{Shinozaki2017}.
  As stated in the original paper
  Ninomiya and Victoir~\cite{ninomiya2008weak},
  it is often more efficiently implementable because analytical solutions
  to the ODEs are available in certain cases, such as the Heston model
  Heston~\cite{heston1993closed}.
  Furthermore, the Feller condition
  Feller~\cite{feller1951two} guarantees the positivity of processes
  discretised by the N--V method.
  
  \par
  Another advantage is that, when implementing high-order
  discretisation methods using discrete random variables, the N--V
  method is economical in the sense that only $d$ random numbers are
  required to approximate an SDE driven by $d$-dimensional Brownian
  motion.
  %%%%%%
  The numbers of random variables involved and the cardinalities of the supports of the resulting measures are summarised in Table~\ref{tab:NUMS}.

\end{rem}    

\subsection{High-order recombination: approximation theorem}
\label{recsec}
In this subsection, we introduce the approximation theorem established in
Theorem~19 of Litterer and Lyons~\cite{litterer2012high}.
\par
\begin{df}\label{UH}
  The vector fields $V_0, V_1, \dots, V_d$ satisfy the uniform
  H\"{o}rmander (UH) condition if there exists $\ell \in \mathbb{N}$
  such that 
  \begin{equation*}
    \inf_{x, y \in \mathbb{R}^{N}, |y|=1} \bigg(\sum_{\alpha \in A^\ast_1(\ell)}
    \langle V_{[\alpha]} I_N (x),\, y\rangle^2\bigg)>0.
  \end{equation*}
\end{df}
%%%%%%%%%%%%%%%%%%%%%%%%%%%%%%%%%%%%%
\begin{thm}\label{ori:est}
  Assume that the vector fields $V_0, V_1, \dots, V_d$ satisfy the UH
  condition.
  Let $\Alg$ be the discretisation algorithm defined by an $m$-moment
  similar family, and let
  $\Delta$ be a partition of $[0, T]$ as
  $0=t_0<t_1<\dots<t_{\sharp\Delta}=T$.
  Let
  \begin{equation*}
    \big(\mu_{t_i}^{(\Alg, \Delta, {\rm Rec})},
    (U_{k}^{(i)})_{k\in \{1,\dots, l_i\}}\big)_{i\in\{0,\dots,\sharp\Delta\}}
  \end{equation*}
  be a set of recombination measures with initial value
  $x_0\in\mathbb{R}^N$ and patch divisions of their supports.
  %%%%%%%%%%%%%%%%%%%%%%%%%%%%%%%%%%%%%%%%%%%%%%%%%% 
  Then the following inequality holds:
  \begin{align*}
  & \big\lvert P^{X}_T f(x_0) -\mu_{T}^{(\Alg, \Delta, {\rm Rec})}(f) \big\rvert \notag\\
  &\qquad\leq \bigg( C_1(T) \bigg( {s_{\sharp\Delta}}^{1/2} +
  \sum_{i=1}^{\sharp\Delta-1} \sum_{j=m}^{m+1}
  \frac{{s_i}^{(j+1)/2}}{(T - t_i)^{j/2}} \bigg) \notag\\
  &\qquad\qquad\qquad\quad + C_2(T) \sum_{i=2}^{\sharp\Delta-1}
  \frac{{u_i}^{m+1}}{(T - t_i)^{m\ell/2}} \bigg) \|f\|_{\rm{Lip}}.
\end{align*}
  where
  $u_i=\max\{
  \radius(U_{k}^{(i)}) \,:\, k=1,\dots, l_i
  \}$
  and $s_i= t_i-t_{i-1}$.
\end{thm}
%%%%%%%%%%%%%%%%%%%%%%%%%%%%%%%%%%%%%%%%%%%%%%%%%%
In the following, we present a sufficient condition that satisfies
the assumptions of the above theorem and can be used to determine the
patch division.
\begin{exa}\label{ori:cond}
  Given a time step as 
  \begin{equation*}
    t_j = T\bigg(1-\bigg(1 - \frac{j}{\sharp \Delta} \bigg)^{\gamma}
    \bigg)
  \end{equation*}
  for some $\gamma>m-1$, as in Example~21 of Litterer and
  Lyons~\cite{litterer2012high}.
  If we set
  \begin{equation*}
    u_i = \bigg(
    \frac{{s_i}^{m+1}}{(T-t_i)^{m(1-\ell)}}
    \bigg)^{1/(2(m+1))},
  \end{equation*}
  then we have
  \begin{equation*}
    \|
    P^{X}_T f-\mu_{T}^{(\Alg,\Delta,{\rm Rec})}(f)
    \|_{\infty}
    \leq \frac{C_2(T)}{(\sharp \Delta)^{(m-1)/2}}\|f\|_{\rm{Lip}}. 
  \end{equation*}
  In other words, the recombination does not compromise the order of the
  discretisation error.
\end{exa}
We refer to the condition introduced above by Example~21 of Litterer--Lyons as the
Litterer--Lyons patch condition, or simply the LL patch condition.
\section{Refined error estimation and patch division algorithms}
\label{sec:refined}
In this section, we present the main results of this paper. First, we
refine Theorem~\ref{ori:est} and Example~\ref{ori:cond} by extending
the argument in Litterer and Lyons~\cite{litterer2012high}. This leads to
Theorem~\ref{REE} and Example~\ref{pel}. The latter provides a sufficient
condition under which the patch division satisfies the required bound
on the recombination error.
\subsection{Refined error estimation}\label{3-2}
Taking the weight of each patch into account yields Theorem~\ref{REE} below,
which refines Theorem~\ref{ori:est} by incorporating not only the radii of the
patches but also their weights.
\begin{thm}\label{REE}
  Assume that the vector fields $V_0, V_1, \dots, V_d$ satisfy the UH
  condition, \rev{and let $\ell$ be the integer in Definition~\ref{UH}.}
  Let $\Alg$ be the discretisation algorithm defined by an $m$-moment
  similar family (see subsection~\ref{nvsec}), and let $\Delta$
  be a partition of $[0, T]$
  as $0 = t_0 < t_1 < \dots < t_{\sharp\Delta} = T$.  
  Let
  \begin{equation*}
    \big(\mu_{t_i}^{(\Alg, \Delta, {\rm Rec})}\big)_{i\in\{0,\dots,\sharp\Delta\}}
  \end{equation*}
  be recombined measures with initial value $x_0\in\mathbb{R}^N$,
  constructed by~\eqref{eq:high-order-recombination}
  with patch divisions $(U_{k}^{(i)})_{k\in\{1,\dots,l_i\}}$ of
  $\support(\nu^{(i)})$, $i\in\{1,\dots,\sharp\Delta-1\}$.
  Then there exist constants $C_1(T), C_2(T)>0$ and an integer
  $J\geq m+1$ such that, with $s_i = t_i - t_{i-1}$,
  the following inequality holds \rev{for all
  $f\in C_b^{\infty}(\mathbb{R}^N;\mathbb{R})$}:
  \begin{align}
    &\big\lvert
      P^{X}_T f(x_0) -\mu_{T}^{(\Alg, \Delta, {\rm Rec})}(f)
      \big\rvert \notag\\
    &\leq \bigg(C_1(T)
      \bigg({s_{\sharp\Delta}}^{1/2} + \sum_{i=1}^{\sharp\Delta - 1}
      \sum_{q=m+1}^{J}\frac{{s_i}^{q/2}}{(T - t_{i})^{(q-1)/2}}\bigg) \notag\\
    &\qquad\qquad\qquad + C_2(T)\sum_{i=1}^{\sharp\Delta - 1}
      \frac{\sum_{k=1}^{l_i}\big(\radius(U_k^{(i)})\big)^{m+1}
      \weight(U_k^{(i)})}{(T - t_i)^{m\rev{\ell}/2}} \bigg)
      \|\nabla f\|_{\infty}. \label{eq:refinedErr}
  \end{align}
\end{thm}
%%%%%%%%%%%%%%%%%%%%%%%%%%%%%%%%%%%%%%%%%%%%%%%%%%
\noindent Hereafter, we refer to the first term of the right-hand side
of~\eqref{eq:refinedErr}
\begin{equation*}
  C_1(T)
  \bigg(s_{\sharp \Delta }^{1/2} + \sum_{i=1}^{\sharp\Delta-1}
  \sum_{q=m+1}^{J} \frac{{s_i}^{q/2}}{(T-t_{i})^{(q-1)/2}}
  \bigg)\|\nabla f\|_{\infty}
\end{equation*}
as the discretisation error and the second term
\begin{equation*}
  C_2(T) \sum_{i=1}^{\sharp\Delta-1}
  \frac{\sum_{k=1}^{l_i}
    \big(\radius(U_k^{(i)})
    \big)^{m+1}\weight(U_k^{(i)})}
       {(T-t_i)^{m\ell/2}}\|\nabla f\|_{\infty}
\end{equation*}
as the recombination error.
%%%%%%%%%%%%%%%%%%%%%%%%%%%%%%%%%%%%%%%%%%%%%%%%%%

\begin{proof}
Write $n=\sharp\Delta$, and let $\{Q_{(s)}\}_{s>0}$ be the $m$-moment
similar family defining $\Alg$.
Then $\hat{Q}^{(\Alg,\Delta)}_{(s_i)}$ is the dual of $Q_{(s_i)}$
(see subsection~\ref{nvsec}).  
By Theorem~\ref{KLNVAppThm}~(a) and the comments following it,
there exist constants $c_1, c_3>0$ and an integer $J\geq m+1$
such that 
\begin{equation}\label{K1}
  \|P^{X}_s f - Q_{(s)}f\|_{\infty} 
  \leq c_1\bigg(s^{{(m+1)}/{2}}\|\nabla f\|_{\infty}
  +\sum_{q=m+1}^{J}s^{{q}/{2}}\|f\|_{V, q}\bigg)
\end{equation}
and  
\begin{equation}\label{K3}
  \|P^{X}_s f - Q_{(s)}f\|_{\infty}
  \leq c_3 s^{1/2}\|\nabla f \|_{\infty}
\end{equation}
hold for $s\in(0,T]$ and $f\in C_b^{\infty}(\mathbb{R}^N;\mathbb{R})$,
where $\|\cdot\|_{V,q}$ is the semi-norm defined before
Theorem~\ref{KLNVAppThm}.
Furthermore, the UH condition implies the UFG condition
(see Kusuoka~\cite{kusuoka2001approximation}).
Hence the estimate in the comments following
Theorem~\ref{KLNVAppThm}, which is based on Theorem~2 of
Kusuoka~\cite{kusuoka2003malliavin}, yields
the following estimate for $t\in[0,T)$ and $q\in\{m+1,\dots,J\}$:
\begin{equation}\label{K2}
  \| P^X_{T-t} f \|_{V, q}
  \leq c_2 \frac{\|\nabla f \|_{\infty}}{(T-t)^{(q-1)/2}}.
\end{equation}
% Added in November 2026: the UFG condition is needed globally for (3.4).
Note that the UFG condition is used here on the whole of
$\mathbb{R}^N$, since $\|\cdot\|_{V,q}$ is defined by suprema over
$\mathbb{R}^N$.
%It is not implied by a UH condition that holds only on a subset of
%$\mathbb{R}^N$, and it is therefore assumed separately in
%Corollary~\ref{locest} below.
\par
Define $\nu^{(i)}=\hat{Q}_{(s_i)}^{(\Alg,\Delta)}
(\mu_{t_{i-1}}^{(\Alg,\Delta,\rm{Rec})})$
as in \eqref{eq:high-order-recombination}.
Since no recombination is performed at the last step, we have
$\mu_{T}^{(\Alg,\Delta,{\rm Rec})}=\nu^{(n)}$.
Moreover, $P^X_{T-t_n}$ is the identity, and
$\nu^{(1)}(g)=(Q_{(s_1)}g)(x_0)$ for every function $g$.
Hence the total error can be decomposed as follows by telescoping.
Here the initial value $x_0$ is fixed, and the measure
$\mu_{T}^{(\Alg, \Delta, {\rm Rec})}$ and the intermediate measures
$\mu_{t_i}^{(\Alg, \Delta, {\rm Rec})}$ and $\nu^{(i)}$ depend on it.
\begin{align}
  &\Big\lvert P^{X}_Tf (x_0)-\mu_{T}^{(\Alg,\Delta,{\rm Rec})}(f)\Big\rvert \notag\\
  &\qquad\qquad\leq\|P^{X}_T f-Q_{(s_1)}P^X_{T-t_1} f\|_{\infty} \notag\\
  &\qquad\qquad\quad +\sum_{i=1}^{n-1}
  \big\lvert
  \mu_{t_i}^{(\Alg,\Delta,{\rm Rec})}(P^X_{T-t_i} f)
  -\nu^{(i+1)}(P^X_{T-t_{i+1}}f)
  \big\rvert \notag\\
  &\qquad\qquad\quad +\sum_{i=1}^{n-1}\big\lvert\nu^{(i)}(P^X_{T-t_i} f)
  -\mu_{t_i}^{(\Alg,\Delta,{\rm Rec})}(P^X_{T-t_i} f)
  \big\rvert. \label{fest}
\end{align}
The first and second terms of the right-hand side of \eqref{fest},
which correspond to the discretisation error, can be estimated as follows.
Since $\mu_{t_i}^{(\Alg,\Delta,{\rm Rec})}$ is a probability measure and
\begin{equation*}
  \nu^{(i+1)}(g)=\mu_{t_i}^{(\Alg,\Delta,{\rm Rec})}(Q_{(s_{i+1})}g),
  \qquad
  P^X_{T-t_i}=P^X_{s_{i+1}}P^X_{T-t_{i+1}},
\end{equation*}
we have
\begin{align*}
  &\|P^{X}_T f-Q_{(s_1)}P^X_{T-t_1}f
  \|_{\infty}
  +\sum_{i=1}^{n-1}
  \big\lvert
  \mu_{t_i}^{(\Alg, \Delta, {\rm Rec})}
  (P^X_{T-t_i} f)
  -\nu^{(i+1)}(P^X_{T-t_{i+1}} f)\big\rvert\\
  &\;\;= \| P^X_{s_1}P^X_{T- t_1} f -Q_{(s_1)}P^X_{T-t_1}f
  \|_{\infty}
  +\sum_{i=1}^{n-1}\big\lvert\mu_{t_i}^{(\Alg, \Delta, {\rm Rec})}
  \big((P^X_{s_{i+1}} - Q_{(s_{i+1})})P^X_{T - t_{i+1}}f\big)
  \big\rvert\\
  &\;\;\leq
  \| (P^X_{s_1} -Q_{(s_1)})P^X_{T-t_1} f
  \|_{\infty}
  +\sum_{i=1}^{n-1}
  \| (P^X_{s_{i+1}} - Q_{(s_{i+1})}
  )P^X_{T - t_{i+1}}f
  \|_{\infty}.
\end{align*}
Then \eqref{K1}, \eqref{K3} and \eqref{K2} yield
\begin{align*}
  &\|
    P^{X}_T f - Q_{(s_1)}P^X_{T-t_1}f
    \|_{\infty}
    +\sum_{i=1}^{n-1}
    \big\lvert \mu_{t_i}^{(\Alg, \Delta, {\rm Rec})}
    (P^X_{T-t_i} f)
    - \nu^{(i+1)}(P^X_{T-t_{i+1}} f) \big\rvert\\
  &\;\;\leq
    \sum_{i=1}^{n-1} c_1\bigg({s_i}^{{(m+1)}/{2}}
    \|\nabla P^X_{T - t_{i}} f\|_{\infty}
    +\sum_{q=m+1}^{J}{s_i}^{{q}/{2}}
    \|P^X_{T-t_i}f\|_{V, q} \bigg)
    +c_3 {s_n}^{{1}/{2}}\|\nabla f\|_{\infty}\\
  &\;\;\leq \sum_{i=1}^{n-1} c_1
    \bigg({s_i}^{{(m+1)}/{2}}
    \|\nabla P^X_{T - t_{i}} f \|_{\infty}
    +c_2\sum_{q=m+1}^{J}{s_i}^{{q}/{2}}
    \frac{\|\nabla f\|_{\infty}}{(T - t_{i})^{{(q-1)}/{2}}}
    \bigg)
    +c_3 {s_n}^{{1}/{2}}\|\nabla f\|_{\infty}\\
  &\;\;\leq C_1(T) \bigg({s_n}^{1/2}
    + \sum_{i=1}^{n-1}\sum_{q=m+1}^{J}
    \frac{{s_{i}}^{q/2}}{(T-t_{i})^{(q-1)/2}}\bigg)\|\nabla f\|_{\infty}.
\end{align*}
%%%%%%%%%%%%%%%%%%%%%%%%%%%%%%%%%%%%%%%%%%%%%%%%%% 
The third term of \eqref{fest} corresponds to the recombination error.
We first estimate the derivatives of $P^X_t f$.
By Corollary~14 of Litterer and Lyons~\cite{litterer2012high},
the UH condition implies that there exists a constant $c_4>0$ such that
\begin{equation}\label{derP}
  \max_{q_1,\dots,q_{m+1}\in\{1,\dots,N\}}
  \bigg\|\frac{\partial}{\partial x_{q_1}}\cdots
  \frac{\partial}{\partial x_{q_{m+1}}}P^X_t f\bigg\|_{\infty}
  \leq c_4\, t^{-m\rev{\ell}/2}\|\nabla f\|_{\infty}
\end{equation}
for all $t\in(0,T]$ and $f\in C_b^{\infty}(\mathbb{R}^N;\mathbb{R})$.
For $t\in(0,1]$, this is the cited corollary.
For $t\in(1,T]$, it follows from the case $t=1$ applied to
$P^X_{t-1}f$, since $P^X_t=P^X_1P^X_{t-1}$ and
$\|\nabla P^X_{t-1}f\|_{\infty}\leq C\|\nabla f\|_{\infty}$
with $C$ depending only on $T$ and the vector fields.
\par
Fix $i\in\{1,\dots,n-1\}$ and $k\in\{1,\dots,l_i\}$, and write
$U=U^{(i)}_k$ and $g=P^X_{T-t_i}f$.
Choose a point $y_0\in U$, and let $q_U$ be the Taylor polynomial of $g$
of degree $m$ at $y_0$.
Since $\lvert y-y_0\rvert\leq\radius(U)$ for every $y\in U$,
Taylor's formula with integral remainder gives
\begin{equation*}
  \sup_{y\in U}\lvert g(y)-q_U(y)\rvert
  \leq\frac{N^{(m+1)/2}}{(m+1)!}\,\radius(U)^{m+1}
  \max_{q_1,\dots,q_{m+1}\in\{1,\dots,N\}}
  \bigg\|\frac{\partial}{\partial x_{q_1}}\cdots
  \frac{\partial}{\partial x_{q_{m+1}}}g\bigg\|_{\infty}.
\end{equation*}
By the construction~\eqref{eq:high-order-recombination},
the probability measures $\nu^{(i)}\vert_{U}$ and
$\tilde{\nu}^{(i)}_k$ are supported in $U$, and they agree on every
polynomial of degree at most $m$, since $G$ spans the space of such
polynomials.
Hence
\begin{equation*}
  \big\lvert\nu^{(i)}\vert_{U}(g)-\tilde{\nu}^{(i)}_k(g)\big\rvert
  =\big\lvert\nu^{(i)}\vert_{U}(g-q_U)
  -\tilde{\nu}^{(i)}_k(g-q_U)\big\rvert
  \leq 2\sup_{y\in U}\lvert g(y)-q_U(y)\rvert.
\end{equation*}
Since
$\nu^{(i)}=\sum_{k=1}^{l_i}\weight(U^{(i)}_k)\,\nu^{(i)}\vert_{U^{(i)}_k}$
and
$\mu_{t_i}^{(\Alg,\Delta,{\rm Rec})}
=\sum_{k=1}^{l_i}\weight(U^{(i)}_k)\,\tilde{\nu}^{(i)}_k$,
these estimates and \eqref{derP} with $t=T-t_i$ yield
\begin{align*}
  &\Big\lvert\nu^{(i)}(P^X_{T-t_i}f)
    -\mu_{t_i}^{(\Alg,\Delta,{\rm Rec})}(P^X_{T-t_i}f)\Big\rvert\\
  &\qquad\leq\frac{2N^{(m+1)/2}c_4}{(m+1)!}\,
    \frac{\sum_{k=1}^{l_i}\radius(U^{(i)}_k)^{m+1}
    \weight(U^{(i)}_k)}{(T-t_i)^{m\ell/2}}\,\|\nabla f\|_{\infty}.
\end{align*}
Summing over $i\in\{1,\dots,n-1\}$, we obtain the second term
of~\eqref{eq:refinedErr} with $C_2(T)=2N^{(m+1)/2}c_4/(m+1)!$.
\end{proof}
% Added in November 2026: extension of Theorem 3.1 to Lipschitz f.
\medskip
\noindent
The right-hand side of~\eqref{eq:refinedErr} depends on $f$ only through
$\|\nabla f\|_{\infty}$, and the recombined measures do not depend on
$f$.
Hence~\eqref{eq:refinedErr} extends to Lipschitz continuous functions
$f$, with $\|\nabla f\|_{\infty}$ replaced by $\|f\|_{\mathrm{Lip}}$.
Indeed, for $R>0$ let $f_R=\max\{-R,\min\{f,R\}\}$, and let
$f_{R,\varepsilon}=f_R\ast\rho_\varepsilon$ with a mollifier
$\rho_\varepsilon$ supported in the ball of radius $\varepsilon$.
Then $f_{R,\varepsilon}\in C_b^{\infty}(\mathbb{R}^N;\mathbb{R})$,
$\|\nabla f_{R,\varepsilon}\|_{\infty}\leq\|f\|_{\mathrm{Lip}}$ and
$\|f_{R,\varepsilon}-f_R\|_{\infty}\leq\varepsilon\|f\|_{\mathrm{Lip}}$.
Applying~\eqref{eq:refinedErr} to $f_{R,\varepsilon}$ and letting
$\varepsilon\to 0$, we obtain the estimate for $f_R$.
Finally, we let $R\to\infty$.
Then $P^X_Tf_R(x_0)\to P^X_Tf(x_0)$ by dominated convergence, since $f$
is of at most linear growth and $X_T$ has a finite first moment, and
$\mu_{T}^{(\Alg,\Delta,{\rm Rec})}(f_R)
=\mu_{T}^{(\Alg,\Delta,{\rm Rec})}(f)$ for large $R$, since
$\mu_{T}^{(\Alg,\Delta,{\rm Rec})}$ has finite support.
%The same applies to Corollary~\ref{locest} and Example~\ref{pel}
%below.

%%%%%%%%%%%%%%%%%%%%%%%%%%%%%%%%%%%%%%%%%%%%%%%%%%
Example~\ref{pel} below stands in the same relation to Theorem~\ref{REE}
as Example~\ref{ori:cond} does to Theorem~\ref{ori:est}.
It thus provides a refined version of the patch condition given in Example~\ref{ori:cond}.
This condition is hereafter referred to as the weighted LL patch condition,
sometimes abbreviated as the WLL patch condition.
\begin{exa}\label{pel}
  For a set of patch divisions
  $(U_k^{(i)})_{\substack{i\in\{1,\dots,\sharp \Delta-1\}\\k\in\{1,\dots,l_i\}}}$
  that satisfy  
  \begin{equation}\label{eq:refcond}
    \radius(U_k^{(i)})\le
    (T-t_i)^{{m\rev{\ell}}/(2(m+1))}\big(
    l_i\weight(U_k^{(i)}
    )\big)^{-{1}/{(m+1)}}{s_i}^{{1}/{2}},  
  \end{equation}
  we have the recombination error of order $O(n^{-(m-1)/{2}})$.
  Here $\ell$ is the integer in Theorem~\ref{REE}.
  \par
  Note that the choice of the scaling in \eqref{eq:refcond} balances the
  factor $(T-t_i)^{-m\ell/2}$, which comes from the estimate~\eqref{derP}
  of the derivatives of order $m+1$ of $P^X_{T-t_i}f$, with the spatial
  truncation error of order $\radius(U_k^{(i)})^{\,m+1}$ produced by the
  Taylor expansion of degree $m$.
  Under~\eqref{eq:refcond}, the recombination error at the $i$th step is
  bounded by $C_2(T){s_i}^{(m+1)/2}\|\nabla f\|_{\infty}$, which yields a
  global error of order $O(n^{-(m-1)/2})$ when $\max_i s_i=O(n^{-1})$.
    Moreover, this order is preserved even if the radius is rescaled by
    a constant factor.
    Indeed, multiplying the right-hand side of \eqref{eq:refcond}
    by a constant
  $\lambda>0$ does not affect the asymptotic balance between the temporal
    contribution and the spatial truncation error, and hence does
    not change the overall convergence order.
    In Section~\ref{subsec:RecErr},
    we introduce a scaling parameter $\lambda$
    in the penalty function \eqref{eq:scaledpel} and vary it to control
    the degree of recombination.
    We then numerically examine the resulting recombination error under
    different degrees of recombination.
\end{exa}

\subsection{Patch division algorithm}\label{patchdivsec}
\rev{In this subsection, we present a recursive algorithm for constructing patch divisions based on the weighted LL patch condition in \eqref{eq:refcond}.
To this end, we consider a measure $\bar{\nu}=\sum_{j=1}^{M}w_j\delta_{p_j}$ on $\mathbb{R}^N$ and set $U^{\bar{\nu}}=\operatorname{supp}(\bar{\nu})$.}
\par
\rev{Since the number of patches $l_i$ in \eqref{eq:refcond} is not known a priori, we introduce a penalty function with a patch-wise factor $\kappa \geq 1$:}
\begin{equation*}
  \Psi_i (\bar{\nu})
  =\radius(U^{\bar{\nu}}) - (T-t_i)^{m\rev{\ell}/(2(m+1))}{s_{i}}^{1/2}
  \big(\rev{\kappa}\weight(U^{\bar{\nu}})\big)^{-{1}/{(m+1)}},
\end{equation*}
\rev{where $\ell$ is as in Theorem~\ref{REE}.
The recursion described below stops either when $\Psi_i$ becomes
non-positive or when the number of support points is at most
$\card(G)$, and the patches of the latter kind are left unchanged by
the recombination (see Appendix~\ref{appendix3}).
Hence, if $\kappa_k$ denotes the factor used for the patch $U^{(i)}_k$
and $R_i$ denotes the set of indices $k$ such that $U^{(i)}_k$ has more
than $\card(G)$ points, the proof of Theorem~\ref{REE} shows that the
recombination error at the $i$th step is bounded by
$C_2(T){s_i}^{(m+1)/2}\|\nabla f\|_{\infty}\sum_{k\in R_i}\kappa_k^{-1}$.
If $\sum_{k\in R_i}\kappa_k^{-1}\leq 1$, this is the bound in
Example~\ref{pel}.
The number of patches $l_i$ is bounded above by the number $M_i$ of
nodes before recombination.
% Revised in November 2026: Theorem 3.1 holds for every patch division;
% what this choice ensures is the condition on the sum above.
%Using this upper bound as $\kappa$ provides a sufficient condition for the recombination error bound in Theorem \ref{REE}, but may lead to an unnecessarily large number of patches.\par
Using this upper bound as $\kappa$ ensures
$\sum_{k\in R_i}\kappa_k^{-1}\leq 1$, but may lead to an unnecessarily
large number of patches.}\par
\rev{Alternatively, a sufficient choice is $\kappa=2^h$, where $h$ denotes the recursion depth of $\bar{\nu}$, with $h=0$ at the root.
% Revised in November 2026: h_k is defined, and the reference to
% Theorem 3.1 is replaced by the condition on the sum above.
%Indeed, the terminal patches of a binary partition tree satisfy $\sum_k 2^{-h_k}\leq 1$, which yields the required bound on the recombination error in Theorem~\ref{REE}.
Indeed, the recursion depths $h_k$ of the terminal patches of a binary
partition tree satisfy $\sum_k 2^{-h_k}\leq 1$.
However, this choice may also result in a substantial increase in the number of patches.
% Revised in November 2026: the resulting condition is stated; it is
% stricter than the LL condition by the factor (T-t_i)^{m/(2(m+1))}.
%Furthermore, when the patch weights are approximately balanced, the factor $2^h$ offsets the decrease in the patch weights, so that the resulting condition becomes comparable to the LL condition.
Furthermore, when the patch weights are approximately balanced, that
is, $\weight(U^{\bar{\nu}})\approx 2^{-h}$, the factor $2^h$ offsets
the decrease in the patch weights, so that the resulting condition
becomes close to the unweighted condition
$\radius(U^{\bar{\nu}})\leq(T-t_i)^{m\ell/(2(m+1))}{s_i}^{1/2}$.
For the numerical experiments, we employ the simplified choice $\kappa=1$, which amounts to multiplying the right-hand side of  \eqref{eq:refcond} by $l_i^{1/(m+1)}$.
% Revised in November 2026: the bound for kappa=1, and a weaker
% statement on the numerical results (cf. Figure 4.1).
%Although this choice does not, in general, ensure the convergence order established under the condition of Example~\ref{pel}, the numerical results in Section~\ref{sec:numerical} exhibit the expected convergence behaviour.}
In this case, $\sum_{k\in R_i}\kappa_k^{-1}$ is the number $L_i$ of
elements of $R_i$, so that the convergence order established in
Example~\ref{pel} is retained if $L_i$ is bounded independently of
$n$, but it is not ensured in general.
Nevertheless, the numerical results in Section~\ref{sec:numerical}
indicate that this choice works well in practice.}

%To this end, we consider a measure
%$\bar{\nu}= \sum_{j=1}^{M}w_j\delta_{p_j}$
%on $\mathbb{R}^N$ and define a penalty function $\Psi_i$ as
%\begin{equation*}
%  \Psi_i (\bar{\nu})
%  =\radius(U^{\bar{\nu}}) - (T-t_i)^{m\rev{\ell}/(2(m+1))}{s_{i}}^{1/2}
%  \big(\weight(U^{\bar{\nu}})\big)^{-{1}/{(m+1)}},
%\end{equation*}
%where $U^{\bar{\nu}} = \support(\bar{\nu}) = \{p_j : j\in\{1,\dots,M\}\}$
%and $\ell$ is as in Theorem~\ref{REE}.
%The penalty function corresponds to the difference between the radius and
%the WLL patch condition with $l_i=1$. 
 
Given a measure $\bar{\nu}$, we perform Steps~1--3 to
divide it into two measures $\bar{\nu}_1$ and $\bar{\nu}_2$, and
continue the procedure recursively for each measure until the penalty
function $\Psi_i(\bar{\nu}_j)$ becomes non-positive: 
\begin{equation*}
  \Psi_i(\bar{\nu}_j) \le 0 \quad (j=1,2).
\end{equation*}
or the number of support points satisfies 
\begin{equation*}
\card\big(\support(\bar{\nu}_j)\big) \le \card(G),
\end{equation*}
whichever occurs first.
$\card(G)$ denotes the number of monomials in the polynomial test
function space used in the construction of the reduced measure in
Definition~\ref{redmeas}, and $\card(G)=\binom{N+m}{m}$
for the discretisation algorithm of degree $m$. 
The second condition states that, once the support size is reduced to
at most $\card(G)$, the recombination procedure cannot further
decrease the number of nodes without violating the moment conditions
required in Definition~\ref{redmeas}.
When one of the stopping criteria is met,
we define the support of $\bar{\nu}_j$ as a patch $U_k^{(i)}$, and
increment $k$ by one.

\begin{enumerate}[label=Step \arabic*.,left=0.2cm]
\item
  Calculate a weighted least-squares projection vector $v\in \mathbb{R}^{N}$.
  We first find a vector $(h_{1},\dots, h_{N-1})\in \mathbb{R}^{N-1}$
  that minimizes a least-squares error
  \[
  \sum_{i=1}^{M} w_i\bigg(p_{i,N}-\sum_{j=1}^{N-1}h_{j}p_{i, j}\bigg)^2,
  \]
  where $p_{i,j}$ denotes the $j$-th coordinate of $p_{i}$.
  That is $p_{i}=(p_{i,1},\dots, p_{i,N})$.
  We define a normalized projection vector $v$ as
  \begin{align*}
    v=&\bigg(
        \frac{\max_{i}\{p_{i,N}\}-\min_{i}\{p_{i,N}\}}
        {h_{1}\Big(\max_{i}\{p_{i,1}\}-\min_{i}\{p_{i,1}\}\Big)},
        \frac{\max_{i}\{p_{i,N}\}-\min_{i}\{p_{i,N}\}}
        {h_{2}\Big(\max_{i}\{p_{i,2}\}-\min_{i}\{p_{i,2}\}\Big)}, \notag\\
      &\qquad \dots,
        \frac{\max_{i}\{p_{i,N}\}-\min_{i}\{p_{i,N}\}}
        {h_{N-1}\Big(\max_{i}\{p_{i,N-1}\}-\min_{i}\{p_{i,N-1}\}\Big)},
        1\bigg). 
  \end{align*}
  The vector $v$ defines a projection direction that incorporates
  coordinate-wise scaling so that each coordinate contributes
  comparably to the resulting inner product.
  If the denominator vanishes,
  $h_j(\max_i\{p_{i,j}\}-\min_i\{p_{i,j}\})=0$, the corresponding component of
  $v$ is set to zero.
\item
  Sort the nodes  $\{(p_j, w_j)\}_{j\in\{1,\dots, M\}}$
  in order of the inner product with the projected vector.
  We first calculate the centre of mass $c \in \mathbb{R}^{N}$ by
  \begin{equation}\label{gravitypoint}
    c =\sum_{i=1}^{M} w_i p_i.
  \end{equation}
  We sort the nodes
  $\big((p_i, w_i)\big)_{i\in\{1,\dots, M\}}$
  in order of the inner product
  \begin{equation}\label{eq:projinnerprod}
    \langle p_i - c, v \rangle
  \end{equation}
  to $\big((p_{s(i)}, w_{s(i)})\big)_{i\in\{1,\dots,M\}}$. 

\item  Divide the nodes into two patches.
  This division is performed subject to one of the following conditions:
  (i) both patches contain the same number of nodes; or
  (ii) the total mass (sum of the weights) in each patch is approximately
  balanced.
  If (i) is used, we divide $\bar{\nu}$ into
  \begin{equation*}
    \bar{\nu}_1 = \sum_{i=1}^{\lfloor M/2 \rfloor} w_{s(i)}\delta_{p_{s(i)}},
    \qquad
    \bar{\nu}_2 = \sum_{i=\lfloor M/2 \rfloor+1}^{M} w_{s(i)}\delta_{p_{s(i)}}.
  \end{equation*}    
\end{enumerate}

%%%%%%%%%%%%%%%%%%%%%%%%%%%%%%%%%%%%%%%%%%%%%%%%%%
\subsubsection{Implementation advantages} 
% Revised in November 2026 (co-author): "that satisfies the condition in
% (3.7)" was removed, since (3.7) is not ensured for kappa=1.
The above algorithm produces a patch division, while minimising the total number of
patches as far as possible. Moreover, due to its recursive nature, the
division can be performed efficiently. This contributes to reducing
the overall computational cost of the recombination method.

\subsubsection{Computational complexity}
We briefly discuss the computational complexity of the recursive
recombination procedure. 
Let $M$ denote the number of nodes contained in a patch to be split. 
The cost of a single split is dominated by two operations: 
(i) solving the weighted least-squares problem in Step~1, and 
(ii) sorting the nodes according to the projection in Step~2. 
Solving the least-squares problem requires forming the weighted normal
equations, which can be done in a single pass over the $M$ nodes.
Since the spatial dimension $N$ is fixed, the cost of this step scales
linearly in $M$ and is of order $O(M)$. 
The subsequent sorting step incurs a computational cost of order
$O(M \log M)$. 
Therefore, the total computational cost per division is $O(M \log M)$.

The overall complexity depends on the recursion depth.
In the worst case, the subdivision continues until each patch
contains at most $\card(G)=\binom{N+m}{m}$ nodes. 
Since a recursive split reduces the number of nodes in a patch by
approximately a factor of two, 
the worst-case depth $a$ is at most on the order of
$\log_2(M/\card(G))$. 
Combining this with the per-division cost $O(M \log M)$ yields
a worst-case computational complexity of order $O(M(\log M)^2)$.
Empirically, however, we observe that the overall computational cost scales approximately as $O(M \log M)$ in typical configurations.
We report a representative behaviour of the recursion depth and wall-clock timings in Section \ref{subsec:stepwise_recomb}.

%%%%%%%%%%%%%%%%%%%%%%%%%%%%%%%%%%%%%%%%%%%%%%%%%%
\subsubsection{Example of patch division}
\begin{figure}\normalsize
  \centering
  \vspace{-3cm}
  \includegraphics[width=11.5cm]{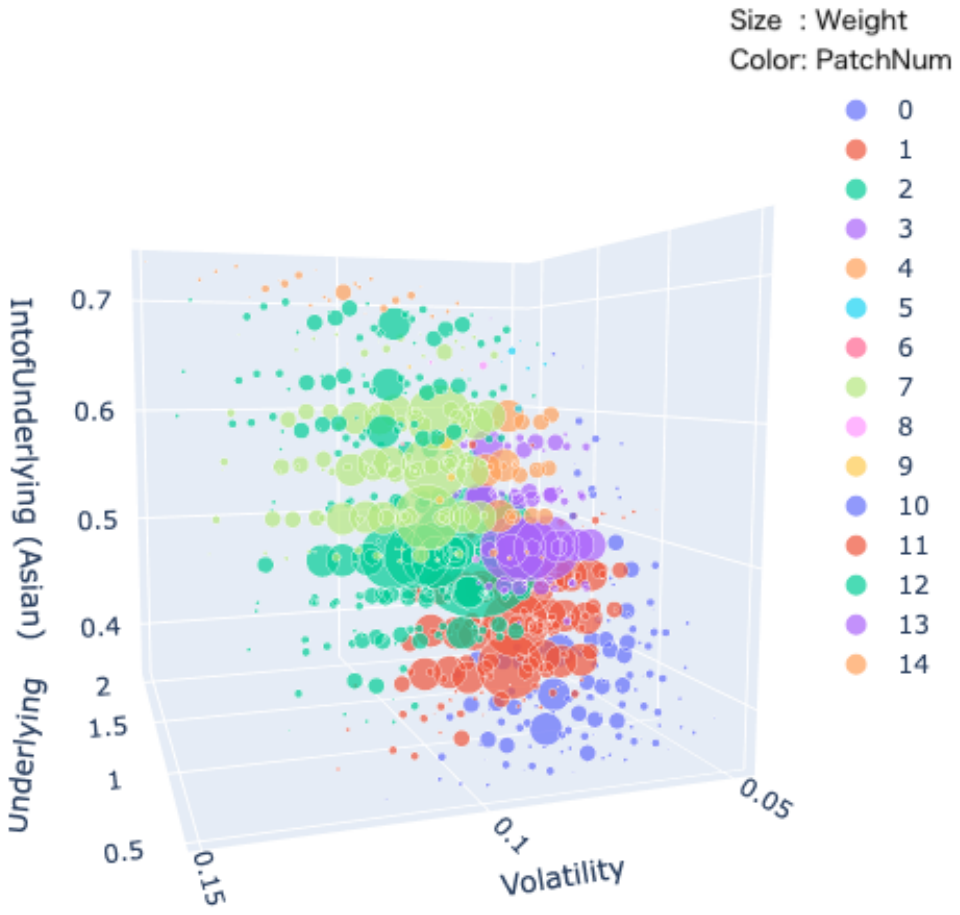}
  \vspace{-4cm}
  \caption{Example of patch division}\label{patchvis1}
\vspace{-3cm}
  \includegraphics[width=11.5cm]{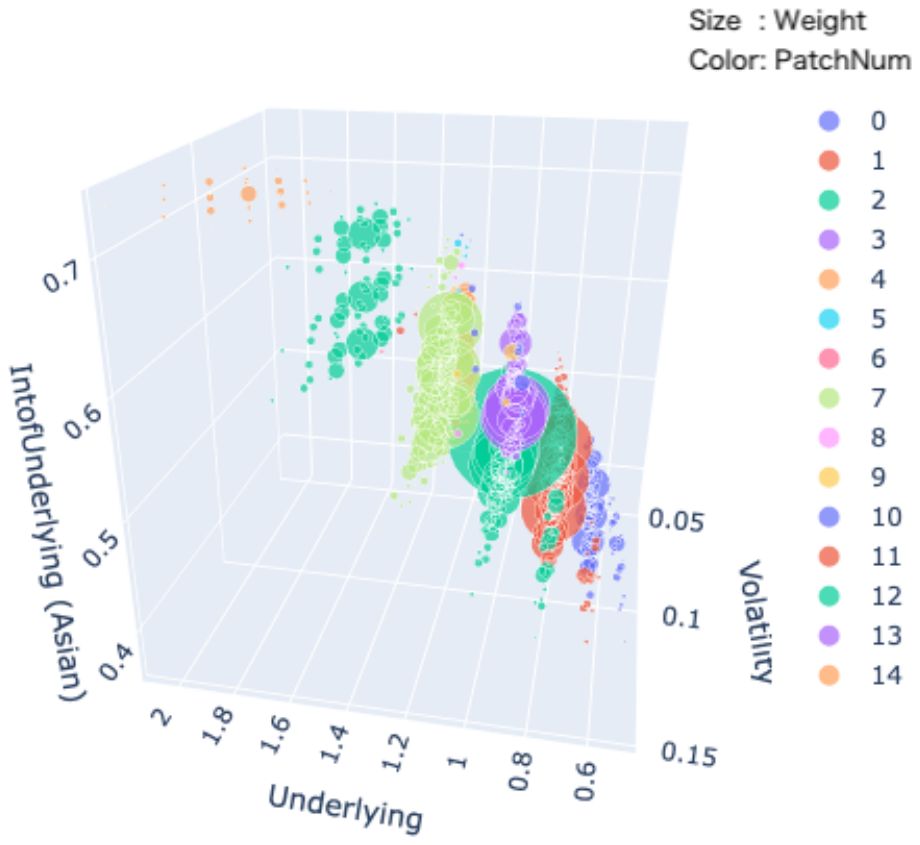}        
  \vspace{-4cm}
  \caption{Example of patch division (from a different angle)}
  \label{patchvis2}
\end{figure}

Figures~\ref{patchvis1} and~\ref{patchvis2} show an example of patch division.
This example corresponds to the pricing of an Asian option under the
Heston model.
It is discussed in Subsection~\ref{validation1}, with
the number of partitions $\sharp \Delta = 6$ and the fifth time step
$t_5 = (5/6)T$.
The size of the nodes indicates their weight and their colour
indicates the patch to which they belong.
We observe that nodes with smaller weights are represented by
smaller-sized circles and are grouped into a single large patch.

\subsection{Pseudocode of the recombination algorithm}
The overall workflow of the proposed high-order recombination method is
summarised in Algorithm~\ref{alg:main}.
At each discretisation step, the nodes are advanced by the
N--V scheme.
At every step except the last one, the nodes at time $t_{i+1}$ are then
% Revised in November 2026 (co-author): with kappa=1, in place of
% "satisfying the WLL condition", which is not ensured for kappa=1.
recursively partitioned into patches \rev{using the penalty function $\Psi_{i+1}$ of subsection~\ref{patchdivsec} with $\kappa=1$}, and
finally reduced via the recombination procedure described in
%%Appendix~1.
Appendix~\ref{appendix1}.
No recombination is performed at the last step
(see subsection~\ref{recsubsec}).
For completeness, the detailed pseudocode of the patch division and
recombination steps is provided in
%%Appendix~3.
Appendix~\ref{appendix3}.

\begin{algorithm}\normalsize
\caption{High-order recombination scheme for SDE discretisation}
\label{alg:main}
\SetKwInOut{Input}{Input}
\SetKwInOut{Output}{Output}
\SetKwFunction{OneStepFwd}{OneStepFwd}
\SetKwFunction{PatchDiv}{PatchDiv}
\SetKwFunction{Recombine}{Recombine}

\Input{Initial state $x_0$, drift and diffusion fields $V_0,\dots,V_d$, \\
       Time partition $\Delta=(t_i)_{i=0,\dots,n}$ with $\sharp \Delta=n$, payoff $f$}
\Output{Approximation of $E[f(X_{T})]$}

Initialize the support: $\nu^{(0)} = \delta_{x_0}$\;

\For{$i=0$ \KwTo $n-1$}{
  % --- Step 1: discretize
  $\nu^{(i+1)} \gets \OneStepFwd(\nu^{(i)},t_i,t_{i+1})$ \Comment*[r]{N--V discretisation}

  \If{$i+1<n$}{
    % --- Step 2: patch division
    $(U^{(i+1)}_k)_{k\in\{1,\dots,l_{i+1}\}}\gets\PatchDiv(\nu^{(i+1)})$
    \Comment*[r]{Section~\ref{patchdivsec}}

    % --- Step 3: local recombination
    $\nu^{(i+1)} \gets \Recombine(\nu^{(i+1)},\{U^{(i+1)}_k\})$
    \Comment*[r]{Appendix~\ref{appendix1}}
  }
}

\Return{$\displaystyle\sum_{x \in \support(\nu^{(n)})}
        w_{x}\, f(x)$}
\end{algorithm}
\FloatBarrier
%%%%%%%%%%%%%%%%%%%%%%%%%%%%%%%%%%%%%%%%%%%%%%%%%%
\section{Numerical results}\label{sec:numerical}
This section presents numerical experiments to examine the practical
efficiency of the high-order recombination method.
All numerical parameters used in this section are summarised in
Appendix~\ref{appendix4}
(Tables~\ref{tab:exp_setting}--\ref{tab:spiky_payoff_setting}).
% Removed in November 2026 (fifth version): the experiments did not
% recombine at the last step, in agreement with Algorithm 1, so the
% following note does not apply.
%% Added in November 2026: recombination at the last step in the
%% experiments.
%\par
%In all experiments of this section, recombination is also performed at
%the last step, unlike Algorithm~\ref{alg:main}
%(see subsection~\ref{recsubsec}).
%Theorem~\ref{REE} does not cover this step, since its proof uses that no
%recombination is performed at the last step.
%For a Lipschitz continuous function $f$, the additional error caused
%by this step satisfies
%\begin{equation*}
%  \big\lvert\nu^{(\sharp\Delta)}(f)
%  -\mu_{T}^{(\NV,\Delta,{\rm Rec})}(f)\big\rvert
%  \leq 2\|f\|_{\mathrm{Lip}}\sum_{k=1}^{l_{\sharp\Delta}}
%  \weight(U^{(\sharp\Delta)}_k)\radius(U^{(\sharp\Delta)}_k),
%\end{equation*}
%where $(U^{(\sharp\Delta)}_k)_{k\in\{1,\dots,l_{\sharp\Delta}\}}$ is the
%patch division used at this step, since the reduced measures are
%supported in the patches and preserve their total weights.
%Moreover, patches on which $f$ coincides with a polynomial of degree at
%most $m$ do not contribute to this error.
%For example, for the payoff $\max(y_3/T-K,0)$ of the Asian call option
%considered below, only patches containing points on both sides of the
%hyperplane $\{y_3=KT\}$ contribute.
%
\subsection{Problem setting}
To verify the theoretical convergence in practical problems, we consider
the pricing of an Asian call option under the Heston model.
This problem is also studied in
Ninomiya and Victoir~\cite{ninomiya2008weak},
Ninomiya and Ninomiya~\cite{ninomiya2009new}
and Shinozaki~\cite{Shinozaki2017}.
Let $Y^{(1)}_{t}$ denote the price process.
Let $Y^{(2)}_{t}$ denote the \rev{variance} process.
They are defined as follows.
\begin{align*}
  \mathrm{d}Y^{(1)}_{t} &= \mu Y^{(1)}_{t}\,\mathrm{d}t+Y^{(1)}_{t}
                           \sqrt{Y^{(2)}_{t}}\,\mathrm{d}B^{(1)}_{t}, \\
  \mathrm{d}Y^{(2)}_{t} &= \alpha ( \theta - Y^{(2)}_{t})\,\mathrm{d}t\\
                         &\quad+\beta\sqrt{Y^{(2)}_{t}}
                           (\rho \,\mathrm{d}B^{(1)}_{t}
                           +\sqrt{1-\rho^{2}}\,\mathrm{d}B^{(2)}_{t}).
\end{align*}
The initial value is given by
$(x_{1},x_{2})\in(\mathbb{R}_{>0})^{2}$.
The correlation parameter $\rho$ satisfies $-1\leq \rho \leq 1$.

The parameters $\alpha$, $\beta$, and $\theta$ are chosen to satisfy
the Feller condition
\begin{equation}\label{eq:Feller}
  2\alpha \theta - \beta^{2} > 0.
\end{equation}
This condition, due to Feller~\cite{feller1951two},
ensures the existence and uniqueness of the solution to
the CIR process $Y^{(2)}_{t}$. 
When the CIR process is discretised using the N--V method, positivity is
guaranteed under the condition
\begin{equation*}
  4\alpha \theta - \beta^{2} > 0.
\end{equation*}
This condition is weaker as a constraint on $\beta$ than
\eqref{eq:Feller}.
It can be derived from the discretisation formula presented in
Ninomiya and Victoir~\cite{ninomiya2008weak}.
%%%%%%%%%%%%%%%%%%%%%%%%%%%%%%%%%%%%%%%%%%%%%%%%%%
The payoff of the Asian call option on this asset with maturity $T$ and strike $K$ is $\max(Y^{(3)}_{T}/T-K, 0)$, where 
\[
Y^{(3)}_{t} = \int^{t}_{0}Y^{(1)}_{s}\;\mathrm{d}s.
\]
Let $Y_{t} =  {}^ t (Y^{(1)}_{t},Y^{(2)}_{t}, Y^{(3)}_{t})$.
For numerical stability, the coordinate of $Y^{(2)}$ is transformed.
The process $Y_{t}$ satisfies the following system of SDEs in
Stratonovich form:
\[
Y_{t}=x+\sum_{i=0}^{2}\int_{0}^{t}W_{i}I_N (Y_{s})\circ \mathrm{d}B^{i}_{s},
\]
where 
\begin{align*}
W_{0} I_N
\begin{pmatrix}
  y_{1} \\
  y_{2} \\
  y_{3}
\end{pmatrix}
&=
\begin{pmatrix}
  y_{1}(\mu - y_2/2)-\beta\rho y_1/4 \\
  \alpha(\theta - y_{2})-\beta^2/4 \\
  y_{1}
\end{pmatrix},
\\
W_{1} I_N
\begin{pmatrix}
  y_{1} \\
  y_{2} \\
  y_{3}
\end{pmatrix}
&=
\begin{pmatrix}
  y_{1}\sqrt{y_{2}} \\
  \beta \rho \sqrt{y_2} \\
  0
\end{pmatrix},
\\
W_{2} I_N
\begin{pmatrix}
  y_{1} \\
  y_{2} \\
  y_{3}
\end{pmatrix}
&=
\begin{pmatrix}
  0 \\
  \beta\sqrt{y_2(1-\rho^2)} \\
  0
\end{pmatrix}.
\end{align*}
The parameters are set as follows:
$T=1, K=1.05, \mu = 0.05, \alpha = 2.0, \beta = 0.1, \theta =
0.09, \rho = 0.3$, and initial value  $x=(x_{1}, x_{2}, 0)=(1.0, 0.09, 0)$.

\subsubsection{Computation of the reference value}
Under the above setup, a reference value with high accuracy is required
to assess the convergence of the discretisation error.

We adopt
\begin{equation*}
  E_{P^{X}}[f(
    Y_{T})]
  =E_{P^{X}}[\max (Y^{(3)}_{T}/T-K, 0)]
  = 0.06068740243939.
\end{equation*}
This value is computed by combining the third-order method
of Shinozaki~\cite{Shinozaki2017} with quasi-Monte Carlo.
The number of partitions is set to $\sharp \Delta = 100$ and the number
of sample points is $M_1 = 10^{9}$.

The $95\%$ confidence interval of this value is
\rev{$[0.0606855, 0.0606896]$}. %$[0.0606824, 0.0606874]$.
It is estimated from $R=20$ independent randomly shifted QMC runs
\cite{rQMCori,rQMCsurvey} with \rev{$M_2 = 5 \times 10^{7}$}.
The across--replication variance $\widehat{\sigma}^2_{M_2}$ is estimated
and extrapolated to $M_1 = 10^{9}$ via
\[
\widehat{\sigma}^2_{M_1}
= \widehat{\sigma}^2_{M_2}\, (M_2/M_1).
\]
%%%%%%%%%%%%%%%%%%%%%%%%%%%%%%%%%%%%%%%%%%%%%%%%%%
\subsection{Convergence of discretisation error}\label{validation1}
\begin{figure}\normalsize
  \begin{center}
    \includegraphics [width=11cm,angle=0,origin=c]{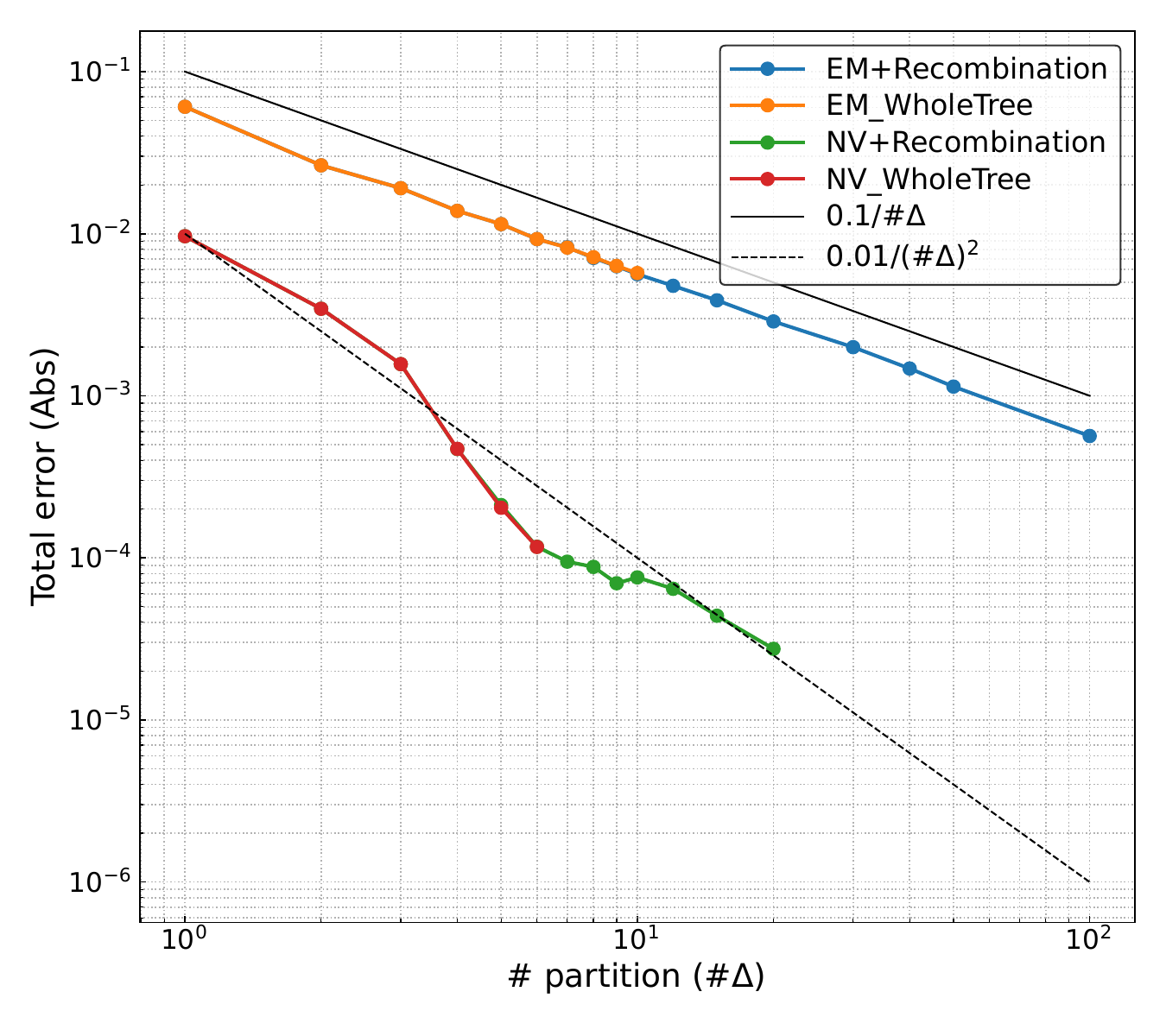}%{HesAsi_discre_reprod.pdf}  
    \caption{Total error for the Heston model}
    \label{HesDiscre}
    
  \end{center}
\end{figure}

Figure~\ref{HesDiscre} shows the convergence of the total error
\begin{equation*}
  \lvert
  E_{\mu^{(\NV, \Delta, {\rm Rec})}}
  [f(Y^{(\sharp \Delta)}_{T})]
  -
  E_{P^{X}}
  [f(Y_{T})]
  \rvert.
\end{equation*}
This quantity represents the difference between the expectations under
the recombined measure $\mu_{T}^{(\NV, \Delta, {\rm Rec})}$ and the
original measure $P^{X}_{T}$.
The error is plotted with respect to the number of partitions
$\sharp \Delta$.

The time partition $\Delta$ is chosen as the uniform partition.
It is given by $s_i = 1/\sharp \Delta$.
% Revised in November 2026: the patches are constructed by the algorithm
% of subsection 3.2 with kappa=1 (Example 3.2 has no factor kappa).
%The set of radii is chosen as in Example~\ref{pel} with $\ell=3, \kappa=1$
\rev{The patches are constructed by the algorithm of
subsection~\ref{patchdivsec} with $\ell=3$ and $\kappa=1$
(see Remark~\ref{hesrem}). 
For E--M and N--V, we set $m=3$ and $m=5$, respectively, so that $\mathrm{card}(G)=20$ and $56$.} 
See Remark~\ref{rem:optimalpartition} below for the reasons for adopting
the uniform partition.

\begin{figure}\normalsize
    \begin{center}
    \includegraphics [width=11cm,angle=0]{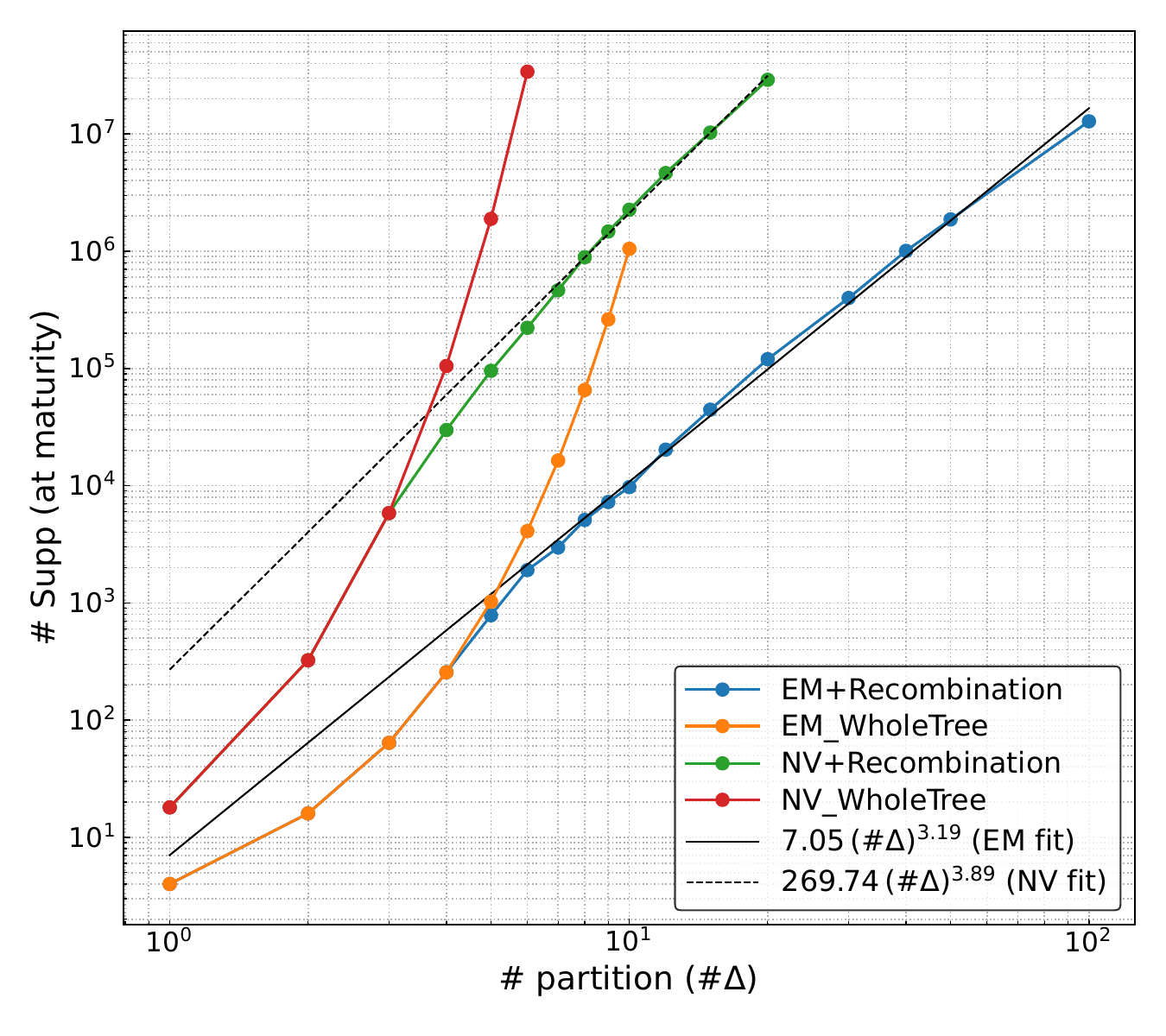}%{HesAsi_numpatch_fit_auto.pdf}
    \caption{Growth of the support size for the Heston model}\label{HesSupp}
  \end{center}
\end{figure}

The errors of the E--M and N--V methods with and without recombination
are plotted.
Figure~\ref{HesSupp} shows the growth of the support size.
As shown in Figure~\ref{HesDiscre}, the theoretical convergence order of
the error is observed for each method, both with and without
recombination.
Without recombination, the whole-tree computation is executable only
up to around $\sharp \Delta = 8$.
This limitation is due to the computational complexity.

With recombination, larger values of $\sharp \Delta$ can be treated and
the theoretical order of convergence is observed.
Figure~\ref{HesSupp} shows that recombination suppresses the exponential
growth of the support size and yields polynomial growth.
Specifically, the fitted slopes in Figure~\ref{HesSupp} indicate that the
support size after recombination grows as
$O(\sharp \Delta^{\alpha})$.
The estimated values, \rev{obtained by fitting the data points with $\sharp\Delta \geq 7$ where recombination becomes effective}, are $\alpha \approx 3.19$ for E--M and
$\alpha \approx 3.89$ for N--V.
This behaviour is in sharp contrast to the exponential growth observed
in the corresponding whole-tree computations.

\begin{figure}\normalsize
  \begin{center}
    \includegraphics [width=11cm,angle=0,origin=c]{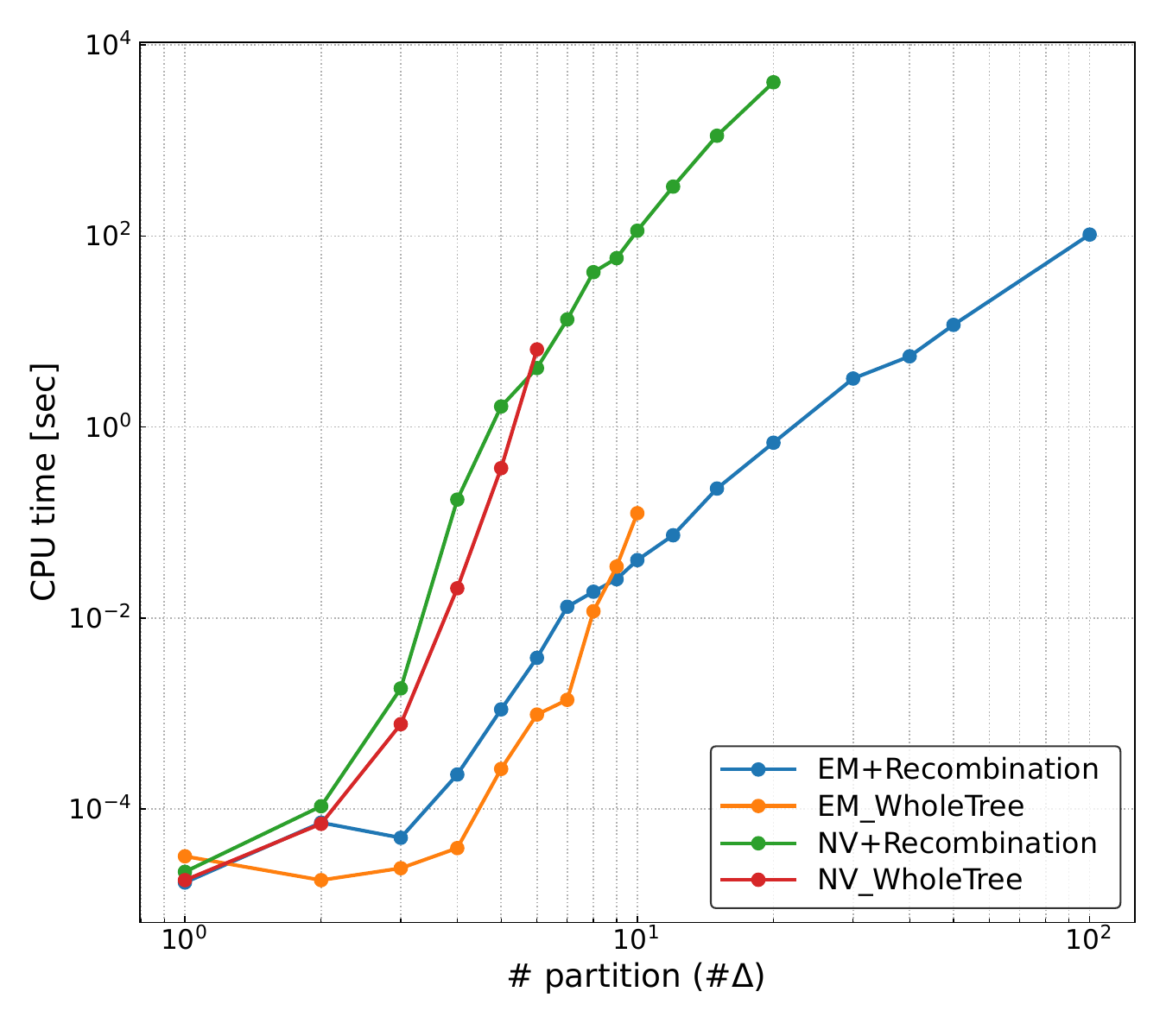}%{HesAsi_CPU.pdf}    
    \caption{CPU time for the Heston model}\label{HesCPU}
  \end{center}
\end{figure}

Figure~\ref{HesCPU} shows that the CPU time exhibits a similar growth
pattern.
This reflects that the runtime is essentially determined by the size of
the recombined support.
%% %%%%%%%%%%%%%%%%%%%%%%%%%%
\begin{rem}\label{hesrem}
  In the Heston model the positivity of the discretised process obtained
  by the N--V method ensures that
  \[
  \support (\nu^{(i)}) \subset (\mathbb{R}_{>0})^3 =
  \{(y_1, y_2, y_3)\;:\; y_1, y_2, y_3 > 0 \},
  \]
  where $\nu^{(i)}$ is an intermediate measure defined in
  \eqref{eq:high-order-recombination}.
  This implies that for any patch division
  $(U_{k}^{(i)})_{k\in\{1,\dots,l_i\}}$ of $\support (\nu^{(i)})$
  there exists $\varepsilon>0$ such that
  \[
  \bigcup_{k=1,\dots,l_i} U_{k}^{(i)} \subset W_{\varepsilon}=
  \{(y_1, y_2, y_3)\;:\; y_1, y_2, y_3 > \varepsilon \}.
  \]
  % Revised in November 2026: the spanning property requires |rho| < 1.
  Moreover, \rev{for $\lvert\rho\rvert<1$,} the set of vectors
  $\{W_1(x), W_2(x), [W_0, W_1](x)\}$ spans the linear space
  $\mathbb{R}^3$ for each $x \in (\mathbb{R}_{>0})^3$ in the Heston model.
  \rev{Thus, the H\"{o}rmander bracket-generating condition holds at every point of $(\mathbb{R}_{>0})^3$.}
  Since $W_1$, $W_2$ and their brackets have no component in the
  direction of $y_3$, only brackets involving $W_0$ can have such a
  component.
  These brackets have order at least $3$, and $[W_0,W_1]$ has order $3$
  and the $y_3$-component $-y_1\sqrt{y_2}\neq 0$.
  Hence this condition holds with $\ell=3$ and with no smaller $\ell$.
  % Revised in November 2026: uniformity on W_epsilon.
  \rev{On $W_{\varepsilon}$, the infimum in Definition~\ref{UH} with
  $\ell=3$, taken over $x\in W_{\varepsilon}$ instead of
  $x\in\mathbb{R}^N$, is positive, but it tends to zero as
  $\varepsilon\to 0$, since the vector fields appearing in this infimum
  have $y_3$-components $0$ or $\pm y_1\sqrt{y_2}$.}
  Accordingly, we use $\ell=3$ in the patch conditions throughout this
  section, namely in~\eqref{eq:refcond}, in~\eqref{eq:scaledpel} and in
  the condition of Example~\ref{ori:cond}.
 \rev{However, this local bracket-generating property does not directly imply the global UH condition assumed in Theorem~\ref{REE}. Moreover, the coefficients of the Heston model are unbounded and are not smooth at $y_2=0$.}
%  However, Corollary~\ref{locest} also assumes the UFG condition, which
%  we do not verify for the Heston model, whose coefficients are unbounded
%  and are not smooth at $y_2=0$.
%  In addition, the value of $\varepsilon$ depends on the discretisation and
%  recombination algorithm.
%  It cannot be determined before executing the algorithm.
  \rev{In addition, the positive lower bound on $W_{\varepsilon}$ depends
  on $\varepsilon$, whose value depends on the discretisation and
  recombination algorithm and cannot be determined before executing the
  algorithm.}
%  Therefore, to justify the proposed algorithm for the Heston model, it is
%  necessary to establish Theorem~\ref{REE} under a weaker condition than
%  the UH condition.
%  One possible candidate is the UFG condition.
  \rev{Therefore, to justify the proposed algorithm for the Heston model,
  it is necessary to extend Theorem~\ref{REE} to vector fields that are
  smooth only on a domain such as $(\mathbb{R}_{>0})^3$ and satisfy a
  weaker, possibly local, condition than the UH condition, such as a
  localised version of the UFG condition.}
  This remains a subject for future work.
\end{rem}

\begin{rem}
  \label{rem:optimalpartition}
  Throughout these numerical experiments, a uniform time partition is
  used.
  It is given by $s_i = T/\sharp \Delta$.
  For the original continuous versions of the N--V and N--N methods, this
  choice is proved to be optimal in Kusuoka~\cite{kusuoka2013gaussian}.
  This optimality result does not apply to the corresponding discrete
  versions, as noted in Remark~\ref{regpart}.
  However, the numerical experiments indicate that the uniform partition
  also performs well for the discretised N--V method.
  Similar observations are reported in
  Ninomiya and Shinozaki~\cite{NinomiyaShinozaki2019}.
\end{rem}
%%%%%%%%%%%%%%%%%%%%%%%%%%%%%%%%%%%%%%%%%%%%%%%%%%
\begin{rem}
In Step~3 of the patch division algorithm, the nodes can be divided
either by equal number or by approximately equal total weight.
Both approaches are tested in the numerical experiments.
No significant difference in recombination error is observed.
Therefore, the equal-number node division criterion, denoted by
option~(i) in Step~3 of the patch division algorithm
(subsection~\ref{patchdivsec}),
is employed for all numerical results reported in this paper.
\end{rem}
%%%%%%%%%%%%%%%%%%%%%%%%%%%%%%%%%%%%%%%%%%%%%%%%%%

%%%%%%%%%%%%%%%%%%%%%%%%%%%%%%%%%%%%%%%%%%%%%%%%%%
\subsection{Step-wise behaviour of the recombination algorithm}
\label{subsec:stepwise_recomb}
\begin{figure}\normalsize
  \begin{center}
    \includegraphics[width=11cm]{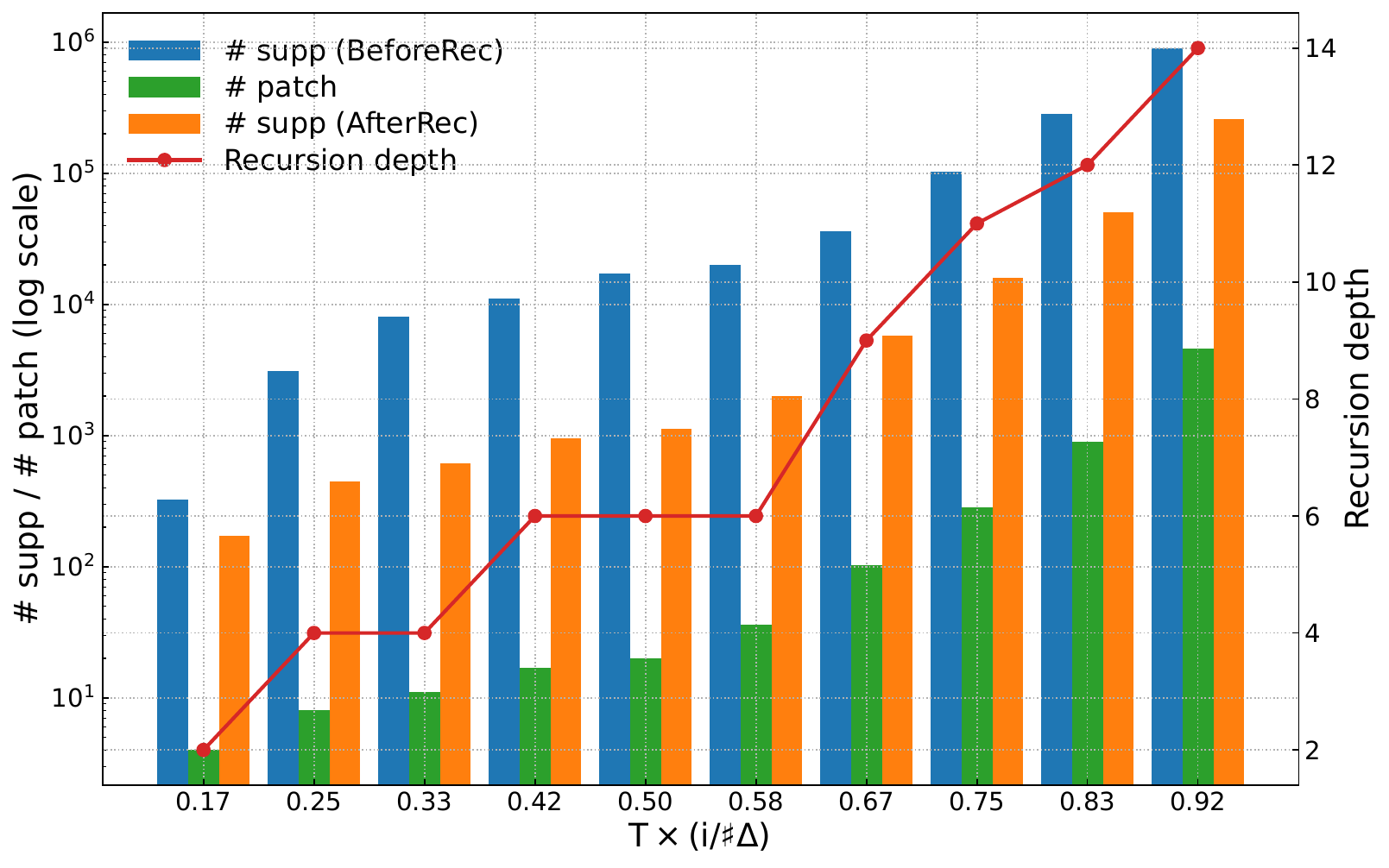}%{recomb_step_summary.pdf}
    \caption{Step-wise behaviour of the recombination algorithm.
      The left axis, shown on a logarithmic scale, indicates the number of
      nodes before recombination and the number of patches after
      recombination.
      The right axis shows the maximum recursion depth.}
    \label{fig:step_recomb}
  \end{center}
\end{figure}

To illustrate how the recombination operates at each time step, the
number of time partitions is fixed to $\sharp \Delta = 12$.
% Revised in November 2026: kappa=1, as in subsection 4.2.
The patches are constructed by the patch division algorithm of
subsection~\ref{patchdivsec} with $\ell=3$ \rev{and $\kappa=1$}.
For each step $i = \rev{2},\dots,\sharp \Delta\rev{-1}$, the following quantities are
recorded:
(i) the number of nodes before recombination,
(ii) the number of patches after recombination, and
(iii) the maximum recursion depth required in the patch division
procedure.
Figure~\ref{fig:step_recomb} summarises these quantities.

The horizontal axis represents the time step
$T \times ( i / \sharp \Delta )$.
The left vertical axis, shown on a logarithmic scale, indicates the
number of nodes before recombination and the number of patches produced
after recombination.
The right vertical axis shows the maximum recursion depth.

The figure shows that the number of nodes before recombination increases
rapidly as time advances.
By contrast, the number of patches after recombination remains
moderate.
The depth of recursive subdivision also remains limited throughout the
simulation.
This step-wise observation indicates that the proposed recombination
algorithm effectively controls the local computational complexity at
each time step.
%%%%%%%%%%%%%%%%%%%%%%%%%%%%%%%%%%%%%%%%%%%%%%%%%%
%%%%%%%%%%%%%%%%%%%%%%%%%%%%%%%%%%%%%%%%%%%%%%%%%%%%
%%%%%%%%%%%%%%%%%%%%%%%%%%%%%%%%%%%%%%%%%%%%%%%%%%
\subsection{Effect of the patch division algorithm on recombination error}
\label{subsec:RecErr}
This subsection examines the effect of the patch division algorithm on
the recombination error.
The recursive patch division algorithms are compared with several
simpler alternatives, including concentric circle division (CC) and
random division (Rand).
The theoretical patch division condition is not necessarily satisfied by
the CC and Rand methods.
The comparison also includes recursive patch division algorithms based
on the Litterer--Lyons condition and its weighted variant.
These methods are denoted by {\rm Recursive\_LL} and
{\rm Recursive\_WLL}, respectively.

In the recursive patch division algorithms, the degree of recombination
can be controlled by changing the number of patches.
% Revised in November 2026: this holds for patch divisions satisfying
% (3.7); for kappa=1, see the sentence after (4.2).
%This adjustment does not affect the overall convergence order, as
%discussed in Example~\ref{pel}.
\rev{For patch divisions satisfying~\eqref{eq:refcond}, this adjustment
does not affect the overall convergence order, as discussed in
Example~\ref{pel}.}
To this end, a scaling parameter $\lambda>0$ is introduced into the
penalty function.
It is defined by
\begin{equation}\label{eq:scaledpel}
  \Psi^{(\lambda)}_i(U)
  = \radius(U)
  - \lambda (T - t_i)^{{m\rev{\ell}}/({2(m+1)})}
    {s_i}^{{1}/{2}}
    \big(\weight(U)\big)^{-{1}/{(m+1)}}.
\end{equation}
As in subsections~\ref{validation1} and~\ref{subsec:stepwise_recomb},
we set $\ell=3$ in~\eqref{eq:scaledpel} (see Remark~\ref{hesrem}).
% Added in November 2026: (4.2) corresponds to kappa=1, and the bound in
% subsection 3.2 acquires the factor lambda^{m+1}.
\rev{Thus \eqref{eq:scaledpel} is the penalty function of
subsection~\ref{patchdivsec} with $\kappa=1$ and with its second term
multiplied by $\lambda$, and the bound on the recombination error given
in subsection~\ref{patchdivsec} is multiplied by $\lambda^{m+1}$.}
The condition of Example~\ref{ori:cond}, on which {\rm Recursive\_LL}
is based, is also taken with $\ell=3$.
% Added in November 2026: the condition of Example 2.14 under the
% uniform partition.
Note that Example~\ref{ori:cond} assumes the non-uniform partition given
there, whereas the uniform partition is used in the experiments
(Remark~\ref{rem:optimalpartition}).
Under the uniform partition, the condition of Example~\ref{ori:cond}
makes the $i$th recombination term in Theorem~\ref{ori:est} equal to
${s_i}^{(m+1)/2}(T-t_i)^{-m/2}$, which is also the $i$th term in
Theorem~\ref{KLNVAppThm}~(b), and the resulting bound on the
recombination error is only of order $O((\sharp\Delta)^{-1/2})$.
By contrast, under~\eqref{eq:refcond} the recombination error is of order
$O((\sharp\Delta)^{-(m-1)/2})$ also for the uniform partition
(Example~\ref{pel}).
% Added in November 2026: Recursive_WLL uses kappa=1.
\rev{Note, however, that {\rm Recursive\_WLL} uses $\kappa=1$, for which
\eqref{eq:refcond} is not ensured and only the bound in
subsection~\ref{patchdivsec} is available.}
Since the N--V method is a Gaussian K-scheme in the sense of
Kusuoka~\cite{kusuoka2013gaussian}, the order $(m-1)/2$ can be expected
also for {\rm Recursive\_LL} under the uniform partition, but a proof
for the discrete N--V method with recombination has not yet been given.

For each patch division algorithm
\[
\text{PatchDiv} \in
\{\text{CC}, \text{Rand},
\text{Recursive\_LLcond}^{(\lambda)},
\text{Recursive\_WLLcond}^{(\lambda)}\},
\]
the corresponding recombined measure at time $T$ is denoted by
\[
\mu_T^{(\NV, \Delta, {\rm Rec}, \text{PatchDiv})}.
\]
The notation
${\rm Recursive\_LLcond}^{(\lambda)}$ and
${\rm Recursive\_WLLcond}^{(\lambda)}$
emphasises the dependence of the recursive algorithms on the scaling
parameter $\lambda$.

For a fixed time partition $\sharp \Delta$, the recombination error is
defined by
\begin{equation}\label{eq:recomb-error}
  \lvert
  \mu_{T}^{(\NV, \Delta, {\rm Rec}, {\rm PatchDiv})}(f)
  -
  \mu_{T}^{(\NV, \Delta)}(f)
  \rvert.
\end{equation}
The second term corresponds to the unrecombined measure
$\mu_{T}^{(\NV, \Delta)}$.
It is defined as the full tree-based measure \eqref{tbcalc} without
applying recombination.

\begin{figure}\normalsize
  \begin{center}
    \includegraphics[width=8.5cm]{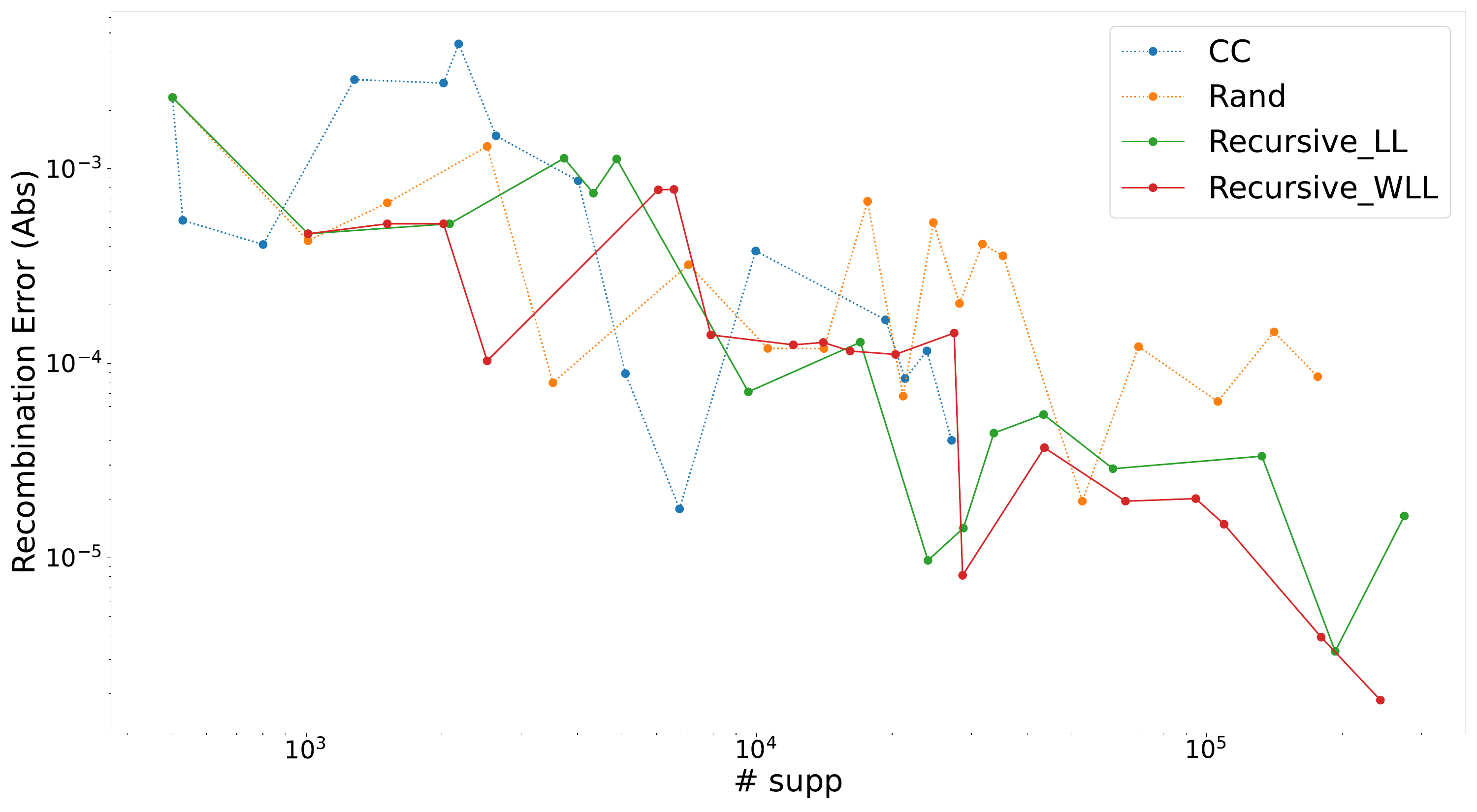}%{RecErr_Asian.pdf}
    \caption{Recombination error for the Heston model}
    \label{recerr}
  \end{center}
\end{figure}

Figure~\ref{recerr} shows the recombination errors for pricing the Asian
option under the Heston model.
The numerical setting is the same as in the previous subsection.
The number of partitions is fixed to $\sharp \Delta = 8$.
Among the methods considered, {\rm Recursive\_WLL} achieves the smallest
recombination error, followed by {\rm Recursive\_LL}.
The random division method yields moderate accuracy.
The CC method results in the largest error among the four.

%% \begin{figure}[H]

To further investigate situations in which the patch conditions have a
greater influence on the recombination error, additional experiments
are conducted using spiky payoff functions.
Define
\[
(\tilde{X}_{T}, \tilde{Y}_{T})
=
\bigg(
\frac{Y^{(3)}_{T}}{Y^{(1)}_{0}},
\frac{Y^{(2)}_{T}}{Y^{(2)}_{0}}
\bigg).
\]
For parameters $(x_c,y_c)$, $r$, and $h$, define
\[
g^{(x_c, y_c, r, h)}(x, y)
=
\max\bigg\{0,\; h - \frac{h}{r}\sqrt{(x-x_c)^2+(y-y_c)^2}\bigg\}.
\]
The parameters $h$, $r$, and $(x_c,y_c)$ are the height, the radius,
and the centre of the cone.
We set $h = 10 r^{-2}$ so that the volume of the cone remains constant.
\par
Eighty cone payoff functions are considered.
The parameters are chosen as follows:
\begin{itemize}\setlength{\itemindent}{-3mm}
\item Radius $r$: 5.0, 4.0, 3.0, 2.0, 1.5, 1.2, 1.0, 0.9
\item Centre $(x_c, y_c)$: (0,0), (1,0), (0,1), (1,1), (2,0), (0,2),
      (2,2), (3,0), (0,3), (3,3)
\end{itemize}

\begin{figure}\normalsize
  \begin{minipage}[b]{0.48\linewidth}
    \includegraphics[width=\linewidth]{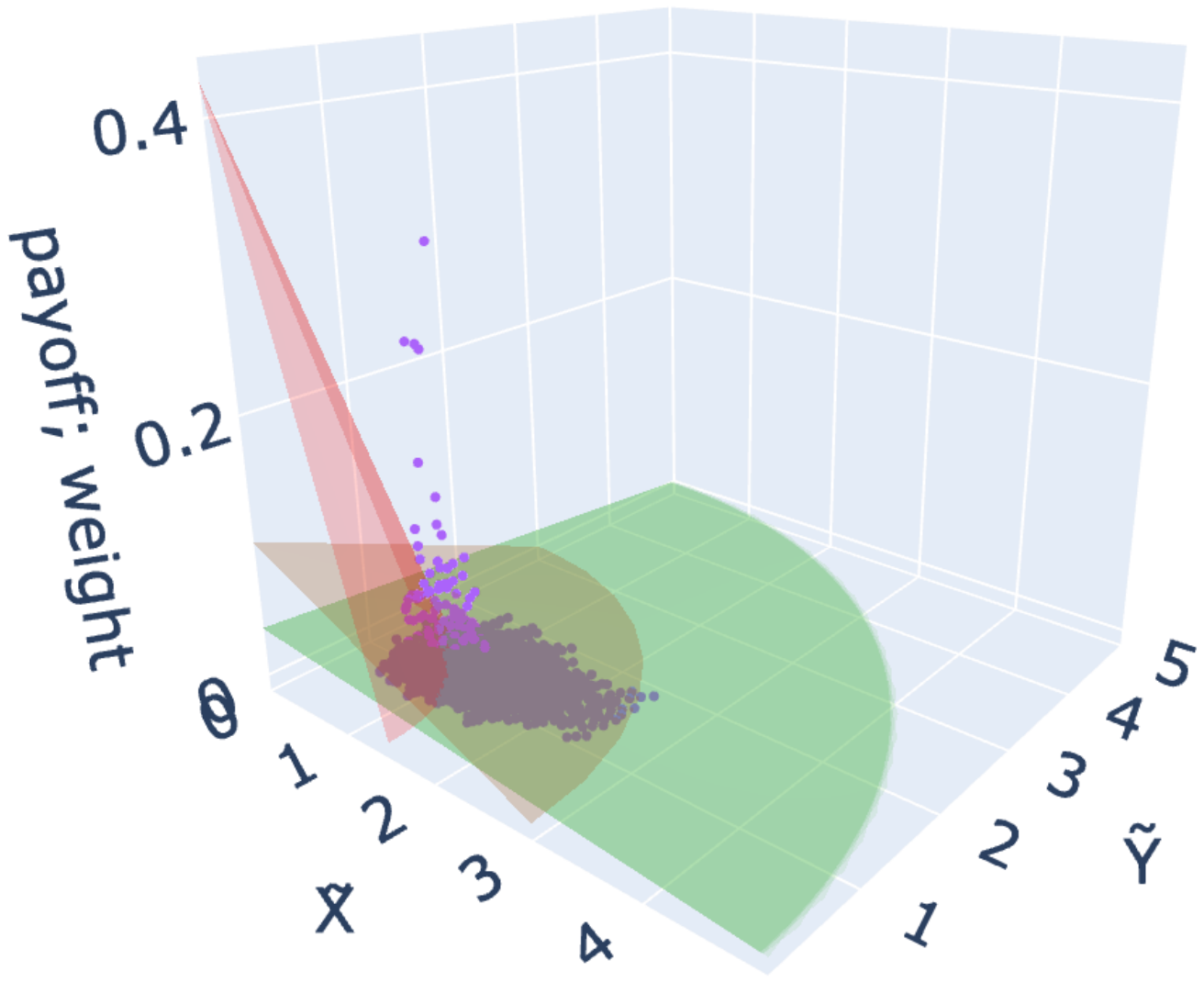}    
  \end{minipage}
  \begin{minipage}[b]{0.48\linewidth}
    \centering
    \includegraphics[width=\linewidth]{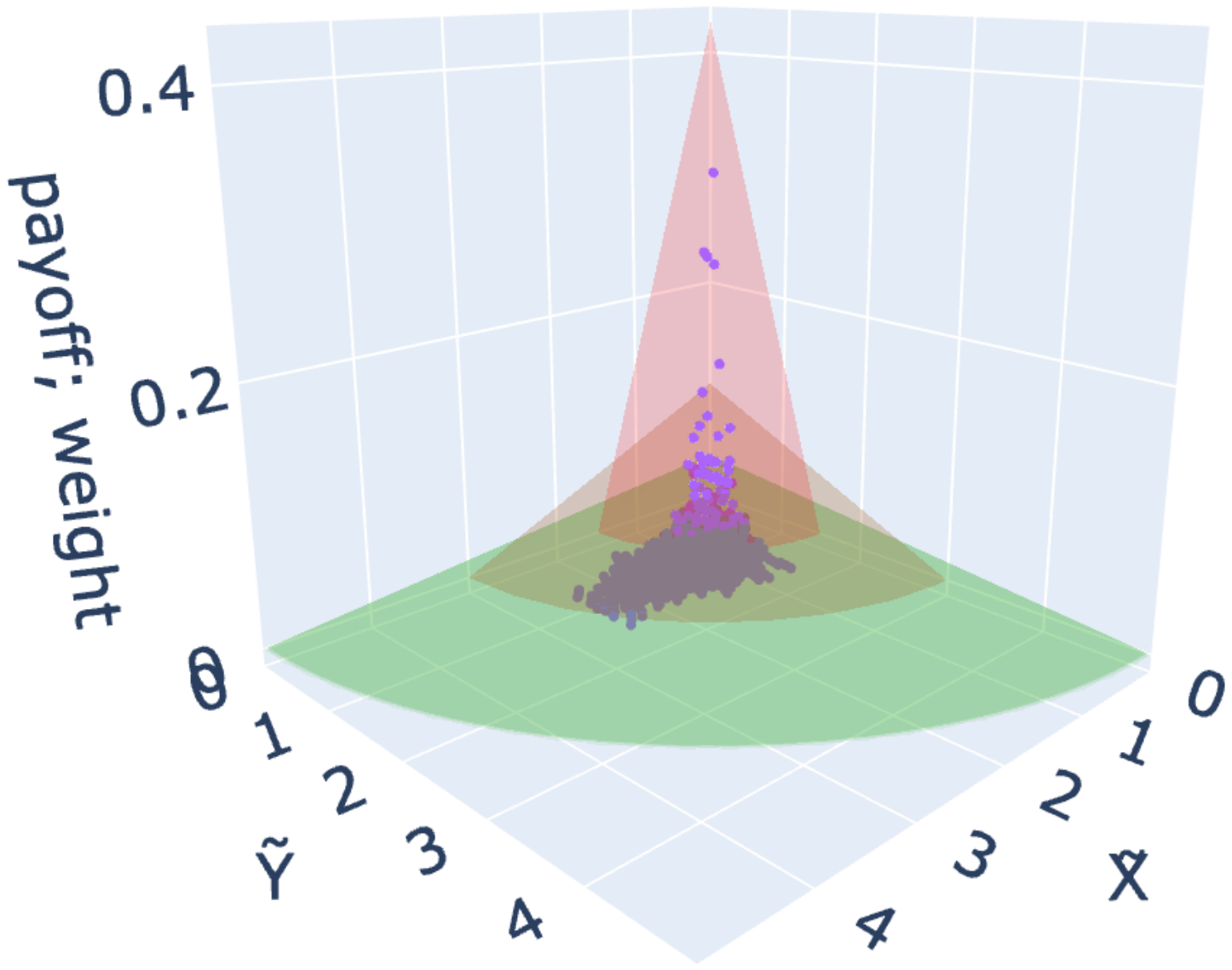}
  \end{minipage}\caption{Spiky payoff}\label{fig:cone}
\end{figure}

Figure~\ref{fig:cone} shows examples of these payoff functions together
with the nodes at maturity.
This demonstrates that the above parameter choices cover a wide range
of configurations.
For $\sharp\Delta=6$ and $8$, we estimate the slope of the log-log regression of the recombination error
\[
  \left|\mu_T^{(\mathrm{NV},\Delta,\mathrm{Rec},\mathrm{PatchDiv})}(g)-\mu_T^{(\mathrm{NV},\Delta)}(g)\right|
\]
against the support size at maturity
\[
  \mathrm{card}\left(\mathrm{supp}\left(\mu_T^{(\mathrm{NV},\Delta,\mathrm{Rec},\mathrm{PatchDiv})}\right)\right)
\]
for each payoff and each patch division method.
%For $\sharp \Delta = 6$ and $8$, the normalised recombination error is
%computed.
%It is defined as the recombination error divided by the support size at
%maturity.
%That is,
%\[
%\frac{
%\lvert
%\mu_{T}^{(\NV, \Delta, {\rm Rec}, {\rm PatchDiv})}(g)
%-
%\mu_{T}^{(\NV, \Delta)}(g)
%\rvert
%}{
%\card\!\big(
%\support\!(
%\mu_{T}^{(\NV, \Delta, {\rm Rec}, {\rm PatchDiv})}
%)
%\big)
%}.
%\]

\begin{figure}\normalsize
  \includegraphics [width=11.5cm]{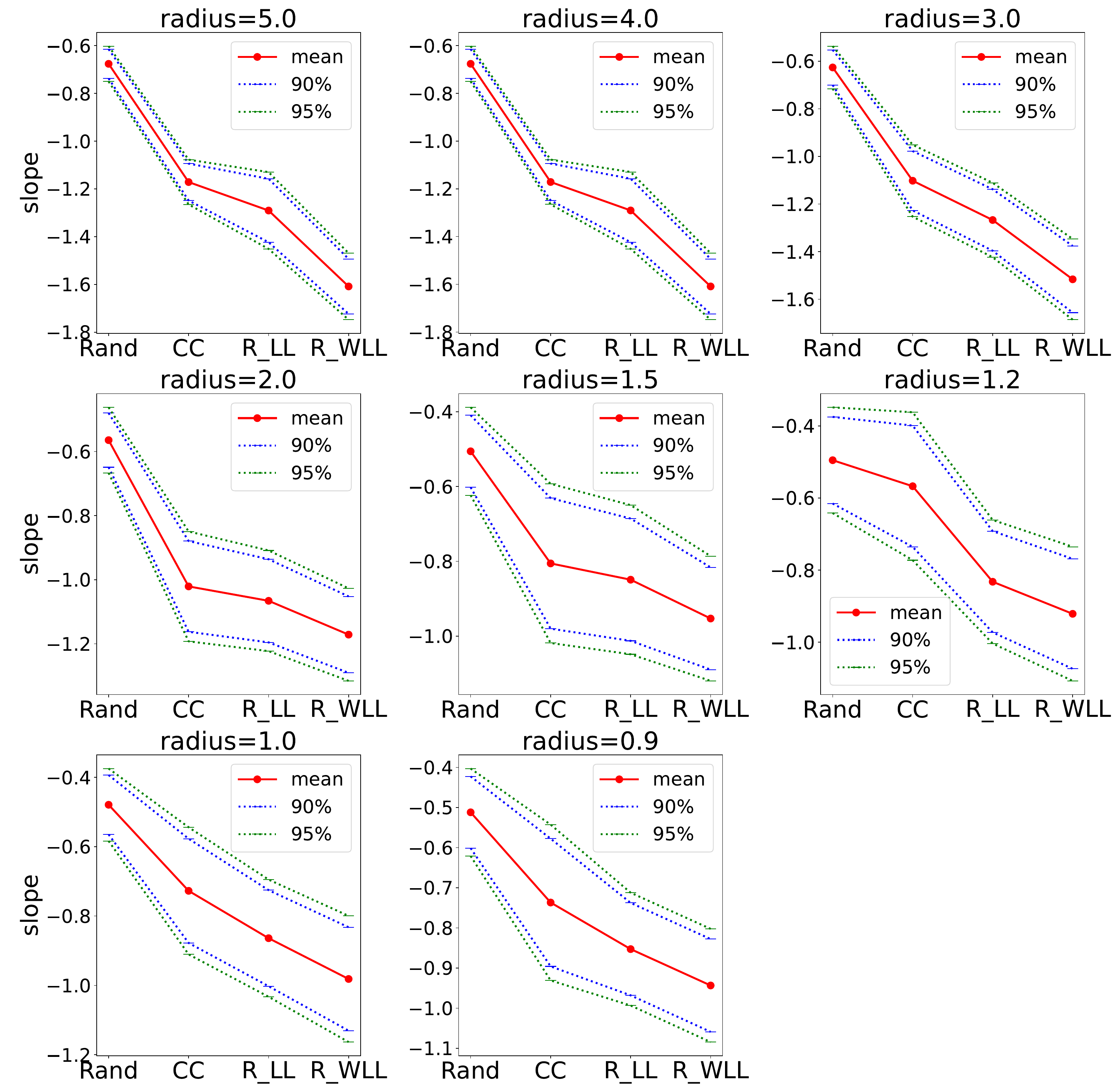}%{average_confinterval_rad.pdf}
  \caption{\rev{Mean slopes of the log-log regressions and their confidence intervals}}\label{slope1}    
\end{figure}

\begin{figure}\normalsize
  \includegraphics [width=11.5cm]{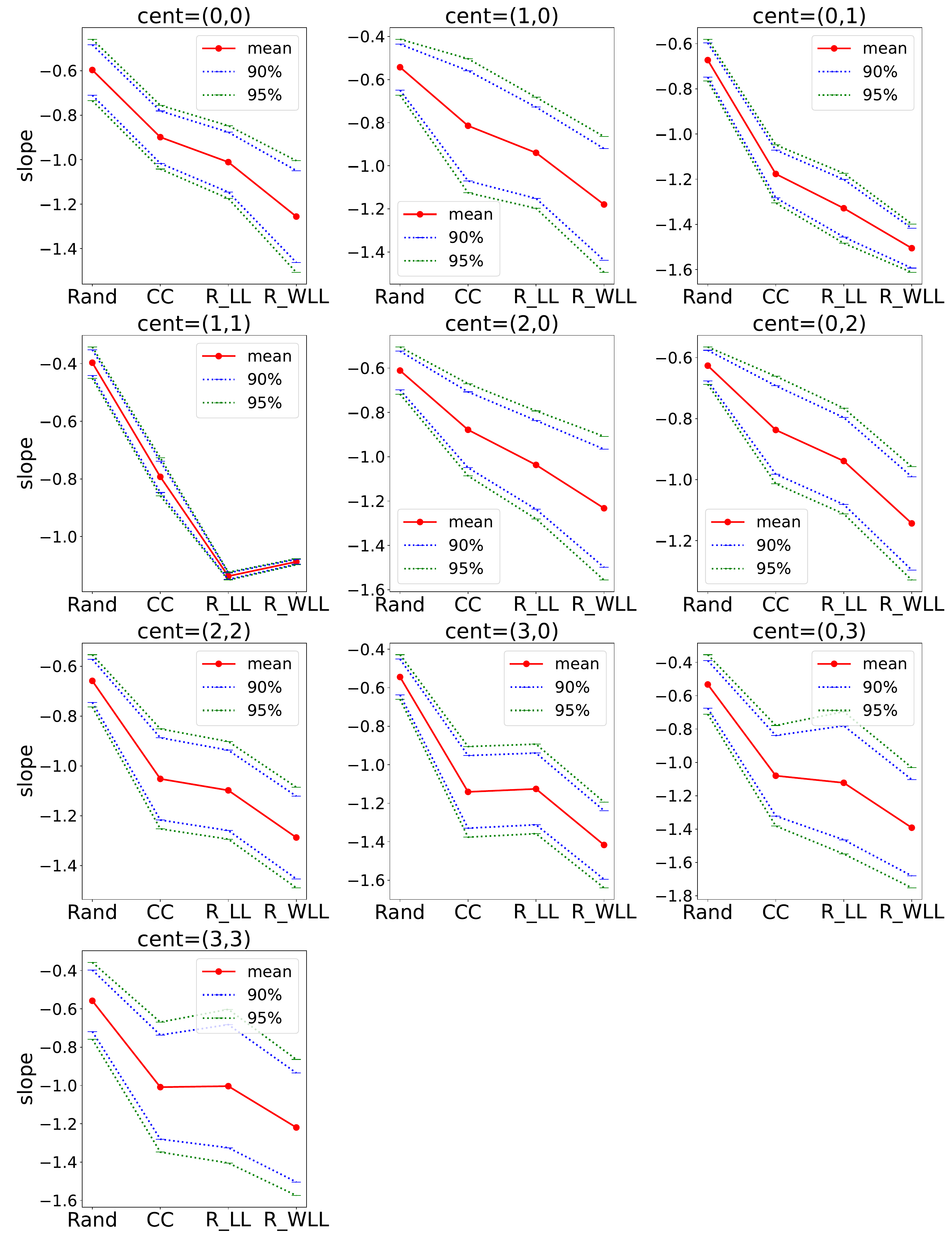}%{average_confinterval_cent.pdf}
  \caption{\rev{Mean slopes of the log-log regressions and their confidence intervals}}\label{slope2}
\end{figure}

Figures~\ref{slope1} and~\ref{slope2} summarise \rev{the mean slopes} together with their confidence
intervals. 
The results are obtained from \rev{136} test cases.
These cases consist of \rev{68} cone payoff functions and two choices of
$\sharp \Delta$.
\rev{The combinations with $(x_c,y_c)\in\{(3,0),(0,3),(3,3)\}$ and $r\leq 1.5$ are excluded, since the supports of the corresponding payoff functions lie largely outside the support of the distribution at maturity.}
Figure~\ref{slope1} corresponds to fixed radii.
Figure~\ref{slope2} corresponds to fixed centre points.
The results indicate that the proposed patch division algorithm performs
well in most cases.

\FloatBarrier
\begin{acknowledgements}
One of the authors gratefully acknowledges the support and hospitality
of the Research Institute for Mathematics (Forschungsinstitut für
Mathematik, FIM), ETH Zurich, during a visiting appointment.
\end{acknowledgements}

\appendix\normalsize
\section{Construction of a reduced measure}\label{appendix1}
Let $\nu=\sum_{j=1}^{M}w_j\delta_{p_j}$ be an induced probability
measure on a patch $U$, defined in \eqref{eq:patch-measure}. 
$M$ denotes the number of nodes included in the patch $U$.
Let $G=\{g_1, g_2,\dots, g_{\card(G)}\}$ be a set of test functions. 
Hereafter, $x^{(j)}$ denotes the $j$th coordinate of $x$.
Thus $x = (x^{(1)},\dots, x^{(N)})$.

\begin{enumerate}
\item[(1)] Define the initial test-function matrix $A_0$
  and the initial weight vector $\beta_0$: 
  We set the $\card(G)\times M$ matrix $A_0$ by evaluating each 
  $g_i$ at each node $p_j$
  \begin{equation*}
    A_0 = \big(g_i(p_j)\big)_{i \in \{1,\dots,\card(G)\},\; j \in \{1,\dots,M\}}.
  \end{equation*}
  We also set the initial weight vector as
  \begin{equation*}
    \beta_0 = \prescript{t}{}{\begin{pmatrix}w_1 & \dots & w_M
    \end{pmatrix}}.
  \end{equation*}
  %%%%%%%%%%%%%%%%
\item[(2)] Compute a basis of the kernel of the test function matrix: 
  Then, we calculate $\ker(A_0)$ using the singular value
  decomposition. 
  Since
  \begin{equation*}
 dim\big(\ker(A_0)\big)=M-\rank(A_0),  
  \end{equation*}
  we can choose a basis for
  $\ker(A_0)$ consisting of
  $M-\rank(A_0)$ elements
  denoted as 
  $\{v_{0,1}, v_{0,2}\dots, v_{0, M-\rank(A_0)}\}\subset\mathbb{R}^{M}$.
  %%% 
  %% 
\item[(3)] Eliminate nodes and update the weights: 
  We determine the set of nodes
 \begin{equation*}  
  (p_{e_i})_{i\in\{1,\dots,M-\rank(A_0)\}}
  \end{equation*}
  to be eliminated by repeating steps (3-1)–(3-3) for each $i$ as follows.
  \begin{itemize}
  \item[(3-1)] Determine the index $(e_i)$ of the node to be eliminated at step $i$: 
    we choose \(e_i\) that attains the minimum value of
    \[
      \frac{\beta_{i-1}^{(j)}}{(v_{i-1,1})^{(j)}}
    \]
    over all indices \(j\) satisfying \((v_{i-1,1})^{(j)} > 0\). 
    We define
    \[
      \alpha_i =
      \frac{\beta_{i-1}^{(e_i)}}
      {(v_{i-1,1})^{(e_i)}}.
    \]
    
  \item[(3-2)] Eliminate the node
    $(p_{e_i}, w_{e_i})$ and update the test function matrix and
    the weight vector: 
    By eliminating the $e_i$th column of $A_{i-1}$, we obtain $A_{i}$
    as an $\card(G) \times (M-i)$ matrix.
    Then, we update the probability measure by defining
    $\beta_{i} \in [0,1]^{M-i}$ as follows: 
    \begin{equation}\label{redbeta}
      \beta_{i}^{(j)} =
      \begin{cases}
	\beta_{i-1}^{(j)}-\alpha_i (v_{i-1,1})^{(j)}
	& \text{for } j=1,\dots,e_i-1,
	\\
	\beta_{i-1}^{(j+1)}-\alpha_i (v_{i-1,1})^{(j+1)}
	& \text{for } j=e_i,\dots,M-i.
      \end{cases}
    \end{equation}
  \item[(3-3)] Calculate the basis of kernel of the updated test function matrix: 
    For $l=1,\dots, M -\rank(A_0)-i$, we define 
%% \begin{equation}\label{redupdate}
%%   \begin{split}
%%     d_{i, l+1} &=\frac{(v_{i-1, l+1})^{(e_i)}}
%%     {(v_{i-1, 1})^{(e_i)}}, \\
%%     v_{i, l} &= v_{i-1, l+1}-d_{i, l+1}v_{i-1, 1}
%%                       \in \mathbb{R}^{M}, \\
%%     (v_{i, l})^{(j)} &=\begin{cases}
%%       (v_{i,l})^{(j)}
%%       & \text{for } j=1,\dots,e_i-1,\\
%%       (v_{i,l})^{(j+1)}
%%       & \text{for } j=e_i,\dots,M-i.
%%     \end{cases}
%%   \end{split}
%% \end{equation}

  \end{itemize}
  \begin{align}
    d_{i, l+1} &= \frac{(v_{i-1, l+1})^{(e_i)}}
                 {(v_{i-1, 1})^{(e_i)}}, \notag \\
    v_{i, l} &= v_{i-1, l+1} - d_{i, l+1} v_{i-1, 1}
               \in \mathbb{R}^{M}, \notag \\
    (v_{i, l})^{(j)} &= \begin{cases}
      (v_{i,l})^{(j)}
      & \text{for } j = 1, \dots, e_i - 1, \\
      (v_{i,l})^{(j+1)}
      & \text{for } j = e_i, \dots, M - i.
    \end{cases}
                       \label{redupdate}
  \end{align}
\end{enumerate}

\section{Derivation of an Illustrative example}
\label{appendix2}
Let $X, Y$ be $2$-dimensional discrete random variables on
$\Omega=\{\omega_0,\dots,\omega_7\}$ such that

\begin{align*}
(X,Y)(\omega_0) &= (0,0),
& (X,Y)(\omega_1) &= (1,0),
& (X,Y)(\omega_2) &= (0,1),
\\
(X,Y)(\omega_3) &= (1,1),
& (X,Y)(\omega_4) &= (2,0),
& (X,Y)(\omega_5) &= (2,1),
\\
(X,Y)(\omega_6) &= (2,2),
& (X,Y)(\omega_7) &= (1,2).
\end{align*}
Let $\tilde{Q}$ be a probability measure on $\Omega$ defined by
$\tilde{Q}(\omega_0)=\dots=\tilde{Q}(\omega_7)=1/8$.
The aim is to construct a reduced measure $Q$ such that the moments of
$(X,Y)$ under $Q$ agree with those under $\tilde{Q}$ up to the second order.
Equivalently, $Q$ is a reduced measure of $\tilde{Q}$ with respect to
$G=\big(g_i(X,Y)\big)_{i\in\{1,\dots,6\}}$,
where
$(g_i)_{i\in\{1,\dots,6\}}\subset\mathbb{R}[x,y]$ is defined as
%% \begin{gather*}
%%   g_1(x,y)=x^0y^0=1,\quad  g_2(x,y)=x,\quad  g_3(x,y)=y, \\
%%   g_4(x,y)=x^2,\quad  g_5(x,y)=y^2,\quad g_6(x,y)=xy.
%% \end{gather*}
\begin{align*}
  g_1(x,y) &= x^0y^0=1, &  g_2(x,y) &= x, &  g_3(x,y) &= y, \\
  g_4(x,y) &= x^2, &  g_5(x,y) &= y^2, &  g_6(x,y) &= xy.
\end{align*}
The problem reduces to finding a basis of $\Kernel{A}$
of the following matrix $A$, where the $i$th column and the $j$th row
correspond to the state $\omega_{i-1}$ and the monomial $g_j$ respectively:
\begin{equation*}
  A=\begin{pmatrix}
  1 & 1 & 1 & 1 & 1 & 1 & 1 & 1 \\
  0 & 1 & 0 & 1 & 2 & 2 & 2 & 1 \\
  0 & 0 & 1 & 1 & 0 & 1 & 2 & 2 \\
  0 & 1 & 0 & 1 & 4 & 4 & 4 & 1 \\
  0 & 0 & 1 & 1 & 0 & 1 & 4 & 4 \\
  0 & 0 & 0 & 1 & 0 & 2 & 4 & 2 
  \end{pmatrix}.
\end{equation*}
Then we see that $\Kernel(A) = \mathbb{R}v_1+\mathbb{R}v_2$
where
%% \begin{equation*}\begin{split}
%%     v_1&={\prescript{t}{}{\begin{pmatrix}
%%           -1 & 2 & 1 & -2 & -1 & 1 & 0 & 0 \end{pmatrix}}}, \\
%%     v_2&={\prescript{t}{}{\begin{pmatrix}
%%           2 & -3 & -2 & 2 & 1 & 0 & -1 & 1 \end{pmatrix}}}.
%% \end{split}\end{equation*}
\begin{align*}
  v_1
  &=
    {\prescript{t}{}{\begin{pmatrix}
      -1 & 2 & 1 & -2 & -1 & 1 & 0 & 0
    \end{pmatrix}}}, \\
  v_2
  &=
    {\prescript{t}{}{\begin{pmatrix}
      2 & -3 & -2 & 2 & 1 & 0 & -1 & 1
    \end{pmatrix}}}.
\end{align*}

If we set
%% \begin{equation*}
%%   \begin{split}
\begin{align*}
    Q&= \tilde{Q}-(1/8)v_1-(1/24)v_2\\
    &=\prescript{t}{}{\begin{pmatrix}
        {1}/{6} & {0} & {1}/{12} & {7}/{24}
        & {5}/{24} & {0} & {1}/{6} & {1}/{12}
    \end{pmatrix}},
\end{align*}
%%   \end{split}
%% \end{equation*}
we see that $Q$ satisfies
%% \begin{equation*}
%%   \begin{split}
%%     &E^{\tilde{Q}}[X]=E^{Q}[X], \;\; E^{\tilde{Q}}[Y]=E^{Q}[Y], \;\; \\
%%     &E^{\tilde{Q}}[X^2]=E^{Q}[X^2],\;\;E^{\tilde{Q}}[XY]=E^{Q}[XY], \;\;E^{\tilde{Q}}[Y^2]=E^{Q}[Y^2].
%%   \end{split}
%% \end{equation*}
\begin{align*}
  E^{\tilde{Q}}[X] &= E^{Q}[X], \qquad
  E^{\tilde{Q}}[Y] = E^{Q}[Y], \\
  E^{\tilde{Q}}[X^2] &= E^{Q}[X^2], \qquad
  E^{\tilde{Q}}[XY] = E^{Q}[XY], \qquad
  E^{\tilde{Q}}[Y^2] = E^{Q}[Y^2].
\end{align*}

\newpage
\section{Pseudocode of the high-order recombination algorithm}
\label{appendix3}
Let $\Delta = \{0=t_0<\dots<t_n=T\}$ be a partition of $[0,T]$.
As in Algorithm~\ref{alg:main}, the patch division of the nodes at time
$t_{i+1}$ uses the penalty function $\Psi_{i+1}$, and no recombination
is performed at the last step.

\SetKwComment{Comment}{/* }{ */}
\begin{algorithm}\normalsize
  \SetKwInOut{Parameter}{Parameter}
  \SetKwInOut{Input}{Input}
  \SetKwInOut{Output}{Output}
  \SetKwFunction{Main}{Main}
  \SetKwFunction{OneStepFwd}{OneStepFwd}
  \SetKwFunction{PatchDiv}{PatchDiv}
  \SetKwFunction{Recombine}{Recombine}
  
  \Input{Number of partitions : $n$,\\
    Parameters of SDE : $x_0$, $V_0,\dots, V_d$,\\
    Time horizon and payoff : maturity $T$, payoff function $f$}
  \Output{$E_{\mu_{t_n}^{({\rm Alg}, \Delta, {\rm Rec})}} \bigl[ f (X^{({\rm NV}, \Delta, {\rm Rec})}_{T}) \bigr]$}
  
  \MyStruct{Node}{
    $value$  \Comment*[r]{State $X^{({\rm NV}, \Delta, {\rm Rec})}_{t_i}(\omega_j)$}
    $weight$ \Comment*[r]{Probability weight $P[\omega_j]$}
    $patchnum$ \Comment*[r]{Patch index of this node}
  }
  
  \SetKwProg{Fn}{Function}{:}{}
  \Fn{\Main{}}{
    $leaves \gets \text{array of \textit{Node}}$\;
    $nextleaves \gets \text{array of \textit{Node}}$\;
    
    $leaves[0].value \gets x_0$\;
    $leaves[0].weight \gets 1.0$\;
    $leaves[0].patchnum \gets 0$\;
    
    \For{$i \gets 0$ \KwTo $n-1$}{
      $nextleaves \gets \OneStepFwd(leaves, i, n, T, V_0,\dots, V_d)$\;
      \If{$i+1<n$}{
        $numpatch \gets 0$\;
        \PatchDiv($nextleaves$, $i+1$, \&$numpatch$)\;
        $nextleaves \gets \Recombine(nextleaves, numpatch)$\;
      }
      $leaves \gets nextleaves$\;
    }
    $value \gets \displaystyle\sum_{\text{node} \in leaves} \text{node.weight}\cdot f(\text{node.value})$\;
    \Return{$value$}\;
  }
\end{algorithm}

\begin{algorithm}\normalsize
  \SetKwFunction{FOneStepFwd}{OneStepFwd}
  \SetKwProg{Fn}{Function}{:}{}
  
  \Fn{\FOneStepFwd{$leaves, i, n, T, V_0,\dots,V_d$}
    \Comment*[r]{Compute $X^{({\rm NV},\Delta,{\rm Rec})}_{t_{i+1}}$ from $t_i$}}{
    
    $nextleaves \gets \text{empty array of \textit{Node}}$\;
    $\Delta t \gets T/n$\;
    
    \ForAll{$node$ in $leaves$
      \Comment*[r]{Generate random variables and solve ODE flows}}{
      
      Generate random variables $(\eta_{i}^{(r)})_{r\in\{1,\dots,2\times 3^d\}}$\;
      \For{$r \gets 1$ \KwTo $2\times 3^d$}{
        $x_{\rm new} \gets \text{N--V-flow}\bigl(node.value, \Delta t, \eta_i^{(r)},
        V_0,\dots,V_d\bigr)$\;
        \Comment*[r]{Implementation of \eqref{eq:nvalg}}
        $newnode.value \gets x_{\rm new}$\;
        $newnode.weight \gets node.weight \times \frac{\phi(\eta_{i}^{(r)})}{2}$\;
        $newnode.patchnum \gets -1$\;
        append $newnode$ to $nextleaves$\;
      }
    }
    \Return{$nextleaves$}\;
  }
\end{algorithm}

\begin{algorithm}\normalsize
  \SetKwFunction{FPatchDiv}{PatchDiv}
  \SetKwFunction{ProjectionVector}{ProjectionVector}
  \SetKwProg{Fn}{Function}{:}{}
  
  \Fn{\FPatchDiv{$leaves, i, \&numpatch$}}{
    
    $M \gets \text{length}(leaves)$\;
    $\bar{\nu} \gets \sum_{j=1}^{M} w_j \delta_{p_j}$ from $leaves$\;
    \Comment*[r]{$p_j =$ node.value, $w_j =$ node.weight}
    
    \If{$\Psi_i(\bar{\nu}) \le 0$ \textbf{or} $M \le \card(G)$}{
      \ForAll{$node$ in $leaves$}{
        $node.patchnum \gets numpatch$\;
      }
      $numpatch \gets numpatch + 1$\;
      \Return\;
    }
    
    $v \gets \ProjectionVector(leaves)$\;
    $c \gets \displaystyle\sum_{j=1}^{M} w_j p_j$
    \Comment*[r]{Centre of mass \eqref{gravitypoint}}
    
    \ForAll{$node$ in $leaves$}{
      $node.order \gets \langle node.value - c, v\rangle$
      \Comment*[r]{Projection \eqref{eq:projinnerprod}}
    }
    \text{QuickSort}($leaves$ by $node.order$)\;
    
    $M_{\rm left} \gets \lfloor M/2 \rfloor$\;
    $leftleaves \gets leaves[1:M_{\rm left}]$\;
    $rightleaves \gets leaves[M_{\rm left}+1:M]$\;
    
    \FPatchDiv{$leftleaves, i, \&numpatch$}\;
    \FPatchDiv{$rightleaves, i, \&numpatch$}\;
  }
\end{algorithm}

\begin{algorithm}\normalsize
  \SetKwFunction{FRecombine}{Recombine}
  \SetKwProg{Fn}{Function}{:}{}
  
  \Fn{\FRecombine{$leaves, numpatch$}}{
    
    $newleaves \gets \text{empty array of \textit{Node}}$\;
    
    \For{$k \gets 0$ \KwTo $numpatch-1$}{
      $U_k \gets \{\, node \in leaves \,:\, node.patchnum = k \,\}$\;
      $M_k \gets \text{length}(U_k)$\;
      
      \If{$M_k \le \card(G)$}{
        append all nodes in $U_k$ to $newleaves$\;
        \textbf{continue}\;
      }
      
      \tcp{(1) Build test-function matrix and weight vector}
      Construct $A_0 \in \mathbb{R}^{\card(G) \times M_k}$ with
      $(A_0)_{ij} = g_i(p_j)$ and
      %%      $\beta_0 = (w_1,\dots,w_{M_k})^\mathsf{T}$\;
      $\beta_0 = \prescript{t}{}{\begin{pmatrix}
          w_1 & \dots & w_{M_k}\end{pmatrix}}$\;
      Compute $\rank (A_0)$ and a basis
      $\{v_{0,1},\dots, v_{0,M_k-\rank(A_0)}\} \subset \mathbb{R}^{M_k}$
      of $\ker(A_0)$ using SVD\;
      
      $A\gets A_0,\ \beta \gets \beta_0,\ r_{\max}\gets M_k-\rank(A_0)$
      \;
      \tcp{(2) Sequential elimination on patch $U_k$}
      \For{$r \gets 1$ \KwTo $r_{\max}$}{
        Let $v_1$ be the first basis vector of the current kernel\;
        Choose index $e_r$ that minimises
        $\beta^{(j)}/v_{1}^{(j)}$ over all $j$ such that $v_1^{(j)} > 0$\;
        $\alpha_r \gets \beta^{(e_r)}/v_1^{(e_r)}$\;
        
        \tcp{Update weights as in \eqref{redbeta}}
        Update $\beta$ by
        \[
        \beta^{(j)} \gets
        \begin{cases}
          \beta^{(j)} - \alpha_r v_1^{(j)}, & j<e_r,\\
          \beta^{(j+1)} - \alpha_r v_1^{(j+1)}, & j\ge e_r,
        \end{cases}
        \]
        and remove column $e_r$ from $A$\;
        
        \tcp{Update kernel basis as in \eqref{redupdate}}
        Update the kernel basis vectors $\{v_{l}\}$ accordingly\;
      }
      
      \tcp{(3) Construct reduced measure on the patch}
      Let the remaining columns of $A$ correspond to nodes
      $\{p_{k,1},\dots,p_{k,\rank(A_0)}\}$ with weights
      $\tilde{w}_{k,1},\dots,\tilde{w}_{k,\rank(A_0)}$ from $\beta$\;
      For each remaining node, create a \textit{Node} with
      $value = p_{k,j}$, $weight = \tilde{w}_{k,j}$,
      $patchnum = k$ and append to $newleaves$\;
    }
    
    \Return{$newleaves$}\;
  }
\end{algorithm}

\FloatBarrier  

\section{Numerical experiment setup (Summary)}\label{appendix4}
This appendix summarises the parameter choices used
in the numerical experiments.
Tables~\ref{tab:exp_setting}--\ref{tab:spiky_payoff_setting} list the setups.
\begin{table}[!t]\normalsize
  \centering
  \caption{Discretisation Errors in Section~\ref{validation1} (Heston model, Asian call option)}
  \label{tab:exp_setting}
  \begin{tabular}{l|l}
    \hline
    Model & Heston model \\
    Payoff & $\max(Y^{(3)}_{T}/T - K, 0)$ \\
    \hline
    Maturity & $T = 1$ \\
    Strike & $K = 1.05$ \\
    Initial price & $x_1 = 1.0$ \\
    Initial variance & $x_2 = 0.09$ \\
    Drift & $\mu = 0.05$ \\
    Mean reversion & $\alpha = 2.0$ \\
    Long-term variance & $\theta = 0.09$ \\
    Vol-of-vol & $\beta = 0.1$ \\
    Correlation & $\rho = 0.3$ \\
    \hline
    Time partition & Even partition: $s_i = T/\sharp\Delta$ \\
    Discretisation schemes & Euler--Maruyama (E--M), N--V  \\
    % Revised in November 2026: the algorithm of subsection 3.2 with kappa=1.
    %Patch condition & \eqref{eq:refcond} with $\ell = 3$ \\
    Patch condition & \rev{$\Psi_i$ of subsection~\ref{patchdivsec} with $\kappa=1$ and} $\ell = 3$ \\
    Reference value & $0.06068740243939$ (95\% CI: $[0.0606855, 0.0606896]$)\\%$[0.0606824, 0.0606874]$) \\
    \hline
  \end{tabular}
\end{table}

\begin{table}[!t]\normalsize
  \centering
  \caption{Recombination Errors in Section~\ref{subsec:RecErr} (Heston model, Asian call option)}
  \label{tab:recerr_asian}
  \begin{tabular}{l|l}
    \hline
    Underlying model & Heston model (same parameters as in Table~\ref{tab:exp_setting}) \\
    Payoff & Asian call: $\max(Y^{(3)}_{T}/T - K, 0)$ \\
    Discretisation & N--V  scheme with recombination \\
    Time partition & Even partition: $\Delta = \{0 = t_0 < \dots < t_{\sharp\Delta} = T\}$, $\sharp\Delta =  8$ \\
    Reference measure & $\mu_{T}^{({\rm NV}, \Delta)}$ (full tree, no recombination) \\
    Recombined measure & $\mu_{T}^{({\rm NV}, \Delta, {\rm Rec}, \lambda)}$ \\
    %Recombined measure & $\mu_T^{({\rm NV}, \Delta, {\rm Rec}, \text{PatchDiv})}$ \\
    Error metric & Recombination error in \eqref{eq:recomb-error} \\
    Patch division methods &
    Recursive\_WLL, Recursive\_LL, CC, Rand \\
    Patch condition & \eqref{eq:scaledpel} and Example~\ref{ori:cond} with $\ell = 3$ \\
    \hline
  \end{tabular}
\end{table}

\begin{table}[!t]\normalsize
  \centering
  \caption{Recombination Errors in Section~\ref{subsec:RecErr} (Heston model, Cone payoffs)}
  \label{tab:spiky_payoff_setting}
  \begin{tabular}{l|l}
    \hline
    Underlying variables &
    $\tilde{X}_{T} = Y^{(3)}_{T}/Y^{(1)}_{0}$,
    $\tilde{Y}_{T} = Y^{(2)}_{T}/Y^{(2)}_{0}$ \\
    Payoff family &
    $g^{(x_c, y_c, r, h)}(x, y )=\max\{0,h-\frac{h}{r}\sqrt{(x-x_c)^2+(y-y_c)^2}\}$\\
    Cone height &
    $h = 10 r^{-2}$ (cone volume kept constant) \\
    Radii $r$ &
    $5.0,\, 4.0,\, 3.0,\, 2.0,\, 1.5,\, 1.2,\, 1.0,\, 0.9$ \\
    Cone centres $(x_c,y_c)$ &
    $(0,0), (1,0), (0,1), (1,1), (2,0), (0,2),$ \\
    & $(2,2), (3,0), (0,3), (3,3)$ \\
    Model and parameters &
    Same Heston setup as in Table~\ref{tab:exp_setting} \\
    Discretisation &
    N--V scheme with recombination \\
    Time partitions &
    $\sharp\Delta = 6$ and $\sharp\Delta = 8$ (even partitions) \\
    Patch division methods &
    Recursive\_WLL, Recursive\_LL, CC, Rand \\
    Patch condition & \eqref{eq:scaledpel} and Example~\ref{ori:cond} with $\ell = 3$ \\
    %% Patch split rule &
    %% Equal-number node division in recursive algorithms \\
    Error metric&
    Mean slopes of the log-log regressions\\
    \hline
  \end{tabular}
\end{table}
\FloatBarrier  
%

%% bibliography 
%%% forced to use bibtex %%% \printbibliography
\bibliographystyle{amsplain}
\bibliography{MyCollection} % .bib

\providecommand{\bysame}{\leavevmode\hbox to3em{\hrulefill}\thinspace}
\providecommand{\MR}{\relax\ifhmode\unskip\space\fi MR }
\providecommand{\MRhref}[2]{%
  \href{http://www.ams.org/mathscinet-getitem?mr=#1}{#2}
}
\providecommand{\href}[2]{#2}
\begin{thebibliography}{10}

\bibitem{bally1996law}
Vlad Bally and Denis Talay, \emph{The law of the euler scheme for stochastic
  differential equations}, Probability theory and related fields \textbf{104}
  (1996), no.~1, 43--60.

\bibitem{rQMCori}
R.~Cranley and T.~N.~L. Patterson, \emph{Randomization of number theoretic
  methods for multiple integration}, SIAM Journal on Numerical Analysis
  \textbf{13} (1976), no.~6, 904--914.

\bibitem{dan2002minimal}
Dan Crisan and Terry Lyons, \emph{{Minimal entropy approximations and optimal
  algorithms}}, Monte Carlo Methods and Applications \textbf{8} (2002), no.~4,
  343--356.

\bibitem{Crisan2012}
Dan Crisan and Konstantinos Manolarakis, \emph{{Solving backward stochastic
  differential equations using the cubature method: Application to nonlinear
  pricing}}, SIAM Journal on Financial Mathematics \textbf{3} (2012), no.~1,
  534--571.

\bibitem{Crisan2014}
\bysame, \emph{{Second order discretization of backward SDEs and simulation
  with the cubature method}}, Annals of Applied Probability \textbf{24} (2014),
  no.~2, 652--678.

\bibitem{feller1951two}
William Feller, \emph{{Two Singular Diffusion Problems}}, Annals of Mathematics
  \textbf{54} (1951), no.~1, 173--182.

\bibitem{ferrucci2024high}
Emilio Ferrucci, Timothy Herschell, Christian Litterer, and Terry Lyons,
  \emph{High-degree cubature on wiener space through unshuffle expansions},
  arXiv preprint arXiv:2411.13707 (2024).

\bibitem{fujiwara}
Takehiko Fujiwara, \emph{{Sixth order methods of Kusuoka approximation}},
  PreprintSeries, Graduate School of Mathematical Sciences, Univ. Tokyo (2006).

\bibitem{gyurko2011efficient}
Lajos Gyurk{\'{o}} and Terry Lyons, \emph{{Efficient and practical
  implementations of cubature on Wiener space}}, Stochastic analysis 2010,
  Springer, 2011, pp.~73--111.

\bibitem{hayakawa2022monte}
Satoshi Hayakawa and Ken’ichiro Tanaka, \emph{Monte carlo construction of
  cubature on wiener space}, Japan Journal of Industrial and Applied
  Mathematics \textbf{39} (2022), no.~2, 543--571.

\bibitem{heston1993closed}
Steven Heston, \emph{{A closed-form solution for options with stochastic
  volatility with applications to bond and currency options}}, Review of
  financial studies \textbf{6} (1993), no.~2, 327--343.

\bibitem{kusuoka2001approximation}
Shigeo Kusuoka, \emph{{Approximation of expectation of diffusion process and
  mathematical finance}}, Proceedings of Final Taniguchi Symposium, Nara,
  Advanced Studies in Pure Mathematics, vol.~31, 2001, pp.~147--165.

\bibitem{kusuoka2003malliavin}
\bysame, \emph{{Malliavin calculus revisited}}, J. Math. Sci. Univ. Tokyo
  \textbf{10} (2003), no.~2, 261--277.

\bibitem{kusuoka2004approximation}
\bysame, \emph{{Approximation of expectation of diffusion processes based on
  Lie algebra and Malliavin calculus}}, Advances in Mathematical Economics
  \textbf{6} (2004), 69--83.

\bibitem{kusuoka2013gaussian}
\bysame, \emph{{Gaussian K-scheme: justification for KLNV method}}, Advances in
  Mathematical Economics \textbf{17} (2013), 71--120.

\bibitem{lee2016adaptive}
Wonjung Lee and Terry Lyons, \emph{The adaptive patched cubature filter and its
  implementation}, Communications in Mathematical Sciences (2016).

\bibitem{LittererDthesis}
Christian Litterer, \emph{The signature in numerical algorithms}, Ph.D. thesis,
  University of Oxford, 2008.

\bibitem{litterer2012high}
Christian Litterer and Terry Lyons, \emph{{High order recombination and an
  application to cubature on Wiener space}}, The Annals of Applied Probability
  \textbf{22} (2012), no.~4, 1301--1327.

\bibitem{lyons2004cubature}
Terry Lyons and Nicolas Victoir, \emph{{Cubature on Wiener Space}}, Proceedings
  of the Royal Society of London A: Mathematical, Physical and Engineering
  Sciences, vol. 460, 2004, pp.~169--198.

\bibitem{rQMCsurvey}
Pierre L’Ecuyer and Christiane Lemieux, \emph{Recent advances in randomized
  quasi-monte carlo methods}, Modeling uncertainty: An examination of
  stochastic theory, methods, and applications (2002), 419--474.

\bibitem{ninomiya2010application}
Mariko Ninomiya, \emph{{Application of the {K}usuoka approximation with a
  tree-based branching algorithm to the pricing of interest-rate derivatives
  under the {HJM} model}}, LMS Journal of Computation and Mathematics
  \textbf{13} (2010), 208--221.

\bibitem{ninomiya2009new}
Mariko Ninomiya and Syoiti Ninomiya, \emph{{A new higher-order weak
  approximation scheme for stochastic differential equations and the
  Runge--Kutta method}}, Finance and Stochastics \textbf{13} (2009), no.~3,
  415--443.

\bibitem{ninomiya2003new}
Syoiti Ninomiya, \emph{{A new simulation scheme of diffusion processes:
  application of the Kusuoka approximation to finance problems}}, Mathematics
  and Computers in Simulation \textbf{62} (2003), no.~3, 479--486.

\bibitem{ninomiya2003partial}
\bysame, \emph{{A partial sampling method applied to the Kusuoka
  approximation}}, Monte Carlo methods and Applications \textbf{9} (2003),
  no.~1, 27--38.

\bibitem{NinomiyaShinozaki2019}
Syoiti Ninomiya and Yuji Shinozaki, \emph{{Higher-order Discretization Methods
  of Forward-backward SDEs Using KLNV-scheme and Their Applications to XVA
  Pricing}}, Applied Mathematical Finance \textbf{26} (2019), no.~3, 257--292.

\bibitem{rec2021}
\bysame, \emph{{O}n implementation of high-order recombination and its
  application to weak approximations of stochastic differential equations},
  Proceedings of the NFA 29th Annual Conference, 2021.

\bibitem{ninomiya2008weak}
Syoiti Ninomiya and Nicolas Victoir, \emph{{Weak approximation of stochastic
  differential equations and application to derivative pricing}}, Applied
  Mathematical Finance \textbf{15} (2008), no.~2, 107--121.

\bibitem{oshima2012new}
Kojiro Oshima, Josef Teichmann, and Dejan Velu{\v{s}}{\v{c}}ek, \emph{{A new
  extrapolation method for weak approximation schemes with applications}},
  Annals of Applied Probability \textbf{22} (2012), no.~3, 1008--1045.

\bibitem{Shinozaki2017}
Yuji Shinozaki, \emph{{Construction of a third-order K-scheme and its
  application to financial models}}, SIAM Journal on Financial Mathematics
  \textbf{8} (2017), no.~1, 901--932.

\bibitem{tchernychova2016caratheodory}
Maria Tchernychova, \emph{Carath{\'e}odory cubature measures}, Ph.D. thesis,
  University of Oxford, 2016.

\bibitem{uchino2020}
Katsuya Uchino, \emph{Application of recombination to higher-order weak
  approximation problems of stochastic differential equation (in japanese)},
  Master thesis, Tokyo Institute of Technology (2020).

\end{thebibliography}
\end{document}